\documentclass[12pt,amsfonts]{article}

\usepackage{latexsym}
\usepackage{amsfonts,amssymb} 

\newtheorem{theorem}{Theorem}[section]
\newtheorem{lemma}[theorem]{Lemma}
\newtheorem{corollary}[theorem]{Corollary}
\newtheorem{proposition}[theorem]{Proposition}

\newtheorem{definition}[theorem]{Definition}
\newcommand{\bd}[1]{\begin{definition}\label{#1}\rm}
\newcommand{\ed}{\end{definition}}
\newcommand{\bt}[1]{\begin{theorem}\label{#1}}
\newcommand{\et}{\end{theorem}}
\newcommand{\bprop}[1]{\begin{proposition}\label{#1}}
\newcommand{\eprop}{\end{proposition}}
\newcommand{\bcor}[1]{\begin{corollary}\label{#1}}
\newcommand{\ecor}{\end{corollary}}

\newcommand{\D}{\displaystyle}
\newcommand{\T}{\textstyle}
\newcommand{\lra}{\longrightarrow}
\newcommand{\Ra}{\Longrightarrow}

\newcommand{\stack}[2]{\raisebox{-2pt} 
{\renewcommand{\arraystretch}{.01} 
\begin{tabular}{c} 
$#2$\\$\scriptscriptstyle #1$ 
\end{tabular} 
}} 

\newcommand{\wstack}[2]{\raisebox{-2pt} 
{\renewcommand{\arraystretch}{.01} 
\begin{tabular}{c} 
$\scriptscriptstyle w$\\$#2$\\$\scriptscriptstyle #1$ 
\end{tabular} 
}} 

\newcommand{\sstack}[2]{\raisebox{-2pt} 
{\renewcommand{\arraystretch}{.01} 
\begin{tabular}{c} 
$\scriptscriptstyle s$\\$#2$\\$\scriptscriptstyle #1$ 
\end{tabular} 
}}

\newcommand{\vp}{\varphi}
\newcommand{\ve}{\varepsilon}
\newcommand{\nid}{\noindent}
\newcommand{\qed}{\hfill$\Box$} 
\def\nn{|||} 
\def\0{\, {\rm 0}\mskip-11mu 0} 
\def\1{\, {\rm I}\mskip-10mu 1} 
\, 

\def\Del{\, \Delta{\hspace{-3mm}\Delta}}

\renewcommand{\t}[1]{\tilde{#1}} 
\newcommand{\bomega}{\mbox{\boldmath${\omega}$}}
\newcommand{\bmu}{\mbox{\boldmath${\mu}$}}
\newcommand{\bnu}{\mbox{\boldmath${\nu}$}} 
 
\newcommand{\sbnu}{\mbox{\scriptsize\boldmath${\nu}$}} 

\begin{document}
\title{Mosco Type Convergence and Weak Convergence for a Fleming-Viot 
type Particle System} 
\par 
\author{J\"org-Uwe L\"obus 
\\ Matematiska institutionen \\ 
Link\"opings universitet \\ 
SE-581 83 Link\"oping \\ 
Sverige 
}
\date{}
\maketitle
{\footnotesize
\noindent
\begin{quote}
{\bf Abstract}
We are concerned with Mosco type convergence for a non-symmetric 
$n$-particle Fleming-Viot system $\{X_1,\ldots ,X_n\}$ in a bounded 
$d$-dimensional domain $D$ with smooth boundary. Moreover, we are interested 
in relative compactness of the $n$-particle processes. It turns out that 
integration by parts relative to the initial measure and the generator 
is the appropriate mathematical tool. For finitely many particles, such 
integration by parts is established by using probabilistic arguments. For 
the limiting infinite dimensional configuration we use a result from 
infinite dimensional non-gaussian calculus. 
\noindent 


\end{quote}
}

\section{Introduction}
\setcounter{equation}{0} 

We consider a system $\{X_1,\ldots ,X_n\}$ of $n$ particles in a bounded 
$d$-dimensional domain $D$ ($d\ge 2$) with smooth boundary $\partial D$. 
During periods none of the particles $X_1,\ldots ,X_n$ hit the boundary 
$\partial D$, the system behaves like $n$ independent $d$-dimensional 
Brownian motions. When one of the particles hits the boundary $\partial 
D$, then it instantaneously jumps back to $D$ and relocates according to 
a probability distribution $\eta\, dx$ on $(D,{\cal B}(D))$. The 
probability distribution $\eta\, dx$ depends on the location of the 
remaining $n-1$ particles in a way that relocation is more likely close 
to one of those particles than elsewhere. Similar models have been 
investigated in \cite{AFG11}, \cite{BP09}, \cite{FM07}, and \cite{Vi11}. 
For the background in the sciences we refer to \cite{CBH05}. 

The situation that a particle which has just jumped relocates at the site 
of one of the remaining $n-1$ particles with probability $1/(n-1)$ has 
been studied in \cite{BHM00}, \cite{GK04}, and \cite{Lo09}.  
\medskip 

For the measure valued processes $X^n_t=\frac1n\sum_{i=1}^n\delta_{(X_i 
)_t}$, $t\ge 0$, we are establishing Mosco type convergence and relative 
compactness in the Skorokhod space $D_E([0,\infty))$. Here $E$ is the 
compact space ${\cal M}_\partial(\overline{D})$ of all equivalence 
classes $\mu$ of probability measures on $\left(\overline{D},{\cal B} 
(\overline{D})\right)$ such that $m_1,m_2\in\mu$ implies $m_1|_D=m_2|_D$. 
For the precise definition of the metric in ${\cal M}_\partial(\overline 
{D})$, see Subsection 2.2. The reason for choosing $E={\cal M}_\partial 
(\overline{D})$ is the crucial technical condition (c7), cf. Subsection 
2.2 and Theorem \ref{Theorem4.13}. Weak convergence of ${\bf X}^n=\left( 
(X^n_t)_{t\ge 0},P_{{\sbnu}_n}\right)$ for a certain class of initial 
measures $\bnu_n$ to ${\bf X}=\left((X_t)_{t\ge 0},P_{\sbnu}\right)$ as 
$n\to\infty$ is now obtained by the the weak convergence of the finite 
dimensional distributions. We get this precisely as in \cite{Ko08}, proof 
of Corollary 2.9, together with our final Remark (1) of Section 7. Similar 
methods are applied for example in \cite{CGLW12}, \cite{GKR07}, \cite{Ki06}, 
\cite{{KU97}, {KU96}}. 
\medskip 

Let ${\cal M}_1(D)$ denote the set of all probability measures on $(D,{\cal 
B}(D))$. The trajectory $v(t,x)\, dx$, $t\ge 0$, with initial value $\mu\in 
{\cal M}_1(D)$ of the limiting process $X_t$, $t\ge 0$, is the solution to 
the non-linear PDE 
\begin{eqnarray*} 
\frac{\partial}{\partial t}v(t,x)=\frac12\Delta v(t,x)-\frac{z'(t)}{z(t)} 
v(t,x)\, ,\quad v(t,x)|_{x\in \partial D}=0,\ v(t,x)\, dx\stack {t\to 0} 
{\Ra}\mu\, , 
\end{eqnarray*}
established in \cite{GK04}. Here $z(t)=\int_{y\in D}\int_{x\in D}p(t,x,y) 
\, \mu(dx)\, dy$ where $p(t,x,y)$, $x,y\in D$, $t\ge 0$, is the Lebesgue 
transition density of the $d$-dimensional Brownian motion killed at 
$\partial D$. This solution can be represented in the form $v(t,y)=\frac 
{1}{z(t)}\int_{x\in D}p(t,x,y)\, \mu(dx)$, $t\ge 0$, $y\in D$. It becomes 
immediate that $v(t,x)\, dx$ stays in ${\cal M}_1(D)$ for $t\ge 0$ whenever 
$\mu=v(0,x)\, dx\in {\cal M}_1(D)$. 

Moreover by definition, the trajectories of the processes ${\bf X}^n=\left 
((X^n_t)_{t\ge 0},P_{{\sbnu}_n}\right)$, $n\in {\Bbb N}$, take values in 
${\cal M}_1(D)$ for all $t\ge 0$, almost surely. Noting that $\bnu_n$, $n 
\in {\Bbb N}$, and $\bnu$ are concentrated on ${\cal M}_1(D)$, weak 
convergence in the Skorokhod space $D_{{\cal M}_\partial (D)}([0,\infty))$ 
implies weak convergence in $D_{{\cal M}_1(D)}([0,\infty))$. 
\bigskip 

We recall that Mosco type convergence together with relative compactness 
in the Skorokhod space is in general stronger than weak convergence in the 
(same) Skorokhod space. A possible consequence are relatively strong 
conditions on the initial measures, cf. for example \cite{{KU97}, {KU96}} 
and \cite{Ko08}. 

In the present paper, the limiting initial measures $\bnu$ are concentrated 
on a certain class of perturbations of $h_1/|h_1\|_{L^1(D)}$, see Section 4. 
Here $h_1$ is the positive version of the first eigenfunction of the 
Dirichlet Laplacian in $L^2(D)$. The class of initial measures $\bnu_n$, $n 
\in {\Bbb N}$, appears rather specific. However, in order to guarantee 
$\bnu_n\Ra\bnu$ (weekly), this class is quite natural. 

The paper \cite{GK04} shows that in a comparable particle system weak 
convergence in an appropriate Skorokhod space does not require particular 
conditions on the initial measures. Among other things, the present paper 
is meaningful from the point of view of strengthening the mode of convergence. 

\medskip
Although the paper primarily addresses to Mosco-type convergence, relative 
compactness, and weak convergence of the processes ${\bf X}^n=\left((X^n_t 
)_{t\ge 0},P_{{\sbnu}_n}\right)$, $n\in {\Bbb N}$, we also prove 
characterizing properties of the system. In particular, these are a limit 
theorem as $n\to\infty$ for the jump-off location of the $n$-particle 
process $(X_1,\ldots ,X_n)$ in Proposition \ref{Proposition4.5} (a) and a 
limit theorem as $n\to\infty$ for the jump distances of the $n$-particle 
process $(X_1,\ldots ,X_n)$ in Proposition \ref{Proposition4.5} (b). 
\medskip 

Among the many technical results in the paper we would like to point to 
finite dimensional integration by parts with respect to the processes 
$X^n_t$, $t\ge 0$, in Proposition \ref{Proposition4.10} (a). This partial 
integration result is used to prove the convergences in Proposition 
\ref{Proposition4.10} (b) and (c) which are powerful tools for the 
Mosco-type convergence. 
\bigskip

As in \cite{Lo14-1} we conclude the introduction with a remark on the 
notation in the paper. The greek letter $\bnu$ comes always in bold. This 
letter is exclusively used to denote probability measures over spaces of 
probability measures, the states of measure valued stochastic processes. 
Those states of stochastic processes are denoted using the greek letter 
$\mu$ non-bold. 

However, if $\mu$ has a subindex, as for example $\mu_x$ or $\mu_\cdot$, 
then it denotes some harmonic measure. In order to better distinguish 
between certain classes of harmonic measures we also use superindices. 

All expressions (functions, measures) that come with a tilde $\ \t{}\ $, 
$\ \widetilde{}\ $ are related to $n\cdot d$-dimensional $n$-particle 
processes $\left((X_1)_t,\ldots ,(X_n)_t\right)$ symmetric in the $n$ 
entries rather than to measure valued processes $X^n_t=\frac1n\sum_{i= 
1}^n\delta_{(X_i)_t}$, $t\ge 0$. 

\section{Preliminaries} 
\setcounter{equation}{0} 

In this section we introduce some basic notation and summarize the facts 
we take from \cite{Lo14-1} and \cite{Lo14-2}. 

\subsection{Mosco Type Convergence}

{\bf Convergences on sequences of $L^2$-spaces.} 
Let $\bnu_n$, $n\in {\Bbb Z}_+$, be mutually orthogonal probability 
measures on some measurable space $(E,{\cal B})$. Suppose that $\bnu 
\equiv\bnu_0$ is a measure with countable base. In addition, assume 
that there are mutually exclusive subsets $E_n$, $n\in {\Bbb Z}_+$, 
of $E$ such that ${\bnu}_n(E\setminus E_n)=0$. Let $\alpha_n$, $n\in 
{\Bbb Z}_+$, be a sequence of positive numbers with $\sum_{n=0}^\infty 
\alpha_n=1$. 
Define ${\Bbb M}:=\sum_{n=0}^\infty\alpha_n{\bnu}_n$. We say that $u 
\in\bigcap_{n\in{\Bbb Z}_+}L^2(E,{\bnu}_n)$ if $u$ is an equivalence 
class consisting of all everywhere defined ${\cal B}$-measurable 
functions satisfying $f_1=f_2$ ${\Bbb M}$-a.e. if $f_1,f_2\in u$ and 
$\int u^2\, d{\bnu}_n <\infty$, $n\in {\Bbb Z}_+$. Let $\langle\cdot 
\, ,\, \cdot\rangle_n$ denote the inner product in $L^2(E,{\bnu}_n)$, 
$n\in {\Bbb N}$, and let $\langle\cdot \, ,\, \cdot \rangle$ denote 
the inner product in $L^2(E,{\bnu})$. 
Introduce 
\begin{eqnarray*}
{\cal D}:=\left\{\vp\in \bigcap_{\, n\in {\Bbb Z}_+}L^2(E,{\bnu}_n): 
\langle \vp\, , \, \vp\rangle_n\stack{n\to\infty}{\lra}\langle \vp 
\, , \, \vp\rangle \right\}\, . 
\end{eqnarray*}
Suppose that there exists a linear subset ${\cal F}$ of ${\cal D}$ 
which is dense in $L^2(E,\bnu)$ and let ${\cal C}$ denote the set of 
all functions $\vp\in {\cal D}$ satisfying the following conditions: 
\begin{itemize} 
\item[(c1)] For each $\vp\in {\cal C}$, there exists a {\it 
representing sequence} $\vp_n\in {\cal F}$, $n\in {\Bbb N}$, 
such that $\vp =\vp_n$, $\bnu_n$-a.e., $n\in {\Bbb N}$. 
\item[(c2)] $\langle \vp\, , \, \psi\rangle_n\stack{n\to\infty}{\lra} 
\langle\vp\, , \, \psi\rangle$ for all $\psi\in {\cal F}$.  
\end{itemize}
\begin{lemma}\label{Lemma2.1} 
(a) ${\cal F}\subseteq {\cal C}$. \\ 
(b) The set ${\cal C}$ is linear. \\ 
(c) The set ${\cal C}$ is dense in $L^2(E,\bnu)$. \\ 
(d) Let $\vp,\psi\in {\cal C}$. We have $\langle \vp\, , \, \psi\rangle_n 
\stack{n\to\infty}{\lra}\langle\vp\, , \, \psi\rangle$. \\ 
\end{lemma}
\begin{definition}\label{Definition2.2} {\rm 
(a) A sequence $\vp_n\in {\cal C}$, $n\in {\Bbb N}$, is said to 
be {\it w-convergent} to $\vp\in L^2(E,{\bnu})$ as $n\to\infty$ 
(in symbols $\vp_n\wstack{n\to\infty}{\lra} \vp$) if 
\begin{itemize} 
\item[(i)] $\langle\vp_n\, , \, \psi\rangle_n\stack{n\to 
\infty}{\lra}\langle\vp\, , \, \psi\rangle$ for all $\psi\in 
{\cal C}$.  
\end{itemize} 
(b) A sequence $\psi_n\in {\cal C}$, $n\in {\Bbb N}$, is said to 
be {\it s-convergent} to $\psi\in L^2(E,{\bnu})$ as $n\to 
\infty$ (in symbols $\psi_n\sstack{n\to\infty}{\lra} \psi$) if 
\begin{itemize} 
\item[(i)] $\psi_n$ $w$-converges to $\psi$ as $n\to\infty$ and 
\item[(ii)] $\langle \psi_n\, , \, \psi_n\rangle_n \stack{n\to 
\infty}{\lra}\langle \psi\, , \, \psi\rangle$.
\end{itemize} 
(c) Speaking of $w$-convergence or $s$-convergence of subsequences 
$\vp_{n_k}\in {\cal C}$ or $\psi_{n_k}\in {\cal C}$, respectively, 
will mean that in (a) or (b) the index $n\in {\Bbb N}$ is replaced 
with $n_k\in {\Bbb N}$. }
\end{definition}
\begin{proposition}\label{Proposition2.3} 
(a) Let $\vp_n\in {\cal C}$, $n\in {\Bbb N}$, be a sequence 
$w$-convergent to $\vp\in L^2(E,{\bnu})$ as $n\to\infty$. Then 
$\langle \vp_n\, , \, \vp_n\rangle_n$, $n\in {\Bbb N}$, is bounded. 
 \\ 
(b) Let $\vp_n\in {\cal C}$, $n\in {\Bbb N}$, be a sequence 
such that $\langle \vp_n\, , \, \vp_n\rangle_n$ is bounded. Then 
there exists a subsequence $\vp_{n_k}\in {\cal C}$, $k\in {\Bbb N}$, 
$w$-convergent to some $\vp\in L^2(E,{\bnu})$ as $k\to\infty$. \\  
(c) Let $\vp_n\in {\cal C}$, $n\in {\Bbb N}$, be a sequence 
$w$-convergent to $\vp\in L^2(E,{\bnu})$ and let $\psi_n\in{\cal C}$, 
$n\in {\Bbb N}$, be a sequence that $s$-converges to $\psi\in L^2(E, 
{\bnu})$ as $n\to\infty$. Then $\langle\vp_n\, ,\, \psi_n\rangle_n 
\stack{n\to\infty}{\lra}\langle\vp\, ,\, \psi\rangle$. 
\end{proposition}
\medskip 

\nid
{\bf Bilinear forms.} Let $(T_t)_{t\ge 0}$ be a strongly continuous 
semigroup of linear operators on $L^2(E,{\bnu})$. Suppose that $(T_t 
)_{t\ge 0}$ is associated with a transition probability function 
$P(t,x,B)$, $t\ge 0$, $x\in E$, $B\in {\cal B}$, i. e., $T_tf=\int 
f(y)\, P(t,\cdot ,dy)$, $t\ge 0$, $f\in L^2(E,{\bnu})$. Assume, 
furthermore, that $P(t,\cdot ,E)=1$ ${\bnu}$-a.e., $t\ge 0$. 
\medskip

Denoting by $(A,D(A))$ the generator of $(T_t)_{t\ge 0}$ and by 
$\langle \cdot \, , \, \cdot \rangle$ the inner product in $L^2(E,{ 
\bnu})$, we introduce now the class of bilinear forms $S$ we are 
interested in. Define 
\begin{eqnarray*}
\ D(S):=\left\{u\in L^2(E,{\bnu}):\, \lim_{t\to 0}\left\langle 
\textstyle{\frac1t}(u-T_tu)\, , \, v\right\rangle \ \mbox{\rm exists  
for all}\ v\in L^2(E,{\bnu})\right\} 
\end{eqnarray*}
and 
\begin{eqnarray*}
S(u,v):=\lim_{t\to 0}\left\langle\textstyle{ \frac1t}(u-T_tu)\, , \, 
v\right\rangle \, , \quad u\in D(S), \ v \in L^2(E,{\bnu}). 
\end{eqnarray*}
We have $D(A)=D(S)$ according to \cite{P83}, Section 2.1 and 
\begin{eqnarray*}
S(u,v)=-\langle Au \, , \, v \rangle\, , \quad u\in D(A), \ v\in 
L^2(E,{\bnu}).
\end{eqnarray*}
In this sense we would like to understand the term {\it bilinear form}. 
However, as it is customary for Mosco (type) convergence, we also set 
$S(u,v):=\infty$ if $u\in L^2(E,{\bnu})\setminus D(S)$ and $v\in L^2(E, 
{\bnu})$. Let $(G_\beta)_{\beta >0}$ be the resolvent associated with $S$, 
i. e., $G_\beta =(\beta -A)^{-1}$, $\beta >0$. If the semigroup $(T_t)_{t 
\ge 0}$ is contractive in $L^2(E,{\bnu})$, because of $\langle T_tu\, ,\, 
u\rangle^2\le\langle T_tu\, ,\, T_tu\rangle\langle u\, ,\, u\rangle $, it 
holds that  
\begin{eqnarray*}
\langle T_tu\, ,\, u\rangle\le\langle u\, ,\, u\rangle\, . 
\end{eqnarray*}
In this way one shows positivity of the form $S$, that is $S(u,u)\ge 0$ 
for all $u\in D(S)$. This observation is crucial for the whole concept of 
Mosco type convergence of sequences $S_n$ of forms on sequences of spaces 
$L^2(E_n,\bnu_n)$ to a limiting form $S$ on $L^2(E,\bnu)$ as $n\to\infty$. 
However, in \cite{Lo14-1} we have developed a framework of Mosco type 
convergence of sequences of forms when contractivity is replaced by a 
technical condition on $A_n'\1$, $n\in {\Bbb N}$, and $A'\1$ where $\1$ is 
the constant function taking the value one and the $'$ refers to the dual 
generator. 
\bigskip 

\nid 
{\bf Mosco type convergence of non-symmetric, positive bilinear forms.} For 
every $n\in {\Bbb N}$, let $(T_{n,t})_{t\ge 0}$ be a strongly continuous 
contraction semigroup in $L^2(E,\bnu_n)$ and let $(T_t)_{t\ge 0}$ be a strongly 
continuous contraction semigroup in $L^2(E,\bnu)$. Denote by $S_n$, $A_n$, 
$(G_{n,\beta})_{\beta >0}$ the bilinear form in the sense of of the above 
paragraph {\it Bilinear forms}, the generator, and the family of resolvents 
associated with $(T_{n,t})_{t\ge 0}$, $n\in {\Bbb N}$. Similarly, let $S$, 
$A$, and $(G_\beta)_{\beta >0}$ the bilinear form, the generator, and the 
family of resolvents associated with $(T_t)_{t\ge 0}$. For $\beta >0$, $n 
\in {\Bbb N}$, $\vp_n\in D(S_n)$, $\psi_n\in L^2(E,{\bnu}_n)$, set $S_{n, 
\beta}(\vp_n,\psi_n):=\beta\langle\vp_n\, ,\, \psi_n \rangle_n + S_n(\vp_n, 
\psi_n)$, and for $\vp\in D(S)$, $\psi\in L^2(E,{\bnu})$, define $S_{\beta} 
(\vp,\psi):=\beta\langle\vp\, ,\, \psi\rangle + S(\vp,\psi)$. 
\begin{definition}\label{Definition2.4} 
{\rm We say that $S_n$, $n\in {\Bbb N}$, {\it pre-converges} to $S$ if 
\begin{itemize} 
\item[(i)] For every $\vp\in L^2(E,{\bnu})$ and every subsequence 
$\vp_{n_k}\in D(S_{n_k})\cap {\cal C}$, $k\in {\Bbb N}$, $w$-converging 
to $\vp$ such that $\sup_{k\in {\Bbb N}}\left\langle A_{n_k}\vp_{n_k}\, 
,\, A_{n_k}\vp_{n_k}\right\rangle_{n_k}<\infty$, we have 
\begin{eqnarray*}
S(\vp,\vp )\le \liminf_{k\to\infty}S_{n_k}(\vp_{n_k},\vp_{n_k} )\, . 
\end{eqnarray*}
\item[(ii)] For every $\psi\in D(S)$, there exists a sequence $\psi_n 
\in D(S_n)\cap {\cal C}$, $n\in {\Bbb N}$, $s$-converging to $\psi$ 
such that $\sup_{n\in {\Bbb N}}\left\langle A_n\psi_n\, ,\, A_n\psi_n 
\right\rangle_n<\infty$ and 
\begin{eqnarray*}
\limsup_{n\to\infty} S_n(\psi_n ,\psi_n )\le S(\psi ,\psi )\, .  
\end{eqnarray*}
\end{itemize} 
}
\end{definition}
\begin{lemma}\label{Lemma2.6} 
Let $S_n$, $n\in {\Bbb N}$, be a sequence of bilinear forms 
pre-convergent to $S$. Furthermore, let $\beta >0$ and let $u_n\in 
D(S_n)\cap {\cal C}$ such that $A_nu_n\in {\cal C}$, $n\in {\Bbb N}$, be a 
$w$-convergent sequence with $\, \sup_{n\in {\Bbb N}}\langle A_nu_n\, ,\, 
A_nu_n\rangle_n<\infty$. Let $u\in D(S)$. Introduce the following 
conditions. 
\begin{itemize}
\item[{\rm (iii)}] Let $u_n$, $n\in {\Bbb N}$, and $u$ as above.  
\begin{eqnarray*}
\lim_{n\to\infty}S_{n,\beta}(u_n,\psi_n)=S_\beta (u,\psi)
\end{eqnarray*}
for all $\psi\in {\cal C}$ and all sequences $\psi_n\in {\cal C}$, $n\in 
{\Bbb N}$, $s$-convergent to $\psi$ yields 
\begin{eqnarray*}
\lim_{n\to\infty}S_{n,\beta}(\psi_n,u_n)=S_\beta (\psi ,u)  
\end{eqnarray*}
for all $\psi\in D(S)$ and all sequences $\psi_n\in D(S_n)\cap {\cal 
C}$, $n\in {\Bbb N}$, $s$-convergent to $\psi$ in the sense of condition 
(ii) in Definition \ref{Definition2.4}. 
\item[{\rm (iv)}] Let $u_n$, $n\in {\Bbb N}$, and $u$ as above. If 
$\beta u_n-A_nu_n\wstack{n\to\infty}{\lra}\beta u-Au$ then $u_n\wstack 
{n\to\infty}{\lra}u$. 
\end{itemize} 
(a) Condition (iv) implies (iii). \\ 
(b) Suppose 
\begin{itemize}
\item[(c3)] 
\begin{itemize}
\item[$(i)$] ${\cal G}:=\{G_\beta g:g\in {\cal C},\ \beta >0\}\subseteq 
{\cal C}$ in the sense that for every $g \in {\cal C}$ and $\beta>0$, 
there is a $u\in {\cal C}$ with $G_\beta g=u$ ${\bnu}$-a.e.
\item[$(ii)$] ${\cal G}_n:=\{G_{n, \beta}g:g\in {\cal C},\ \beta >0\} 
\subseteq {\cal C}$, $n\in {\Bbb N}$, in the sense that for every $g 
\in {\cal C}$, $\beta>0$, and every $n\in {\Bbb N}$, there exists a $v 
\in {\cal C}$ such that $G_{n,\beta}g=v$ ${\bnu}_n$-a.e. 
\end{itemize} 
\end{itemize} 
Then (iii) implies (iv).
\end{lemma}
\medskip 

\nid
{\bf Definition 2.4 continued } Let $S_n$, $n\in {\Bbb N}$, be a 
sequence of bilinear forms pre-convergent to $S$. If, in addition, 
condition (iii) is satisfied, then we say that $S_n$, $n\in {\Bbb N}$, 
{\it converges} to $S$. 
\medskip 

\begin{theorem}\label{Theorem2.7}
Let $\beta>0$, suppose that conditions (c1)-(c3) as well as (c3) for 
$G_\beta'$ and $G_{n,\beta}'$ are satisfied. Assume that $S_n$, $n\in 
{\Bbb N}$, converges to $S$ in the sense of Definition \ref{Definition2.4}. 
 \\ 
(a) For all ${\cal C}\ni f_n\wstack{n\to\infty}{\lra}f\in L^2(E,\bnu)$
we have $G_{n,\beta}f_n\wstack{n\to\infty}{\lra}G_\beta f$ and $G_{n, 
\beta}'f_n\wstack{n\to\infty}{\lra}G_\beta'f$.  \\ 
(b) For all ${\cal C}\ni g_n\sstack{n\to\infty}{\lra}g\in L^2(E,\bnu)$ 
we have $G_{n,\beta}g_n\sstack{n\to\infty}{\lra}G_\beta g$ and $G_{n,\beta}' 
g_n\sstack{n\to\infty}{\lra}G_\beta' g$. 
\end{theorem}
\medskip 

\noindent 
{\bf Mosco type convergence of non-symmetric, non-positive bilinear forms.} 
Let us drop the assumption that the semigroups $(T_{n,t})_{t\ge 0}$, $n\in 
{\Bbb N}$, and $(T_t)_{t\ge 0}$ are contractive. Anything else for the 
semigroups remains as introduced above. 

In particular, let us assume that there are Markov processes associated 
with the semigroups $(T_{n,t})_{t\ge 0}$, $n\in {\Bbb N}$, and $(T_t)_{t 
\ge 0}$: For $n\in {\Bbb N}$, let $X^n=((X^n_t)_{t\ge 0},(P^n_\mu)_{\mu 
\in E_n})$ be a process taking values in $E_n$ which corresponds to 
the semigroup $(T_{n,t})_{t\ge 0}$ and the form $S_n$. Furthermore, 
let $X=((X_t)_{t\ge 0},(P_\mu)_{\mu\in E})$ be a process associated 
with the semigroup $(T_t)_{t\ge 0}$ and the form $S$ which takes 
values in some subset of $E$. Suppose that the paths of the processes 
$X^n$, $n\in {\Bbb N}$, and $X$ are cadlag. For $\beta >0$, introduce 
$G_{n,\beta}g_n:=\int_0^\infty e^{-\beta t}T_{n,t}g_n\, dt$, $g_n\in 
L^\infty (E,\bnu_n)$, $n\in {\Bbb N}$, $G_{\beta}g:=\int_0^\infty e^{ 
-\beta t}T_{t}g\, dt$, $g\in L^\infty (E,\bnu)$. Since the semigroups 
$(T_{n,t})_{t\ge 0}$, $n\in {\Bbb N}$, and $(T_{t})_{t\ge 0}$ are not 
necessarily contractive, the associated families of resolvents $(G_{n, 
\beta})_{\beta >0}$, $n\in {\Bbb N}$, and $(G_\beta)_{\beta >0}$ may 
not directly be well-defined on the corresponding $L^2$-spaces. 
\begin{itemize}
\item[(c6)] 
\begin{itemize} 
\item[$(i)$] $\1\in D(A')$ and $\1\in D(A'_n)$, $n\in {\Bbb N}$, and $\sup_{ 
n\in {\Bbb N}}\|A_n'\1\|_{L^\infty (E,\sbnu_n)}<\infty$. 
\item[$(ii)$] $T_t'\1\in L^\infty (E,\bnu)$ and the limit $A'\1=\lim_{t\to 
0}\frac1t(T_t'\1-\1)$ exists in $L^\infty (E,\bnu)$. 
\item[$(iii)$] There exist $N_n\in {\cal B}(E_n)$ with $\bnu_n(N_n)\stack{n 
\to\infty}{\lra}0$ such that for $\vp\in {\cal C}$ there is $\Phi_n\equiv 
\Phi_n(\vp)\in {\cal C}$ with $\Phi_n=A_n'\1\cdot\vp$ on $E_n\setminus N_n$. 
Furthermore, $\Phi_n\sstack{n\to\infty}{\lra}A'\1\cdot\vp$. 
\item[$(iv)$] For $n\in {\Bbb N}$, there exists a set $D_n\subseteq D(S_n) 
\cap L^\infty (E,\bnu_n)$ which is dense in $D(S_n)$ with respect to the 
norm $\|f\|_{D_n}:=\left(\langle f\, , \, f\rangle_n+\langle A_nf\, , \, 
A_nf\rangle_n\right)^{1/2}$. 
\end{itemize}
\end{itemize}
Define $D(\hat{S}_n):=D(S_n)$, $n\in {\Bbb N}$, $D(\hat{S}):=D(S)$, and 
\begin{eqnarray*}
\hat{S}_n(u_n,v_n):=S_n(u_n,v_n)+{\textstyle\frac12}\langle A_n'\1\cdot u_n 
\, , \, v_n\rangle_n\, , \quad u_n\in D(\hat{S}_n),\ v_n\in L^2(E,\bnu_n),\ 
n\in {\Bbb N}, \ \ 
\end{eqnarray*} 
and 
\begin{eqnarray*}
\hat{S}(u,v):=S(u,v)+{\textstyle\frac12}\langle A'\1\cdot u\, , \, v\rangle 
\, , \quad u\in D(\hat{S}),\ v\in L^2(E,\bnu). 
\end{eqnarray*} 
\begin{lemma}\label{Lemma2.11} 
Suppose (c6). We have $\hat{S}_n(u_n,u_n)\ge 0$, $u_n\in D(S_n)$, $n\in 
{\Bbb N}$, and $\hat{S}(u,u)\ge 0$, $u\in D(S)$. 
\end{lemma} 
Set $D(\hat{A}_n):=D(A_n)$, $n\in {\Bbb N}$, $D(\hat{A}):=D(A)$, and 
$\hat{A}_nu_n:=A_nu_n-{\textstyle\frac12} A_n'\1\cdot u_n$, $u_n\in 
D(\hat{A}_n)$, $n\in {\Bbb N}$, and 
\begin{eqnarray}\label{2.15}
\hat{A}u:=Au-{\textstyle\frac12} A'\1\cdot u\, , \quad u\in D(\hat{A}),  
\end{eqnarray} 
\begin{eqnarray*}
\hat{T}_tv={\Bbb E}_\cdot\left(\exp\left\{-\int_0^t{\textstyle\frac12}A' 
\1(X_s)\, ds\right\} v(X_t)\right)\, ,\quad v\in L^2(E,\bnu),\ t\ge 0. 
\end{eqnarray*} 
Similarly define $(\hat{T}_{n,t})_{t\ge 0}$. 
In addition, let $(\hat{G}_{n,\beta})_{\beta\ge 0}$, denote the resolvent 
associated with $\hat{A}_n$, $\hat{S}_n$, $(\hat{T}_{n,t})_{t\ge 0}$, $n\in 
{\Bbb N}$, and let $(\hat{G}_\beta)_{\beta\ge 0}$, denote the resolvent 
associated with $\hat{A}$, $\hat{S}$, $(\hat{T}_{t})_{t\ge 0}$. 
\begin{itemize} 
\item[(c3')] 
\begin{itemize} 
\item[$(i)$] 
If ${\cal C}\subseteq L^\infty (E,\bnu)$ then $\{G_\beta g:g\in L^\infty 
(E,\bnu)\, , \ \beta >0\}\subseteq {\cal C}$ in the sense that for every 
$g\in L^\infty (E,\bnu)$, there is a $u\in {\cal C}$ with $G_\beta g=u$ 
${\bnu}$-a.e. Otherwise, $D(S)\subseteq {\cal C}$. 
\item[$(ii)$] 
If, for $n\in {\Bbb N}$, ${\cal C}\subseteq L^\infty (E,\bnu_n)$ then $\{ 
G_{n,\beta} g:g\in L^\infty (E,\bnu_n)\, , \ \beta >0\}\subseteq {\cal C}$ 
in the sense that for every $g\in L^\infty (E,\bnu_n)$, there is a $u\in 
{\cal C}$ with $G_{n,\beta} g=u$ ${\bnu_n}$-a.e. Otherwise, $D(S_n)\subseteq 
{\cal C}$. 
\end{itemize} 
\end{itemize} 
\begin{lemma}\label{Lemma2.13} 
Suppose (c6). (a) If (c3'(i)) then condition (c3(i)) holds for $\hat{G}_\beta$ 
in place of $G_\beta$. \\ 
(b) If (c3'(ii)) then condition (c3(ii)) holds for $\hat{G}_{n,\beta}$, in 
place of $G_{n,\beta}$, $n\in {\Bbb N}$. 
\end{lemma} 
\begin{theorem}\label{Theorem2.14} 
Let $\, C:={\T\frac12}\|A'\1\|_{L^\infty(E,\sbnu)}\vee\, \sup_{n\in {\Bbb N}} 
{\T\frac12}\|A_n'\1\|_{L^\infty (E,\sbnu_n)}$. Suppose (c1)-(c3) as well as (c3) 
for $G_{n,\beta}'$, $n\in {\Bbb N}$, and $G_\beta'$, and (c6). 
Furthermore, suppose (c3) for $\hat{G}_{n,\beta}$, $n\in {\Bbb N}$, and $\hat 
{G}_\beta$ as well as $\hat{G}_{n,\beta}'$, $n\in {\Bbb N}$, and $\hat{G}_\beta'$. 
Assume that the forms $\hat{S}_n$, $n\in {\Bbb N}$, converge to the form $\hat
{S}$ in the sense of Definition \ref{Definition2.4}. 
\medskip 

For all $\beta>C$, the operators $G_{n,\beta}$ and $G_\beta$ can be continuously 
extended to operators $G_{n,\beta}:L^2(E,\bnu_n)\to L^2(E,\bnu_n)$, $n\in {\Bbb 
N}$, and $G_\beta :L^2(E,\bnu)\to L^2(E,\bnu)$, respectively. For all ${\cal C} 
\ni f_n\wstack{n\to\infty}{\lra}f\in L^2(E,\bnu)$, all ${\cal C}\ni g_n\sstack{n 
\to\infty}{\lra}g\in L^2(E,\bnu)$, we have $G_{n,\beta}f_n\wstack{n\to\infty} 
{\lra}G_\beta f$, $G_{n,\beta}'f_n\wstack{n\to\infty}{\lra}G_\beta'f$ and $G_{n, 
\beta}g_n\sstack{n\to\infty}{\lra}G_\beta g$,  $G_{n,\beta}'g_n\sstack{n\to\infty} 
{\lra}G_\beta g$. 
\end{theorem} 
{\bf Remark} (1) It is sufficient to require (c3) (for $G_{n,\beta}$, $n\in {\Bbb N}$, 
and $G_\beta$) and (c3) for $\hat{G}_{n,\beta}$, $n\in {\Bbb N}$, and $\hat{G}_\beta$ 
if one is not interested in the convergence of $G_{n,\beta}'$. 
\medskip 

\nid 
(2) Suppose (c1)-(c3) and ${\cal T}:=\{T_t g:g\in {\cal C},\ \beta >0\}\subseteq 
{\cal C}$, ${\cal T}_n:=\{T_{n,t}g:g\in {\cal C},\ \beta >0\}\subseteq {\cal C}$, 
$n\in {\Bbb N}$, in the sense of condition (c3).

It has been demonstrated and mentioned in \cite{Lo14-1}, Remark (10) to Theorem 
2.14 (c), that for all $g\in {\cal C}$, ${\cal C}\ni g_n\sstack{n\to\infty}{\lra} 
g$, and $\beta> {\T\frac12}\|A'\1\|_{L^\infty(E,\sbnu)}\vee\, \sup_{n\in {\Bbb N}} 
{\T\frac12}\|A_n'\1\|_{L^\infty (E,\sbnu_n)}$, we have $G_{n,\beta}g_n\sstack{n\to 
\infty}{\lra}G_\beta g$ iff $G_{n,\beta}g\sstack{n\to\infty}{\lra}G_\beta g$ iff 
$\ T_{n,t}g\sstack{n\to\infty}{\lra} T_tg$ iff $\ T_{n,t}g_n\sstack{n\to\infty} 
{\lra} T_tg$. 
  
\subsection{Relative Compactness} 

{\bf Specification of Subsection 2.1, notation for the remainder.} For this, 
let $D$ be a bounded $d$-dimensional domain for some $d\in {\Bbb N}$. Let 
${\cal M}_\partial(\overline{D})$ be the set of all equivalence classes $\mu$ 
such that $m_1,m_2\in\mu$ implies $m_1|_D=m_2|_D$. From now on, throughout the 
paper, we will assume that $E$ is the space ${\cal M}_\partial(\overline{D})$. 
Here, we identify all points belonging to $\partial D$ with each other. By 
$r(x,y):=|x-y|\wedge\left(\inf_{b\in\partial D}|b-x|+\inf_{b\in\partial D}|b-y| 
\right)$ if $x,y\in D$ and $r(x,\partial D)=r(\partial D,x):=\inf_{b\in\partial 
D}|b-x|$ if $x\in D$, as well as $r(\partial D,\partial D):=0$ the space $(D 
\cup\partial D,r)$ becomes a separable, complete, and compact metric space, cf. 
also \cite{Ma89}. Furthermore, continuity on $D$ with respect to $r$ coincides 
with continuity with respect to the Euclidean metric and $\{f\in C(\overline{D})
:f$ constant on $\partial D\}$ is the set of all continuous functions on $(D\cup 
\partial D,r)$. 

Let ${\cal M}_\partial(\overline{D})$ be endowed with the Prokhorov metric. We 
note that in this way ${\cal M}_\partial(\overline{D})$ is a separable, complete, 
and compact space.
 
In addition, for $n\in {\Bbb N}$, let $E_n'$ be the set of all measures $\mu$ 
in $E$ of the form $\mu=\frac{1}{n}\sum_{i=1}^n\delta_{z_i}$ where $z_1,\ldots 
,z_n \in\overline{D}$ and $\delta_z$ denotes the Dirac measure concentrated at 
$z$. Furthermore, let $E_1:=E_1'$ and $E_{n+1}:=E_{n+1}'\setminus\bigcup_{i=1 
}^nE_i$, $n\in {\Bbb N}$, and $E_0:=E\setminus\bigcup_{n=1}^\infty E_n$. 
According to the basic setting of Subsection 2.1 $E_n$ and $E_n'$ differ by 
$\bnu_n$-zero set, $n\in {\Bbb N}$. It is therefore reasonable to {\it identify} 
$L^p(E,\bnu_n)$ {\it with both} $L^p(E_n,\bnu_n)$ and $L^p(E_n',\bnu_n)$, $1\le 
p\le\infty$, $n\in {\Bbb N}$. 

To be consistent with \cite{Lo13}, we will keep on writing $C_b(E)$ for $C(E)$. 
Choose ${\cal F}:=C_b(E)$ and note that therefore ${\cal C}$ is now the space of 
all functions $\vp\in {\cal D}$ satisfying the following.  
\begin{itemize} 
\item[(c1')] $\vp$ is bounded and continuous on $E_n$, $n\in {\Bbb N}$. 
\item[(c2')] $\langle \vp\, , \, \psi\rangle_n\stack{n\to\infty}{\lra} 
\langle\vp\, , \, \psi\rangle$ for all $\psi\in C_b(E)$.  
\end{itemize}
We observe that with ${\cal F}=C_b(E)$, 
\begin{eqnarray}\label{2.2} 
{\cal F}\subseteq{\cal D}\quad\mbox{\rm if and only if}\quad\bnu_n\stack{n 
\to\infty}{\Ra}\bnu\, . 
\end{eqnarray} 
Let us assume the latter and note that (c1') and (c2') are now the defining 
properties of ${\cal C}\subseteq {\cal D}$. 
\medskip 

As in Subsection 2.1, let us assume that there are Markov processes $X^n_t  
=\frac1n\sum_{i=1}^n\delta_{(X_j)_t}$ and $X_t$, $t\ge 0$, associated with 
the semigroups $(T_{n,t})_{t\ge 0}$, $n\in {\Bbb N}$, and $(T_t)_{t\ge 0}$. 

Define the measures $P_{{\sbnu}_n}:=\int_E P^n_\mu\, {\bnu}_n(d\mu)$, 
$n\in{\Bbb N}$, and $P_{\sbnu}:=\int_E P_\mu(\mu )\, {\bnu}(d\mu)$, and 
introduce the processes ${\bf X}^n=((X^n_t)_{t\ge 0},P_{{\sbnu}_n})$ 
and ${\bf X}=((X_t)_{t\ge 0},P_{\sbnu})$. Moreover, let $\, {\Bbb E}^n_{ 
\mu}$ be the expectation corresponding to $P^n_\mu$, $\mu\in E_n$, and 
let $\, {\Bbb E}_{{\sbnu}_n}$ be the expectation corresponding to $P_{ 
{\sbnu}_n}$, $n\in {\Bbb N}$. Let us introduce the set of test functions 
we are going to work with in this section. Suppose the following.
\begin{itemize}
\item[(c7)] There exists an algebra $\t C_b(E)\subseteq C_b(E)$ of 
everywhere on $E$ defined functions with $\t C_b(E)\subseteq{\cal G}$ 
in the sense that, for every $f\in \t C_b(E)$, there is a $g\equiv g(f) 
\in {\cal C}$ and a $\beta>0$ with $f=\beta G_\beta g$ ${\bnu}$-a.e. 
$\t C_b(E)$ contains the constant functions and separates points in $E$.
\end{itemize}

For $f\in\t C_b(E)$, $g=g(f)\in {\cal C}$, and a given sequence $\ve_n>0$, 
$n\in {\Bbb N}$, introduce 
\begin{eqnarray}\label{3.1}
B:=\bigcup_{n=1}^\infty\left\{\mu \in E_{n}:\left|\beta G_{n,\beta} 
g(\mu )-f(\mu )\right|\ge\ve_n\|g\|\right\}\, . 
\end{eqnarray}
Furthermore, let $\tau_{B^c}\equiv\tau_{B^c}^n(g)$ denote the first exit 
time of ${\bf X}^{n}$ from the set $B^c\cap E_{n}$, $n\in {\Bbb N}$. 
Let $T>0$ and set 
\begin{eqnarray*}
\gamma_n\equiv\gamma_n(f):=\sup_{s\in [0,T+1]}\left|\beta G_{n,\beta} 
g(X^n_s)-f(X^n_s)\right|\, , \quad n\in {\Bbb N}. 
\end{eqnarray*}
\begin{itemize}
\item[(c8)] There is a sequence $\ve_n>0$, $n\in {\Bbb N}$, with $\ve_n\stack 
{n\to\infty}{\lra}0$ such that with $B\equiv B((\ve_n)_{n\in {\Bbb N}})$ 
defined in (\ref{3.1})
\begin{eqnarray*}
\, {\Bbb E}_{{\sbnu}_n}\left(e^{-\beta\tau_{B^c}}\right)\stack{n\to\infty} 
{\lra}0
\end{eqnarray*}
whenever $\langle f-\beta G_{n,\beta}g\, ,\, f-\beta G_{n,\beta}g\rangle_{n}
\stack{n\to\infty}{\lra}0$. 
\end{itemize}
\begin{theorem}\label{Theorem3.2} 
(a) Let the following be satisfied: 
\begin{itemize} 
\item[(i)] Conditions (c3), (c7), and (c8) hold. 
\item[(ii)] We have the hypotheses of Theorem \ref{Theorem2.14}, namely 
\begin{itemize} 
\item[{}] (c3) for $\hat{G}_{n,\beta}$, $n\in {\Bbb N}$, and $\hat{G_\beta}$ 
in place of $G_{n,\beta}$, $n\in {\Bbb N}$, and $G_\beta$, 
\item[{}] (c6), 
\item[{}] the forms $\hat{S}_n$, $n\in {\Bbb N}$, converge to the form 
$\hat{S}$ in the sense of Definition \ref{Definition2.4}. 
\end{itemize}
\end{itemize}
Then, for $f\in \t C_b(E)$, the family of processes $f({\bf X}^n)=((f( 
X^n_t))_{t\ge 0},P_{{\sbnu}_n}\circ f^{-1})$, $n\in{\Bbb N}$, is 
relatively compact with respect to the topology of weak convergence of 
probability measures over the Skorokhod space $D_{[-\|f\|,\|f\|]}([0, 
\infty))$. \\ 
(b) The family of processes ${\bf X}^n=((X^n_t)_{t\ge 0}$, $P_{{\sbnu}_n} 
)$, $n\in {\Bbb N}$, is relatively compact with respect to the topology 
of weak convergence of probability measures over the Skorokhod space $D_E 
([0,\infty))$. 
\end{theorem}

\subsection{Infinite Dimensional Integration by Parts} 

Let $F$ denote the space of all finite signed measures on $\overline{D}$. 
Let $h_1,h_2,\ldots\, $ be the eigenfunctions of the Dirichlet Laplacian 
on $D$ corresponding to the eigenvalues $0>2\lambda_1>2\lambda_2\ge 
\ldots\, $, normalized in $L^2(D)$ such that $h_1>0$. For $t\in [0,1]$ 
define the space 
\begin{eqnarray*} 
H(t)=\left\{h\, dx\in F:\sum_{n=1}^\infty\lambda_n^2e^{-2t\lambda_n}\cdot 
(h_n,h)^2<\infty\right\}
\end{eqnarray*}
which becomes a Hilbert space with the norm $\|h\, dx\|_{H(t)}=(\sum_{n 
=1}^\infty\lambda_n^2e^{-2t\lambda_n}\cdot (h_n,h)^2)^{1/2}$, $h\, dx 
\in H(t)$. Abbreviate $H\equiv H(0)$. We mention that for all $h\, dx\in 
H$ it holds that  
\begin{eqnarray}\label{2.4} 
\frac12\Delta h=\sum\lambda_i(h_i,h)_{L^2(D)}h_i\quad\mbox{\rm in}\quad 
L^2(D)\, . 
\end{eqnarray}
If $\mu\equiv h\, dx\in H$, we also will write $\Delta\mu$ for $\Delta h 
\, dx$. 
\medskip 

Denote by $p(t,x,y)$, $t\ge 0$, $x,y\in D$, the transition density 
function of a Brownian motion on $D$ killed when hitting $\partial D$. 
For $h\in L^2(D)\setminus\{0\}$, $h\, dx\equiv\mu$ and $t\ge 0$, set 
$|\mu|(D):=\int |h|\, dx$ and 
\begin{eqnarray*}
u(t,y):=\int_{x\in D}p(t,x,y)\, \mu(dx)\, , \quad z(t)\equiv z(\mu ,t):= 
\frac{1}{|\mu|(D)}\int_{y\in D}u(t,y)\, dy\, . 
\end{eqnarray*}
By convention, we have $z(\mu,0)=\int h\, dx/|\mu|(D)$. For $h\in L^2(D) 
\setminus\{0\}$, $h\, dx\equiv\mu$, and $t\ge 0$ define 
\begin{eqnarray*} 
v(t,y)\equiv v(\mu,t,y):=\frac{1}{z(\mu ,t)}u(t,y)\, ,  
\quad y\in D, 
\end{eqnarray*}
which is 
\begin{eqnarray*}
v(t,y)\, dy=\frac{1}{z(\mu ,t)}\int_{x\in D}p(t,x,y)\, \mu(dx)\, dy\, . 
\end{eqnarray*}
Furthermore, if $h=0$ we set $u(t,z)=v(t,y)=0$ and $z(\mu ,t)=1$, $y\in 
D$, $t\ge 0$. Let $\mu\in E\cap H$, $t\ge 0$, and $y\in D$. Comparing with 
\cite{GK04}, we observe that $v(t,y)$ satisfies 
\begin{eqnarray*}
\frac{\partial}{\partial t}v(t,x)=\frac12\Delta v(t,x)-\frac{z'(\mu ,t)} 
{z(\mu,t)}v(t,x)\, ,\quad v(t,x)|_{x\in \partial D}=0,\ v(t,x)\, dx\stack 
{t\to 0}{\Ra}\mu\, , 
\end{eqnarray*}
where $\Delta$ is the Laplace operator on $D$ and $\Ra$ indicates weak 
convergence of finite signed measures in the sense of $\mu_n{\Ra}\mu$ if 
$\int f\,d\mu_n\stack{n\to\infty}{\lra}\int f\,d\mu$ for all bounded and 
continuous test functions $f$. In particular we mention that, for $t=0$, 
$z'(0)$ is the right derivative. Furthermore for $\mu\equiv h\, dx\in E 
\cap H$, 
\begin{eqnarray}\label{2.5} 
z'(0)\equiv z'(\mu,0)&&\hspace{-.5cm}=\int_{y\in D}\left.\frac{d}{dt} 
\right|_{t=0}\int_{x\in D}p(t,x,y)\, h(x)\, dx\, dy\nonumber \\ 
&&\hspace{-.5cm}=\int_{y\in D}\frac12(\Delta h)(y)\, dy  
\end{eqnarray}
which gives 
\begin{eqnarray*}
\left|z'(\mu,0)\right|\le\frac12|D|^{\frac12}\left\|\Delta h\right\|_{ 
L^2(D)}<\infty\, , 
\end{eqnarray*}
$|D|$ denoting the Lebesgue measure of $D$. 
\medskip

Let $C_b(F)$ denote the space of all bounded continuous real functions 
on $F$ and recall that $C_b(E)$ denotes the space of all continuous real 
functions on the compact space $E$. Moreover, let $C_b^{2,1}(F,E)$ (with 
respect to $\bnu$) denote the set of all $f\in C_b(F)$ with the following 
properties. For $\bnu$-a.e. $\mu\in E$ and all $h\in L^v(D)$ the 
directional derivative 
\begin{eqnarray*}
\frac{\partial f}{\partial h}(\mu):=\lim_{t\to 0}\frac1t\left(f(\mu+th 
\, d\lambda)-f(\mu)\right)
\end{eqnarray*}
exists, a measurable map $Df:E\to L^2(D)$ exists such that 
\begin{eqnarray*}
(Df(\mu),h)=\frac{\partial f}{\partial h}(\mu)\quad\mbox{\rm for }\bnu 
\mbox{\rm -a.e. }\mu\in E\ \mbox{\rm and all }h\in L^2(D),  
\end{eqnarray*}
and $Df\in L^2(E,\bnu;L^2(D))$. Define $C_b^{2,1}(E):=\{f|_E:f\in C_b^{ 
2,1}(F,E)\}$. Denote $\Phi(f):=\{\vp\in C_b^{2,1}(F,E):f=\vp|_E\}$, $f 
\in C_b^{2,1}(E)$, and let $A\in L^2(E,\bnu;L^2(D))$ be a vector field 
with $(A,\1)=0$ $\bnu$-a.e. According to \cite{Lo14-2}, Section 2, the 
expression $(D\vp,A)\in L^1(E,\bmu)$ is for fixed $f\in C_b^{2,1}(E)$ 
independent of $\vp\in\Phi(f)$. In this sense we call $Df$ the {\it 
gradient} of $f\in C_b^{2,1}(E)$. 
\medskip 

As discussed in \cite{Lo14-2}, Section 4, the following three conditions 
are reasonable.  
\begin{itemize}
\item[(i)] $\bnu$ is concentrated on all $\mu\in E\cap H(1)$ for which 
$\sum_{j=1}^\infty e^{-\lambda_j}(h_j,\mu)h_j\ge 0$. 
\item[(ii)] $\left\|\mu\right\|^2_{H(1)}\in L^2(E,\bnu)$.  
\item[(iii)] $(z(\cdot ,1))^{-2}\in L^2(E,\bnu)$. 
\end{itemize} 
According  to \cite{Lo14-2}, Sections 1 and 4, there is the following 
flow $U^+(t,\mu)$ for $\bnu$-a.e. $\mu$ and $t\ge -1$. In particular, 
for $t\ge 0$ and $\bnu$-a.e. $\mu$, 
\begin{eqnarray*}
U^+(t,\mu)(dy)=v(\mu,t,y)\, dy =(z(\mu,t))^{-1}\sum_{j=1}^\infty e^{ 
\lambda_j t}(h_j,\mu)h_j\, dy 
\end{eqnarray*} 
and, for $t\in [-1,0]$ and $\bnu$-a.e. $\mu$, there is a $\hat{\mu}\in E 
\cap L^2(D)$ such that $U^+(-t,\hat{\mu})=\mu$. By setting $U^+(t,\mu):= 
\hat{\mu}$ for $t\in [-1,0]$ and $\bnu$-a.e. $\mu$, the map $U^+(t,\cdot 
)$ is then given on $t\ge -1$ $\bnu$-a.e. Finally, we refer to the {\it 
flow property}. For $s,t\ge -1$ such that $s+t\ge -1$ and $\bnu$-a.e. 
$\mu$, the composition $U^+(t,\cdot)\circ U^+(s,\mu)$ is well defined and 
we have 
\begin{eqnarray*}
U^+(t,\cdot)\circ U^+(s,\mu)=U^+(s+t,\mu)\, . 
\end{eqnarray*} 
We note that for $\mu\in E\cap H(1)$ and $t\ge -1$ the derivative $\frac 
{d}{dt}\left(U^+(t,\mu)(dx)/dx\right)$ exists in $L^2(D)$ where for $t= 
-1$ this is a derivative from the right. We set $A^f\mu:=\left[\frac{d} 
{dt}\left(U^+(t,\mu)(dx)/dx\right)\right]\, dx$. For the subsequent 
quasi-invariance result we introduce the following conditions. 
\begin{itemize}
\item[(j)] $\D\frac{d A^f}{dx}\in L^2(E,\bnu;L^2(D))$. 
\item[(jj)] There exists a unique {\it divergence} of the vector field 
$A^f$ relative to $\bnu$ and the gradient $D$, that is an element $\delta 
(A^f)\in L^2(E,\bnu)$ satisfying 
\begin{eqnarray*}
\int (Df,A^f)\, d\bnu=-\int f\, \delta (A^f)\, d\bnu\, ,\quad f\in C^{2,1 
}_b(E). 
\end{eqnarray*}
\item[(jjj)] $\delta (A^f)\in L^\infty(E,\bmu)$. 
\end{itemize} 

The specification to the particular non-linear PDE above allows a more direct 
formulation of the subsequent result than in \cite{Lo14-2}, Section 2. 
\begin{proposition}\label{Proposition4.1} 
Assume (i)-(iii) and (j)-(jjj) of the present subsection. (a) All measures 
$\bnu\circ U^+(-t,\mu)$, $0\le t\le 1$, are equivalent, note that $\bnu\circ 
U^+(0,\mu)=\bnu$. The Radon-Nikodym derivatives have a version such that 
\begin{eqnarray*} 
[0,1]\ni t\to r_{-t}:=\frac{d\bnu\circ U^+(-t,\mu)}{d\bnu} 
\end{eqnarray*}
is $\bnu$-a.e. absolutely continuous on $([0,1],{\cal B}([0,1]))$. We have 
\begin{eqnarray*}
r_{-t}(\mu)=\exp\left\{-\int_{s=0}^t\delta (A^f)(U^+(-s,\mu))\, ds\right\} 
\, ,\quad 0\le t\le 1,\quad\bnu\mbox{\rm -a.e.} 
\end{eqnarray*}
(b) For $f\in L^\infty(E,\bnu)$, we have 
\begin{eqnarray*} 
\left.\frac{d}{dt}\right|_{t=0}\int f(U^+(t,\mu))\, \bnu (d\mu )=-\int f(\mu) 
\delta (A^f)(\mu )\, \bnu (d\mu)\, . 
\end{eqnarray*}
\end{proposition}

Introduce  
\begin{eqnarray}\label{4.3}
\t C_b^2(E):=\left\{f(\mu)=\vp((h_1,\mu),\ldots,(h_r,\mu)),\ \mu\in E: 
\vphantom{\t f}\vp\in C_b^2({\Bbb R}^r),\ r\in {\Bbb N}\right\}\, . 
\end{eqnarray}
Let $C_0^2({\Bbb R})$ denote the set of all $f\in C^2({\Bbb R})$ which 
have compact support. Furthermore define $K$ to be the set of all 
non-negative $k\in C(D)\cap L^2(D)$ such that $\lim_{D\ni u\to v}k(u)= 
\infty$, $v\in\partial D$. 
\begin{eqnarray}\label{4.4}
&&\hspace{-.5cm}\t C_0^2(E):=\left\{f(\mu)=\vp((h_1,\mu),\ldots,(h_r, 
\mu))\cdot\vp_0((k,\mu)),\ \mu\in E:\vphantom{\t f}\right. \nonumber \\ 
&&\hspace{5.5cm}\left.\vphantom{\t f}\vp\in C_b^2 ({\Bbb R}^r),\ r\in 
{\Bbb N},\ \vp_0\in C_0^2({\Bbb R}),\ k\in K\right\}\, . 
\end{eqnarray}
In particular, $k\in C(D)\cap L^2(D)$, $k\in K$, implies that 
discontinuities of $(E,\pi)\ni\nu\to (k,\nu)$ can only occur when $|(k, 
\nu)|\to\infty$. Thus, $\t C_0^2(E)\subseteq C_b(E)$.  
\begin{lemma}\label{Lemma4.2} 
Assume (i)-(iii) and (j)-(jjj) of the present subsection. (a) Then the 
flow $U^+(t,\mu)$, $t\ge 0$, $\mu\in E\cap H$, is associated with a 
strongly continuous semigroup $(T_t)_{t\ge 0}$ on $L^2(E,\bnu)$ given by 
\begin{eqnarray*}
T_tf(\mu)=f(U^+(t,\mu))\, , \quad \mu\in E\cap H,\ t\ge 0,\ f\in L^2(E,\bnu). 
\end{eqnarray*}
(b) Let $A$ denote its generator. We have $\t C_b^2(E)\cup\t C_0^2(E)\subseteq 
D(A)$. If $f\in\t C_b^2(E)$ then with $f$ and $\vp$ related as in (\ref{4.3}), 
\begin{eqnarray*}
Af=\sum_{i=1}^r\frac{\partial\vp}{\partial x_i}\cdot (h_i,\cdot)\left( 
\lambda_i-z'(\cdot ,0)\right)\, .  
\end{eqnarray*}
If $f\in\t C_0^2(E)$ then with $f$ and $\vp$ as well as $\vp_0$ related as 
in (\ref{4.4}), 
\begin{eqnarray*}
Af&&\hspace{-.5cm}=\sum_{i=1}^r\frac{\partial\vp}{\partial x_i}\vp_0\cdot 
(h_i,\cdot)\left(\lambda_i-z'(\cdot ,0)\right)+\vp\vp'_0\cdot\sum_{j=1 
}^\infty (h_j,k)(h_j,\cdot)\left(\lambda_j-z'(\cdot ,0)\right) \\ 
&&\hspace{-.5cm}=\sum_{i=1}^r\frac{\partial\vp}{\partial x_i}\vp_0\cdot 
\lambda_i(h_i,\cdot)+\vp\vp'_0\cdot\left(k,{\T\frac12}\Delta\, \cdot\, 
\right)-z'(0)\cdot\left(\sum_{i=1}^r\frac{\partial\vp}{\partial x_i}\vp_0 
\cdot (h_i,\cdot)+\vp\vp'_0\cdot(k,\cdot)\right)\, . 
\end{eqnarray*}
(c) Both spaces, $\t C_b^2(E)$ as well as $\t C_0^2(E)$, are dense in the 
complete space $D(A)$ with respect to the graph norm $(\langle f\, ,\, f\rangle 
+\langle Af\, ,\, Af\rangle)^{1/2}$. 
\end{lemma} 
\begin{corollary}\label{Corollary4.3} 
Assume (i)-(iii) and (j)-(jjj) of the present subsection. \\ 
(a) We have $A'\1=-\delta (A^f)$. For all $f,g\in D(A)$ it holds that 
\begin{eqnarray}\label{4.23}  
\langle -Af,g\rangle+\langle -Ag,f\rangle=\left\langle\delta (A^f)\cdot f,g 
\right\rangle\, . 
\end{eqnarray} 
(b) We have $D(A)=D(A')$ and $A'f=-Af+A'\1\cdot f$, $f\in D(A)$. 
\end{corollary} 

\section{Asymptotic Properties of the Particle System} 
\setcounter{equation}{0}

For $i\in\{1,\ldots ,n\}$ and $z_1,\ldots ,z_{i-1},z_{i+1},\ldots ,z_n 
\in D$ let 
\begin{eqnarray*}
z\{i\}:=\left\{z_1,\ldots ,z_{i-1},z_{i+1},\ldots ,z_n\vphantom{l^1}\right\} 
\end{eqnarray*}
and
\begin{eqnarray*}
z(i):=\left\{(z_1,\ldots ,z_{i-1},y_i,z_{i+1},\ldots ,z_n):y_i\in\partial D 
\vphantom{l^1}\right\}\, . 
\end{eqnarray*}
Set $\partial^{(1)} D^n:=\bigcup_{z\in D^n}\bigcup_i\, z(i)$ and note that 
$\partial^{(1)} D^n=\partial D\times D\times\ldots\times D\cup\ldots\cup 
D\times\ldots\times D\times\partial D$. 

Let $B=((B_t)_{t\ge 0},(P_x)_{x\in {\Bbb R}^{n\cdot d}})$ be an ${n\cdot d} 
$-dimensional standard Brownian motion. Let $\tau$ denote the first exit 
time of $B$ from $D^n$ and let $B^{D^n}=((B_t^{D^n})_{t\ge 0},(P_x)_{x\in 
D^n})$ be obtained by stopping $B$ at time $\tau$. Let us consider a process 
$X\equiv X^n=((X_t)_{t\ge 0},(Q_x)_{x\in D^n})$ which is constructed as 
follows: $X_0$ is a random point on $(D^n,{\cal B}(D^n))$ with distribution 
$\t \bnu_n$. Inside the state space $D^n$, the process $X$ behaves like  
an $N\cdot d$-dimensional Brownian motion. Then, whenever the process $X$ 
hits the boundary $\partial^{(1)} D^n$, say at some point $y\in z(i)$ for 
some $i\in\{1,\ldots ,n\}$ and $z_1,\ldots ,z_{i-1},z_{i+1},\ldots ,z_n\in 
D$, it instantaneously jumps to $z=(z_1,\ldots ,z_n)$ where 
\begin{itemize}
\item[(k)] $z_i$ is a random point on $(D,{\cal B}(D))$ with distribution 
$\eta_{n,z\{i\}}(x)\, dx$. 
\begin{itemize}
\item[(k1)] For $\eta\equiv\eta_{n,z\{i\}}\in C^1_b(D)$ we assume  that there 
exists $c_1>1$ neither depending on $n\in {\Bbb N}$ nor on $z\{i\}$ such that 
\begin{eqnarray*}
c_1^{-1}h_1\le\eta_{n,z\{i\}}\le c_1h_1\, . 
\end{eqnarray*}
\item[(k2)] For $\bnu$-a.e. $\mu$, we assume $\eta_{n,z\{i\}}\, dx\stack{n\to 
\infty}{\Ra}\mu$ whenever $\frac1n\sum_{j=1}^n\delta_{z_j}\stack{n\to\infty} 
{\Ra}\mu$ where every $z_j$ may depend on $n$.  
\end{itemize}
\end{itemize}
We note that (k1) yields the following. For all $n\in {\Bbb N}$, 
\begin{eqnarray}\label{(k3)}
\inf_{z\{i\},x\in K}\eta_{n,z\{i\}}(x)>0\quad\mbox{\rm for every compact set}\ 
K\subset D 
\end{eqnarray}
and 
\begin{eqnarray}\label{(k4)} 
\sup_{{n,z\{i\}}}\|\eta_{n,z\{i\}}\|<\infty\, . 
\end{eqnarray}
Furthermore, we stress that $\eta_{n,z\{i\}}(x)$ is independent of $y\in z(i)$.  
\medskip 

We note that two independent $d$-dimensional Brownian particles can almost never 
hit $\partial D$ at the same time. This implies that $X$ $Q_{\t \sbnu_n}$-almost 
never reaches $\partial D^n\setminus\partial^{(1)}D^n$. This property allows to 
identify $\partial D^n$ with $\overline{\partial D\times D\times\ldots\times D} 
\cup\ldots\cup\overline{D\times\ldots\times D\times\partial D}$. 
Well-definiteness of the system can now be proved by using condition (k1) and 
the argument of \cite{BP09}, Lemma 1, applied to each of the compact components 
$\overline{\partial D\times D\times\ldots\times D}$, $\ldots\, $, $\overline{D 
\times\ldots\times D\times\partial D}$. 
\medskip 

Recalling that $\t \bnu_n$ is obtained from the 
measure $\bnu_n$ by the map $\frac1n\sum_{i=1}^n\delta_{z_i}\to(z_1,\ldots 
,z_n)$, $z_1,\ldots ,z_n\in\overline{D}$, taking into consideration 
invariance under permutations of $(z_1,\ldots ,z_n)$; for details see 
\cite{Lo09}, (sub)sections 2.4-2.6 and 6. For the measure $\t \bnu_n$ 
we shall assume the following. 
\begin{itemize}
\item[(kk)] $\t \bnu_n$ admits a density $\t m_n$ with respect to the 
Lebesgue measure on ${\Bbb R}^{n\cdot d}$ and $\t m_n$ is symmetric with 
respect to the $n$ $d$-dimensional components of ${\Bbb R}^{n\cdot d}$. 


$\t m_n$ is supported by $\overline{D^n}$ such that 
\begin{eqnarray*}
0<\t m_n\in C^1(D^n)\, ,\quad\nabla\t m_n\in C_b(D^n;{\Bbb R}^{n\cdot d}) 
\end{eqnarray*}
where $C_b(D^n;{\Bbb R}^{n\cdot d})$ denotes the space of all bounded 
continuous functions on $D^n$ with values in ${\Bbb R}^{n\cdot d}$. 
Furthermore, there exist positive constants $C_1,C_2$ which may depend on $n 
\in {\Bbb N}$ but are independent of $z\in D^n$ such that 
\begin{eqnarray*}
C_1\, h_1(z_1)\ldots h_1(z_n)\le\t m_n(z)\le C_2\, h_1(z_1)\ldots h_1(z_n)\, . 
\end{eqnarray*}
We note that the latter implies $\t m_n=0$ on $\partial D^n$ in the sense of 
$\lim_{D^n\ni x\to y}\t m_n(x)=0$, $y\in\partial D^n$. 
\item[(kkk)] Suppose that $-{\T\frac12}\Delta\t m_n$ exists in the sense that
\begin{eqnarray*} 
\int{\T\frac12}\Delta\vp\, d\t \bnu_n=-\int\vp\left(-{\T\frac12}\Delta\t m_n 
\right)\, dx\, ,\quad\vp \in \{\psi\in C^2_b(D^n):\psi(y)=0\, ,\ y\in\partial 
D^n\},   
\end{eqnarray*} 
and belongs to $L^2(D^n)$. 
\item[(kw)]Furthermore, let us assume that 
\begin{eqnarray*} 
\bnu_n\stack{n\to\infty}{\Ra}\bnu\, . 
\end{eqnarray*} 
\end{itemize} 

Denote by $\sigma$ the Lebesgue surface measure on $(\partial D^n,{\cal B} 
(\partial D^n))$ and by $s$ the Lebesgue surface measure on $(\partial D, 
{\cal B}(\partial D))$. In the following, let us arrange the notation 
according to $z=(z_1,\ldots ,z_n)$ where $z_1,\ldots ,z_n$ $\in D$, $y:= 
(z_1,\ldots ,z_{i-1},y_i,z_{i+1},\ldots ,z_n)$, $y_i\in\partial D$, $i\in 
\{1,\ldots ,n\}$. The transition probability function $Q$ of the process 
$X=(X_1,\ldots ,X_n)$ satisfies 
\begin{eqnarray}\label{4.24}
Q_x(X_t\in A)&=&P_x(B^{D^n}_t\in A)\nonumber \\ 
&&\hspace{-2cm}+\int_{z\in D^n}\int_{v=0}^t\sum_{i=1}^n\int_{y_i\in\partial D} 
\frac{dP_x(B_\tau\in dy,\, \tau\in dv)}{\sigma(dy)}\, s(dy_i)\cdot\eta (z_i) 
\, Q_z(X_{t-v}\in A)\, dz\, , \qquad\   
\end{eqnarray} 
$x\in D^n$, $t\ge 0$, $A\in {\cal B}(D^n)$. We mention that this 
representation involves an ordering among the $n$ particles. In order to 
represent this relation in terms of the empirical process $X^n_t:=\frac1n 
\sum_{i=1}^n\delta_{(X_t)_i}$, $t\ge 0$, for which the ordering is irrelevant, 
we refer to a similar situation in \cite{Lo09}, (sub)sections 2.4-2.6 and 6. 

For $f\in C_b(E)$ let $\t f\equiv\t f_n$ be defined by 
\begin{eqnarray*}
\t f(z_1,\ldots ,z_n):=f(\mu)\quad\mbox{\rm where}\quad\mu=\frac1n\sum_{i 
=1}^n\delta_{z_i}\in E_n\, . 
\end{eqnarray*} 
For $\beta\ge 0$ and let $\mu^{n,\beta}_x:=E_x\left(e^{-\beta\tau}\t 
1_{\cdot}(B_\tau)\right)$ and write $\mu^n_x$ instead of $\mu^{n,0}_x$, 
$x\in {D^n}$. For $z_1,\ldots ,z_n\in D$ and $z=(z_1,\ldots ,z_n)$ introduce 
\begin{eqnarray*}
m^{n,\beta}_x(dz):=\sum_{i=1}^n\int_{y_i\in\partial D}\frac{d\mu^{n, 
\beta}_x}{d\sigma}(z_1,\ldots ,z_{i-1},y_i,z_{i+1}, \ldots ,z_n)\, 
s(dy_i)\cdot\eta(z_i)\, dz\, ,\quad\beta\ge 0. 
\end{eqnarray*} 
This defines measures $m^{n,\beta}_x$ on $(D^n,{\cal B}(D^n))$, $x\in 
D^n$. We will also write $m^n_x$ instead of $m^{n,0}_x$, $x\in {D^n}$. 
\medskip 

\nid
{\bf Remarks} (1) For $u\in L^\infty(D^n)$ let us define the resolvent 
relative to the transition probabilities (\ref{4.24}), 
\begin{eqnarray*}
G^n_\beta u:=\int_{t=0}^\infty e^{-\beta t}\int u(y)\, Q_\cdot(X_t\in dy) 
\, dt\, ,\quad \beta>0. 
\end{eqnarray*}

Define $C_{b,n}(D^n):=C_b(D^n\cup\partial^{(1)} D^n)$ and let $C_r 
(D^n)$ denote the space of all $u=v|_{D^n}$ where $v\in C_{b,n}(D^n)$ 
and $v(y)=\int_{z_i\in D}v(z)\, \eta(z_i)\, dz_i$ if $y=(z_1,\ldots ,z_{ 
i-1},y_i,z_{i+1},\ldots ,z_n)\in\partial^{(1)}D^n$ and $z_1,\ldots ,z_{ 
i-1},z_{i+1},\ldots ,z_n\in D$, $y_i\in\partial D$, $i\in\{1,\ldots ,n\} 
$. As a straight forward adaption of (\ref{4.24}) and \cite{Lo09}, 
Proposition 1, one obtains 
\begin{eqnarray*}
\t f=G^n_\beta u\in C_r(D^n)\quad\mbox{\rm for all }u\in L^\infty(D^n) 
\end{eqnarray*}
and that for $y=(z_1,\ldots ,z_{i-1},y_i,z_{i+1},\ldots ,z_n)\in 
\partial^{(1)}D^n$ with $z_1,\ldots ,z_{i-1},z_{i+1},\ldots ,z_n\in 
D$ and $y_i\in\partial D$, $i\in\{1,\ldots ,n\}$, we have  
\begin{eqnarray}\label{4.25}
\t f(y)=\int_{z_i\in D}\t f(z)\, \eta(z_i)\, dz_i 
\end{eqnarray}
where $z_1,\ldots ,z_{i-1},z_{i+1}, \ldots ,z_n$ have been fixed in the 
notation $z=(z_1,\ldots ,z_n)$. This yields  
\begin{eqnarray}\label{4.26}
\int_{D^n}\t f\, dm^{n,\beta}_x=\int_{\partial D^n}\t f\, d\mu_x^{n, 
\beta}\, ,\quad\t f=G^n_\beta u,\ ,u\in L^\infty(D^n),\ \beta\ge 0.   
\end{eqnarray}
(2) Introduce the space $C_{r,c}(D^n):=\{u\in C_r(D^n):u=U|_{D^n},\, U\in 
C (\overline{D^n})\}$ and endow it with the $\sup$-norm. 

Let us keep (\ref{4.24})-(\ref{4.26}) in mind and follow \cite{Lo09}, 
Sections 3 and 4 word for word. Except for a slight modification of Step 
1 of the proof of Lemma 3 all arguments can be taken over to the present 
situation. One arrives at the following version of Theorem 2 in 
\cite{Lo09}. 

Let $\beta>0$ and $D(A^n):=\{G^n_\beta u:u\in C_{r,c}(D^n)\}$ and $A^n 
u:=\frac12\Delta u$, $u\in D(A^n)$. The operator $(A^n,D(A^n))$ is the 
generator of a strongly continuous contraction semigroup $(T_{n,t}^r)_{ 
t\ge 0}$ in $C_{r,c}(D^n)$. For $u\in C_{r,c}(D^n)$, we have $T_{n,t}^r 
u(x)=\int u(y)\, Q_x(X_t\in dy)$, $t\ge 0$, $x\in D^n$. Furthermore, 
$Q_x(X_t\in dy)$ is the unique representing probability measure on $(D^n, 
{\cal B}(D^n))$ of $T_{n,t}^r u(x)$, $u\in D(A^n)$, $t>0$, $x\in D^n$. 
As an immediate consequence we obtain (\ref{4.25}) and (\ref{4.26}) for 
$\t f\in C_{r,c}(D^n)$. 
\medskip 

\nid 
(3) For the proof of Lemma \ref{Lemma4.4} (a) below it is important to 
note that $D(A^n)$ is dense in $L^2(D^n,\t \bnu_n)$. According to Remark 
(2), for this it is sufficient to demonstrate that $C_{r,c}(D^n)$ is 
dense in $L^2(D^n,\t \bnu_n)$ and by (kk) that $C_{r,c}(D^n)$ is dense 
in $L^2(D^n,dz)$ where $dz$ symbolizes in this remark the Lebesgue 
measure on $(D^n,{\cal B}(D^n))$. If this was not true there would exist 
an $G\in L^2(D^n,\t \bnu_n)$ with $\int_{D^n} FG\, dz=0$ for all $F\in 
C_{r,c}(D^n)$. But $G\, dz\cdot\chi_{D^n}+0\cdot\chi_{\partial D^n}$ is 
a finite signed measure on $(\overline{D^n},{\cal B}(\overline{D^n}))$ 
which does not belong to the closed linear hull of ${\cal M}_{r,c}:=\{ 
(\eta_{n,z\{i\}}\, dx-\delta_{y_i})\cdot\chi_{z_1\times\ldots\times z_{i 
-1}\times\overline{D}\times z_{i+1}\times\ldots\times z_n} :i\in\{1, 
\ldots ,n\}$, $z_1,\ldots ,z_n\in D$, $y_i\in\partial D\}$ with respect 
to the convergence $\t \bomega_n\stack{n\to\infty}{\lra}\t \bomega$ if 
$\int f\, d\t \bomega_n\stack{n\to\infty}{\lra}\int f\, d\t \bomega$, 
$f\in C_{r,c}(D^n)$. On the other hand, $C_{r,c}(D^n)=\{F\in C(\overline 
{D^n}):\int F\, d\mu=0$ for all $\mu\in {\cal M}_{r,c}\}$, here $F\in 
C_{r,c}(D^n)$ identifying with its unique extension to $\overline{D^n}$. 
\medskip 

\nid 
(4) Consider a Markov chain on $D^n$ with transition probability kernel 
$m^n_\cdot\equiv m^{n,0}_\cdot$. It is an immediate consequence of 
condition (k) and \cite{DMT95}, Theorem 2.1, that the Markov chain is 
geometrically ergodic in the sense of this reference and that there 
exists a unique invariant probability measure. To verify this, define 
for $\delta >0$ the sets $D_\delta:=\{z\in D:|d-z|>\delta$ for all $d 
\in\partial D\}$, $C_\delta:=D_\delta^n$, and 
\begin{eqnarray*}
C_\delta^{(k)}&&\hspace{-.5cm}:=\left\{z=(z_1,\ldots ,z_n)\in D^n:\, 
\mbox{\rm there exist mutually distinct numbers}\right. \\ 
&&\hspace{.1cm}\left. i_1,\ldots , i_k\in \{1,\ldots ,n\}\ \mbox{such 
that}\ z_i\in D_\delta^c\ \mbox{\rm if}\ i\in\{i_1,\ldots ,i_k\}\ 
\mbox{\rm and}\right. \\ 
&&\hspace{.1cm}\left. z_i\in D_\delta\ \mbox{\rm if}\ i\not\in\{i_1, 
\ldots ,i_k\}\vphantom{D^n}\right\}\, ,\quad k\in\{1,\dots ,n\}.  
\end{eqnarray*}
It follows now from (\ref{(k3)}) that for all $\delta<\delta'$ for some 
$\delta'>0$, $C:=C_\delta$ is a small set in the sense of \cite{DMT95} 
since (4) in \cite{DMT95} is satisfied for the one-step transition 
probability kernel with $\mu$ being the equi-distribution on $\left(D^n, 
{\cal B}(D^n)\right)$. Furthermore, by the definition of the transition 
probability kernel $m^n_\cdot$ and property (\ref{(k4)}) there is a 
$\delta>0$ such that (7) of \cite{DMT95} is satisfied for $V=1$ on $C= 
C_\delta$ and, for example, for $V=\sum_{i=n-k}^{n-1}3^{ni}$ on 
$C_\delta^{(k)}$, $k\in\{1,\dots ,n\}$. To show the latter we note that 
the verification of (7) in \cite{DMT95} for $x\in C_\delta^{(k)}$ 
reduces to a combinatorial problem. One just has to check on which side 
of $\partial D_\delta$ the $n-k$ $d$-dimensional particles which start 
in $D_\delta$ are situated at the time of the next jump of $X$. 
\medskip 

\nid 
(5) Let ${\bf M}^n$ denote the invariant probability measure of the Markov 
chain on $D^n$ with transition probability kernel $m^n_\cdot\equiv m^{n,0 
}_\cdot$ considered in Remark (4). Define by  
\begin{eqnarray*}
{\bf N}^n(dx):=\int K^{D^n}(x,y)\, {\bf M}^n(dy)\, dx 
\end{eqnarray*}
a measure on $(D,{\cal B}(D^n))$, where $K^{D^n}$ denotes the Greenian 
kernel relative to $\frac12\Delta$ on $D^n$. Let us show that ${\bf N}^n$ 
is an invariant measure for the process $X$ by verifying $\int A^nu\, d{\bf 
N}^n=0$ for all $u\in D(A^n)$. We have for all $u\in D(A^n)$ 
\begin{eqnarray*} 
\int A^nu\, d{\bf N}^n&&\hspace{-.5cm}=\int{\T\frac12}\Delta u\, d{\bf N}^n 
 \\ 
&&\hspace{-.5cm}=\int_{x\in D^n}\left({\T\frac12}\Delta u\right)(x)\int_{y 
\in D^n} K^{D^n}(x,y)\, {\bf M}^n(dy)\, dx \\ 
&&\hspace{-.5cm}=\int_{y\in D^n}\int_{x\in D^n}\left({\T\frac12}\Delta u 
\right)(x)K^{D^n}(x,y)\, dx\, {\bf M}^n(dy) \\ 
&&\hspace{-.5cm}=\int (h_u-u)\, d{\bf M}^n, 
\end{eqnarray*}
$h_u$ denoting the harmonic function that satisfies $h_u=u$ on $\partial D^n$. 
Now we are going to use the fact that ${\bf M}^n$ is invariant with respect 
to $m^n_\cdot\equiv m^{n,0}_\cdot$. Furthermore, we will use (\ref{4.26}). We 
verify in this way for all $u\in D(A^n)$ 
\begin{eqnarray*} 
\int A^nu\, d{\bf N}^n&&\hspace{-.5cm}=\int h_u\, d{\bf M}^n-\int u\, d{\bf 
M}^n \\ 
&&\hspace{-.5cm}=\int h_u\, d{\bf M}^n-\int_{x\in D^n}\int_{z\in D^n} u(z)\, 
m^n_x (dz)\, {\bf M}^n(dx) \\ 
&&\hspace{-.5cm}=\int h_u\, d{\bf M}^n-\int_{x\in D^n}\int_{y\in\partial D^n} 
u(y)\, \mu^n_x(dy)\, {\bf M}^n(dx) \\ 
&&\hspace{-.5cm}=0\, . \vphantom{\int}
\end{eqnarray*}
\medskip

For $r\in {\Bbb N}$ and $f=\vp\left((h_1,\cdot),\ldots ,(h_r,\cdot ) 
\right)\in\t C_b^2(E)$, let us introduce  
\begin{eqnarray*}
{\textstyle\frac12}\Del f(\mu)&&\hspace{-.5cm}:=\sum_{i=1}^r\lambda_i 
\cdot\frac{\partial\vp}{\partial x_i}\cdot(h_i,\mu)+\frac{1}{2n}\sum_{ 
i,j=1}^r\frac{\partial^2\vp}{\partial x_i\partial x_j}\cdot (\nabla h_i 
\cdot\nabla h_j,\mu)\, , \quad \mu\in E_n, 
\end{eqnarray*}
where $\nabla h_i\cdot\nabla h_j$ is the scalar product of $\nabla h_i$ 
and $\nabla h_j$ in ${\Bbb R}^d$. Also, let us recall that, for $n\in 
{\Bbb N}$ and $\mu=\frac1n\sum_{k=1}^n\delta_{z_k}\in E_n$, we have the 
representation 
\begin{eqnarray}\label{4.27}
{\T\frac12}\Delta\vp\left(\frac{1}{n}\sum_{k=1}^nh_1(z_k),\ldots , 
\frac{1}{n}\sum_{k=1}^nh_r(z_k)\right)={\T\frac12}\Del f(\mu)\, , 
\end{eqnarray} 
cf. \cite{Lo09} and \cite{Lo13}. Recall also Lemma \ref{Lemma4.2} 
and for $f\in\t C_b^2(E)$ and $\vp$ related as in (\ref{4.3}) set 
\begin{eqnarray*}
Bf:=\sum_{i=1}^r\frac{\partial\vp}{\partial x_i}\cdot \lambda_i(h_i,\cdot) 
\quad\mbox{\rm and}\quad Cf:=-z'(0)\cdot\sum_{i=1}^r\frac{\partial\vp}{ 
\partial x_i}\cdot (h_i,\cdot)\, . 
\end{eqnarray*}
Similarly, if $f\in\t C_0^2(E)$ with $f$ and $\vp$ as well as $\vp_0$ are 
related as in (\ref{4.4}) set 
\begin{eqnarray*}
Bf&&\hspace{-.5cm}:=\sum_{i=1}^r\frac{\partial\vp}{\partial x_i}\vp_0\cdot 
\lambda_i(h_i,\cdot)+\vp\vp'_0\cdot\sum_{j=1}^\infty\lambda_j (h_j,k)(h_j, 
\cdot) \\ 
&&\hspace{-.5cm}=\sum_{i=1}^r\frac{\partial\vp}{\partial x_i}\vp_0\cdot 
\lambda_i(h_i,\cdot)+\vp\vp'_0\cdot\left(k,{\T\frac12}\Delta\, \cdot\, 
\right)
\end{eqnarray*}
and 
\begin{eqnarray*}
Cf&&\hspace{-.5cm}:=-z'(0)\cdot\left(\sum_{i=1}^r\frac{\partial\vp} 
{\partial x_i}\vp_0\cdot (h_i,\cdot)+\vp\vp'_0\cdot\sum_{j=1}^\infty(h_j, 
k)(h_j,\cdot)\right) \\ 
&&\hspace{-.5cm}=-z'(0)\cdot\left(\sum_{i=1}^r\frac{\partial\vp}{\partial 
x_i}\vp_0\cdot (h_i,\cdot)+\vp\vp'_0\cdot(k,\cdot)\right)\, . 
\end{eqnarray*}
We have $Af=Bf+Cf$ for all $f\in\t C_b^2(E)\cup\t C_0^2(E)$. Let us 
furthermore set 
\begin{eqnarray*} 
\t {\bf m}_n(dx):=-{\T\frac12}\Delta\t m_n(x)\, dx\, . 
\end{eqnarray*}

The following lemma is not just a collection of technicalities used in the 
paper. It is also the $L^2$-counterpart to Remark (2) of this section. In 
particular, part (b) of Lemma \ref{Lemma4.4} below is compatible with the 
domain $D(A^n)$ of the $C_{r,c}(D^n)$-infinitesimal operator relative to 
the transition probability function $Q_\cdot$. 
\begin{lemma}\label{Lemma4.4} 
Suppose (k)-(kw). (a) The process $X\equiv X^n=((X_t)_{t\ge 0},
(Q_x)_{x\in D^n})$ is associated with a strongly continuous semigroup on 
$L^2(D^n,\t \bnu_n)$. \\ 
(b) Let $f\in C_b(E_n)$ such that $\t f$ has a (unique) continuous extension 
to $D^n$ which belongs to $C_b^2(D^n)$. Then, using the notation $z=(z_1, 
\ldots ,z_n)$, $y\equiv (z_1,\ldots ,z_{i-1},y_i,z_{i+1},\ldots ,z_n)\in 
\partial D^n$ where $y_i\in\partial D$, $i\in\{1,\ldots ,n\}$, $z_1,\ldots 
,z_n\in D$, we have $f\in D(A_n)$ if and only if 
\begin{eqnarray*}
\int_{z_i\in D}\left(\t f(z)-\t f(y)\right)\eta(z_i)\, dz_i=0\quad\mbox{\rm 
for }\sigma\mbox{\rm -a.e. } y\in\partial D^n. 
\end{eqnarray*} 
In this case $A_nf=\frac12\Del f$. \\ 
(c) Let $f\in C_b(E_n)$ such that $\t f$ has a continuous extension to $D^n$ 
which belongs to $C_b^2(D^n)$ and let $g\in L^2(E,\bnu_n)$ such that $\t g\in 
C^2(D^n)$ and $\t g\cdot\t m_n=0$ on $\partial D^n$. We have 
\begin{eqnarray}\label{4.28}
&&\hspace{-.5cm}\lim_{t\downarrow 0}\frac1t\int\t g(x)\int\left(\t f(y)-\t 
f(x)\right)\, Q_x(X_t\in dy)\, \t \bnu_n(dx)\nonumber \\ 
&&\hspace{.5cm}=\left\langle{\T\frac12}\Delta\t f\, ,\, \t g\right\rangle_n 
+\int\t f\, d\int_{x\in D^n}(m^n_x-\mu^n_x)\left({\T -\frac12\Delta}(\t g\t 
m_n)\right)(x)\, dx\, .  
\end{eqnarray} 
(d) For $n\in {\Bbb N}$, let $D^2_b(A_n)$ denote the set of all $f\in D(A_n) 
\cap C_b(E_n)$ such that $\t f$ has a continuous extension to $D^n$ which 
belongs to $C_b^2(D^n)$. The set $D^2_b(A_n)$ is dense in the complete space 
$D(A_n)$ with respect to the graph norm $(\langle f\, ,\, f\rangle_n+\langle 
A_nf\, ,\, A_nf\rangle_n)^{1/2}$. 
\end{lemma} 
Proof. For the Brownian motion $\left((B_t)_{t\ge 0},(P_x)_{x\in {\Bbb R}^{ 
n\cdot d}}\right)$ on ${\Bbb R}^{n\cdot d}$ we will consider its components 
$\left((B_{i,t})_{t\ge 0},(P_{x_i})_{x_i\in {\Bbb R}^d}\right)$ on ${\Bbb R 
}^d$, $i\in \{1,\ldots ,n\}$. Let $\tau_i$ denote the first exit time of 
$B_{\cdot ,i}$ from $D$, $i\in\{1,\ldots ,n\}$.  
\medskip

\nid
{\it Step 1 } We verify (a). By (kk) we obtain for $f\in L^\infty(D^n, 
\t \bnu_n)$ with $f\ge 0$, 
\begin{eqnarray}\label{4.29}
&&\hspace{-.5cm}\left\|\int\t f(x)\, P_\cdot(B_t^{D^n}\in dx)\right\|_{L^1 
(D^n,\t \sbnu_n)}=\int_{x\in D^n}\t f(x)\int_{z\in D^n}P_z(B_t^{D^n}\in dx) 
\t m_n(z)\, dz\nonumber \\ 
&&\hspace{.5cm}\le C_2\int_{x\in D^n}\t f(x)\int_{z\in D^n}P_z(B_t^{D^n} 
\in dx)h_1(z_1)\ldots h_1(z_n)\, dz\nonumber \\ 
&&\hspace{.5cm}=C_2\int_{x\in D^n}\t f(x)e^{n\lambda_1t}h_1(x_1)\ldots h_1 
(x_n)\, dx\nonumber \\ 
&&\hspace{.5cm}\le\frac{C_2e^{n\lambda_1t}}{C_1}\int_{x\in D^n}\t f(x)\t m_n 
(x)\, dx\equiv\frac{C_2e^{n\lambda_1t}}{C_1}\int f\, d\bnu_n\, , \quad t>0.
\end{eqnarray} 

Next we use the fact that, for fixed $y\equiv (z_1,\ldots ,z_{i-1},y_i,,z_{ 
i+1}\ldots ,z_n)\in\partial D^n$, 
\begin{eqnarray}\label{4.30} 
\frac{P_{x_j}\left(\tau_j\ge t,\, B_{j,t}\in dz_j\right)}{dz_j}=\sum_{k=1 
}^\infty e^{\lambda_kt}h_k(x_j)h_k(z_j)
\end{eqnarray} 
as well as, with ${\bf n}$ being the inner normal vector on $\partial D$, 
\begin{eqnarray}\label{4.31} 
\frac{P_{x_i}\left(\tau_i\in dt,\, B_{i,t}\in dy_i\right)}{s(dy_i)}=\frac 
12\sum_{k=1}^\infty e^{\lambda_kt}h_k(x_i)\frac{\partial h_k}{\partial{\bf 
n}}(y_i)\, dt\, . 
\end{eqnarray} 
Let $H_k$, $k\in {\Bbb N}$, denote the eigenfunctions of the Dirichlet 
Laplacian on functions defined on $D^n$ and normed in $L^2(D^n)$ and let 
$0>2\lambda_{n,1}\ge2\lambda_{n,2}\ldots\, $ denote the corresponding 
eigenvalues. Using again the convention $y=(z_1,\ldots ,z_{i-1},y_i,z_{i 
+1},\ldots ,z_n)$, $z_1,\ldots z_n\in D$, $y_i\in\partial D$, $i\in\{1, 
\ldots ,n\}$, we get 
\begin{eqnarray}\label{4.32} 
&&\hspace{-.5cm}\int_{x\in D^n}\sum_{i=1}^n\int_{y_i\in\partial D}\frac 
{P_x(B_\tau\in dy,\, \tau\in dt)}{\sigma(dy)\times dt}\, s(dy_i)\cdot 
\eta_{n,z\{i\}}(z_i)\t m_n(x)\, dx\nonumber \\ 
&&\hspace{.5cm}\le C_2\int_{x\in D^n}\sum_{i=1}^n\int_{y_i\in\partial D} 
\frac{P_x(B_\tau\in dy,\, \tau\in dt)}{\sigma(dy)\times dt}\, s(dy_i) 
\cdot\eta(z_i)H_1(x)\, dx\nonumber \\ 
&&\hspace{.5cm}=C_2\int_{x\in D^n}\sum_{i=1}^n\prod_{j\neq i}\frac{P_{x_j} 
\left(\tau_j\ge t,\, B_{j,t}\in dz_j\right)}{dz_j}\times\nonumber \\ 
&&\hspace{3.5cm}\times\int_{y_i\in\partial D}\frac{P_{x_i}\left(\tau_i\in 
dt,\, B_{i,t}\in dy_i\right)}{s(dy_i)\times dt}\, s(dy_i)\cdot\eta(z_i)H_1 
(x)\, dx\nonumber \\ 
&&\hspace{.5cm}=C_2\int_{x\in D^n}\sum_{k=1}^\infty -\lambda_{n,k}e^{ 
\lambda_{n,k}t}H_k(x)\sum_{i=1}^n\left(\int_{z_i\in D}H_k(z)\, dz_i\cdot 
\eta(z_i)\right)H_1(x)\, dx\nonumber \\ 
&&\hspace{.5cm}=-C_2\lambda_{n,1}e^{\lambda_{n,1}t}\sum_{i=1}^n\int_{z_i 
\in D}H_1(z)\, dz_i\cdot\eta(z_i)\nonumber \\ 
&&\hspace{.5cm}\le -C_2\lambda_{n,1}e^{\lambda_{n,1}t}n\|h_1\|_{L^1(D)} 
\cdot c_1\, H_1(z)\vphantom{\sum_{i=1}^n}\nonumber \\ 
&&\hspace{.5cm}\le -\frac{nc_1\|h_1\|_{L^1(D)}C_2}{C_1}\lambda_{n,1}e^{ 
\lambda_{n,1}t}\cdot\t m_n (z)\vphantom{\sum_{i=1}^n}\, , \quad t>0, 
\end{eqnarray} 
where we have used (kk) for the first and the last $``\le"$ sign, 
(\ref{4.30}) as well as (\ref{4.31}) and Gauss' theorem for the second $``= 
"$ sign, and (k1) for the second last $``\le"$ sign. 
\medskip 

Let $\gamma:=-(nc_1\|h_1\|_{L^1(D)}C_2/C_1)\lambda_{n,1}$ and choose $f\in 
L^\infty(E,\bnu_n)$ with $f\ge 0$. From (\ref{4.24}), (\ref{4.32}), and 
Gronwall's inequality we obtain 
\begin{eqnarray*}
&&\hspace{-.5cm}\int\t f(x)\, Q_{\t \sbnu_n}(X_t\in dx)\le\int\t f(x)\, P_{ 
\t\sbnu_n}(B_t^{D^n}\in dx)+\gamma\int_{v=0}^t\int\t f(x)\, Q_{\t \sbnu_n}( 
X_v\in dx)\, dv \\ 
&&\hspace{.5cm}\le\int\t f(x)\, P_{\t\sbnu_n}(B_t^{D^n}\in dx)+\gamma\int_{ 
v=0}^te^{\gamma(t-v)}\int\t f(x)\, P_{\t\sbnu_n}(B_v^{D^n}\in dx)\, dv\, . 
\end{eqnarray*} 
By using the Schwarz' inequality we find with (\ref{4.29}) and monotone 
convergence that for $f\in L^2(E,\bnu_n)$
\begin{eqnarray}\label{4.33}
&&\hspace{-.5cm}\int\left(\int\t f(x)\, Q_\cdot (X_t\in dx)\right)^2\, d\t 
\bnu_n\le\int\t f^2(x)\, Q_{\t \sbnu_n}(X_t\in dx) \nonumber \\ 
&&\hspace{.5cm}\le\int\t f^2(x)\, P_{\t\sbnu_n}(B_t^{D^n}\in dx)+\gamma 
\int_{v=0}^te^{\gamma(t-v)}\int\t f^2(x)\, P_{\t\sbnu_n}(B_v^{D^n}\in dx)\, 
dv\nonumber \\ 
&&\hspace{.5cm}\le \frac{C_2e^{n\lambda_1t}}{C_1}\int f^2\, d\bnu_n+\gamma 
\frac{C_2}{C_1}\int_{v=0}^te^{\gamma(t-v)}\, e^{n\lambda_1v}\int f^2\, d 
\bnu_n\, dv\, . 
\end{eqnarray} 
This says that the process $X\equiv X^n=((X_t)_{t\ge 0},(Q_x)_{x\in D^n})$ 
is associated with a semigroup on $L^2(D^n,\t \bnu_n)$. Let $(T_{n,t})_{t 
\ge 0}$ denote the corresponding semigroup on $L^2(E,\bnu_n)$. It remains 
to demonstrate that $(T_{n,t})_{t\ge 0}$ is strongly continuous in $L^2(E, 
\bnu_n)$. As an immediate consequence of (\ref{4.33}) there is a constant 
$C_3>0$ such that for all $f\in L^2(E,\bnu_n)$ and all $t\in [0,1]$ 
\begin{eqnarray*}
\left\langle T_{n,t}f\, ,\, T_{n,t}f\right\rangle_n\le C_3\left\langle f,f 
\right\rangle_n\, . 
\end{eqnarray*} 
But this together with Remarks (2) (3) of this section guarantees strong 
continuity. The proof is just an adaption of \cite{Lo13}, proof of 
Proposition 4.5 (a), Step 1 from (4.26) on. 
\medskip 

\nid
{\it Step 2 } We prove (b) and (c). Denote by $(y)_i$ the $i$-th 
$d$-dimensional component of $y\equiv (z_1,\ldots ,z_{i-1},y_i,z_{i+1} 
,\ldots ,z_n)\in\partial D^n$, $i\in\{1,\ldots ,n\}$, $z_1,\ldots ,z_n 
\in D$. Introduce the notation 
\begin{eqnarray*} 
\int_i\int\equiv\int\limits_{z_1,\ldots ,z_{i-1}\in D,}\int\limits_{z_{ 
i+1},\ldots, z_n\in D}
\end{eqnarray*} 

Let us recall Remark (1) of this section, particularly $\t f=G^n_\beta u 
\in C_r(D^n)$ for all $u\in L^\infty(D^n)$. Recalling that for Lipschitz 
domains the Euklidean boundary is identifiable with the Martin boundary 
by the result of \cite{HW70}, for each such $\t f$ there is a harmonic 
function $h_{\t f}$ on $D^n$ with $h_{\t f}(z)\to h_{\t f}(y)$ if $D^n\ni 
z\to y\in\partial^{(1)}D^n$. This is motivation for the following. 

Let $f\in C_b(E_n)$ such that $\t f\in C_{b,n}(D^n)=C_b(D^n\cup 
\partial^{(1)} D^n)$. By (\ref{4.24}) we have for $x\in D^n$ 
\begin{eqnarray}\label{4.34}
&&\hspace{-.5cm}\lim_{t\downarrow 0}\frac1t\int\left(\t f(y)-\t f(x) 
\right)\, Q_x(X_t\in dy)\nonumber \\ 
&&\hspace{.5cm}=\lim_{t\downarrow 0}\frac1t\left(\int\left(\t f-h_{\t f} 
\right)(y)\, P_x(B^{D^n}_t\in dy)-\left(\t f-h_{\t f}\right)(x)\right) 
\nonumber \\ 
&&\hspace{1.0cm}+\lim_{t\downarrow 0}\frac1t\left(\vphantom{\sum_1^1} 
\left(\int\t f(y)\, P_x(B^{D^n}_t\in dy)-\int\t f(y)\, P_x(B_t\in dy) 
\right)\right. \nonumber \\ 
&&\hspace{1.0cm}+\sum_{i=1}^n\int_i\int\int_{z_i\in D}\t f(z)\left. 
\int_{y_i\in\partial D}\frac{P_x(B_\tau\in dy,\, \tau< t)}{\sigma(dy)}\, 
s(dy_i)\, \eta(z_i)\, dz\vphantom{\sum_1^1}\right) \nonumber \\ 
&&\hspace{.5cm}=\lim_{t\downarrow 0}\frac1t\left(\int\left(\t f-h_{\t f} 
\right)(y)\, P_x(B^{D^n}_t\in dy)-\left(\t f-h_{\t f}\right)(x)\right) 
\nonumber \\ 
&&\hspace{1.0cm}+\lim_{t\downarrow 0}\frac1t\sum_{i=1}^n\int_i\int 
\int_{z_i\in D}\int_{y_i\in\partial D}\left(\t f(z)-\t f(y)\right)\times 
\nonumber \\ 
&&\hspace{4.0cm}\times\frac{P_x(B_\tau\in dy,\, \tau <t)}{\sigma(dy)}\, 
s(dy_i)\, \eta(z_i)\, dz 
\end{eqnarray} 
and are interested in conditions on $f$ such that this limit exists in 
$L^2(D^n,\t \bnu_n)$. According to (\ref{4.30}) and (\ref{4.31}) the 
derivative 
\begin{eqnarray*}
&&\hspace{-.5cm}\lim_{t\downarrow 0}\frac1t\frac{P_x(B_\tau\in dy,\, \tau 
<t)}{\sigma(dy)} \\ 
&&\hspace{.5cm}=\lim_{t\downarrow 0}\frac1t\int_{v=0}^t\prod_{j\neq i} 
\frac{P_{x_j}\left(\tau_j\ge v,\, B_{j,v}\in dz_j\right)}{dz_j}\, \frac 
{P_{x_i}\left(\tau_i\in dv,\, B_{i,v}\in dy_i\right)}{s(dy_i)} \\ 
&&\hspace{.5cm}=\prod_{j\neq i}\sum_{k=1}^\infty h_k(x_j)h_k(z_j)\cdot 
\frac12\sum_{l=1}^\infty h_l(x_i)\frac{\partial h_l}{\partial {\bf n}}(y_i) 
\end{eqnarray*} 
exists in the distributional sense. We also remind of $\mu^n_x(dy)/\sigma 
(dy)=\frac12\sum_{m=1}^\infty$ $-\lambda_{n,m}^{-1}H_m(x)$ $\cdot\left( 
\partial H_m/\partial{\bf n}_n(y)\right)$ where ${\bf n}_n$ denotes the 
inner normal vector on $\partial D^n$. In fact, for $\t \vp\in C(D^n)$ with 
$\t \vp=0$ on $\partial D^n$ such that $-\frac12\Delta\t \vp$ exists in the 
sense of condition (kkk) and $\t \psi$ satisfying the conditions for $\t f$, 
we have 
\begin{eqnarray}\label{4.35}
&&\hspace{-.5cm}\sum_{i=1}^n\int_i\int\int_{z_i\in D}\int_{y_i\in\partial 
D}\int_{x\in D^n}\lim_{t\downarrow 0}\frac1t\frac{P_x(B_\tau\in dy,\, \tau 
<t)}{\sigma(dy)}\times\nonumber \\  
&&\hspace{4.0cm}\times\vphantom{\sum_{m=1}^\infty}\t \vp(x)\t \psi(z)\, 
dx\, s(dy_i)\, \eta(z_i)\, dz\nonumber \\ 
&&\hspace{.5cm}=\frac12\sum_{i=1}^n\int_i\int\int_{z_i\in D}\int_{y_i\in 
\partial D}\int_{x\in D^n}\sum_{m=1}^\infty H_m(x)\frac{\partial H_m} 
{\partial {\bf n}_n}(y)\times\nonumber \\ 
&&\hspace{4.0cm}\vphantom{\sum_{m=1}^\infty}\times\t \vp(x)\t \psi(z)\, dx 
\, s(dy_i)\, \eta(z_i)\, dz\nonumber \\ 
&&\hspace{.5cm}=\frac12\sum_{i=1}^n\sum_{m=1}^\infty\left({\T -\frac12} 
\Delta\t \vp,H_m\right)_{L^2(D^n)}\int_i\int\int_{y_i\in\partial D}- 
\lambda_{n,m}^{-1}\cdot\frac{\partial H_m}{\partial {\bf n}_n}(y)\, s(dy_i) 
\times\nonumber \\ 
&&\hspace{4.0cm}\vphantom{\sum_{m=1}^\infty}\times\int_{z_i\in D}\t \psi(z)
\, \eta(z_i)\, dz\nonumber \\ 
&&\hspace{.5cm}=\left({\T -\frac12}\Delta\t \vp\, ,\, \sum_{i=1}^n\int_i 
\int\int_{y_i\in\partial D}\frac12\sum_{m=1}^\infty -\lambda_{n,m}^{-1}\cdot 
\frac{\partial H_m}{\partial {\bf n}_n}(y)H_m(\cdot)\, s(dy_i)\right. 
\times\nonumber \\ 
&&\hspace{4.0cm}\left.\vphantom{\sum_{i=1}^n}\times\int_{z_i\in D}\t \psi(z) 
\, \eta(z_i)\, dz\right)_{L^2(D^n)} \nonumber \\ 
&&\hspace{.5cm}=\left({\T -\frac12}\Delta\t \vp\, ,\, \int_{z\in D^n}\t 
\psi(z)\, m^n_\cdot (dz)\right)_{L^2(D^n)}\, . 
\end{eqnarray} 
Well-definiteness for $\t \vp$ and especially $\t \psi$ as specified above 
follows from the third line together with property (\ref{(k4)}) and Gauss' 
formula. Similarly, one shows that 
\begin{eqnarray}\label{4.36}
&&\hspace{-.5cm}\sum_{i=1}^n\int_i\int\int_{z_i\in D}\int_{y_i\in\partial 
D}\int_{x\in D^n}\lim_{t\downarrow 0}\frac1t\frac{P_x(B_\tau\in dy,\, \tau 
<t)}{\sigma(dy)}\times\nonumber \\ 
&&\hspace{4.0cm}\times\vphantom{\sum_{m=1}^\infty}\t \vp(x)\t \psi(y)\, 
dx\, s(dy_i)\, \eta(z_i)\, dz\nonumber \\ 
&&\hspace{.5cm}=\left({\T -\frac12}\Delta\t \vp\, ,\, \int_{y\in\partial 
D^n}\t\psi(y)\, \mu^n_\cdot (dy)\right)_{L^2(D^n)}\, .  
\end{eqnarray} 
It follows now from (\ref{4.34}) that, for $f\in C_b(E_n)$ such that $\t f 
\in C_{b,n}(D^n)=C_b(D^n\cup\partial^{(1)} D^n)$, the limit 
\begin{eqnarray*} 
&&\hspace{-.5cm}\lim_{t\downarrow 0}\frac1t\int\left(\t f(y)-\t f(x) 
\right)\, Q_x(X_t\in dy) \\ 
&&\hspace{.5cm}=\lim_{t\downarrow 0}\frac1t\left(\int\left(\t f-h_{\t f 
}\right)(y)\, P_x(B^{D^n}_t\in dy)-\left(\t f-h_{\t f}\right)(x)\right) 
 \\ 
&&\hspace{1.0cm}+\int_{y\in\partial {D^n}}\sum_{i=1}^n\chi_{(y)_i\in 
\partial D}(y)\int_{z_i\in D}\left(\t f(z)-\t f(y)\right)\times \\ 
&&\hspace{3.0cm}\times\left.\frac{d^+}{dt}\right|_{t=0}\frac{P_x(B_\tau 
\in dy,\, \tau <t)}{\sigma(dy)}\, \eta(z_i)\, dz_i\, \sigma(dy) 
\end{eqnarray*} 
exists in $L^2(D^n,\t \bnu_n)$ if and only if 
\begin{eqnarray*} 
\int_{z_i\in D}\left(\t f(z)-\t f(y)\right)\eta(z_i)\, dz_i=0\quad\mbox 
{\rm for}\ \sigma\mbox{\rm-a.e.}\ y\in\partial D^n 
\end{eqnarray*} 
and the limit 
\begin{eqnarray*} 
\lim_{t\downarrow 0}\frac1t\left(\int\left(\t f-h_{\t f}\right)(y)\, 
P_x(B^{D^n}_t\in dy)-\left(\t f-h_{\t f}\right)(x)\right) 
\end{eqnarray*} 
exists in $L^2(D^n,\t \bnu_n)$. This gives part (b). Furthermore, 
(\ref{4.34})-(\ref{4.36}) imply also part (c). 
\medskip 

\nid
{\it Step 3 } We prepare the proof of part (d). Firstly, we remind of 
\cite{Lo14-2}, Step 2 of the proof of Lemma 4.5. In fact, the following 
is proved there. Let ${\cal A}$ be a densely defined closed operator in 
some Hilbert space ${\cal H}$ and let ${\cal C}\subseteq D({\cal A})$ 
be a set dense in ${\cal H}$ with $\{{\cal A}f:f\in {\cal C}\}\subseteq 
D({\cal A}')$. Then ${\cal C}$ is also dense in $D({\cal A})$ with 
respect to the graph norm $(\langle f\, ,\, f\rangle_{\cal H}+\langle 
{\cal A}f\, ,\, {\cal A}f\rangle_{\cal H})^{1/2}$. In Step 4 below we 
will apply this to ${\cal A}=A_n$ and ${\cal H}=L^2(E,\bnu_n)$. 
 
As a consequence of part (a) of the present lemma, there is an $M 
\ge 1$ and an $\omega\ge 0$ such that for the semigroup $(\t T_{n, 
t})_{t\ge 0}$ in on $L^2(D^n,\t \bnu_n)$ with generator $\t A_n$, 
and associated with $X\equiv X^n=((X_t)_{t\ge 0},(Q_x)_{x\in D^n} 
)$, we have $\|\t T_{n,t}f\|_{L^2(D^n,\t \bnu_n)}\le M\cdot e^{ 
\omega t}\|f\|_{L^2(D^n,\t \bnu_n)}$ for all $f\in L^2(D^n,\t  
\bnu_n)$. Below the notation $L_s$ will indicate that we are just 
count in the functions $w\equiv w(z,\ldots ,z_n)$ which are 
symmetric in $z_1,\ldots ,z_n\in D$ in the respective $L$-spaces. 
We shall, secondly, demonstrate that for $\beta>\omega$ the set 
$\{f:\t f=G^n_\beta u\in C_r(D^n),\ u\in L^\infty_s(D^n)\}$ is 
dense in $D(A_n)$ with respect to the graph norm $(\langle f\, ,\, 
f\rangle_n+\langle A_nf\, ,\, A_nf\rangle_n)^{1/2}$. 

For this assume that for a moment that there is a $v\in L^2_s(D^n, 
\t \bnu_n)\setminus\{0\}$ such that 
\begin{eqnarray*} 
\left(1+\beta^2\right)(G^n_\beta)' G^n_\beta v-\beta(G^n_\beta)'v=
\beta G^n_\beta v-v\, . 
\end{eqnarray*} 
This assumption would yield 
\begin{eqnarray*} 
\left\langle \left(1+\beta^2\right)G^n_\beta v-\beta v\, ,\, \t f
\right\rangle_n-\left\langle\beta G^n_\beta v-v\, ,\, \beta\t f
-\widetilde{A_n f}\right\rangle_n=0\, ,\quad f\in D(A_n). 
\end{eqnarray*} 
Thus $\langle G^n_\beta v\, ,\, \t f\rangle_n+\langle\t A_nG^n_\beta 
v\, ,\, \widetilde{A_n f}\rangle_n=0$ for all $f\in D(A_n)$ which 
would imply $v=0$. Now we can suppose that, for $v\in L_s^2(D^n,\t 
\bnu_n)\setminus\{0\}$, we have $\left(1+\beta^2\right)(G^n_\beta)' 
G^n_\beta v-\beta(G^n_\beta)'v-\beta G^n_\beta v+v\neq 0$. It  
follows  that 
\begin{eqnarray*}
\left\langle \left(1+\beta^2\right)(G^n_\beta)'G^n_\beta v-\beta 
(G^n_\beta)'v-\beta G^n_\beta v+v\, ,\, u\right\rangle_n=0  
\end{eqnarray*} 
cannot hold simultaneously for all $u\in L_s^\infty(D^n)$. In other 
words, $\langle\t A_nG^n_\beta v\, ,\, \t A_nG^n_\beta u\rangle_n 
+\langle G^n_\beta v\, ,\, G^n_\beta u\rangle_n=\langle\beta 
G^n_\beta v-v\, ,\, \beta G^n_\beta u-u\rangle_n+\langle G^n_\beta 
v\, ,\, G^n_\beta u\rangle_n=0$ cannot be true simultaneously for 
all $u\in L_s^\infty(D^n)$. This says that $\{f:\t f=G^n_\beta u 
\in C_r(D^n),\ u\in L_s^\infty(D^n)\}$ is dense in $D(A_n)$ with 
respect to the graph norm. 
\medskip 

\nid
{\it Step 4 } We prove part (d). Denote by $p^{D^n}(t,x,y)$ the 
transition density function relative to $P_x(B^{D^n}_t\in dy)$ and 
recall that it is symmetric in $x$ and $y$ for all $t>0$. 

Let $f\in \{\vp:\t \vp=G^n_\beta u\in C_r(D^n),\ u\in L_s^\infty(D^n) 
\}$ and $\t g\in C^2(D^n)$ such that $\t g=0$ on $\partial D^n$ and 
$\nabla\t g\in C_b(D^n;{\Bbb R}^{n\cdot d})$. It follows from (kk) 
that $\t g\t m_n=0$ as well as $\frac{\partial}{\partial{\bf n}_n} 
(\t g\t m_n)=0$ on $\partial D^n$. We have 
\begin{eqnarray*} 
&&\hspace{-.5cm}\lim_{t\downarrow 0}\frac1t\langle f\, ,\, T_{n,t}'g 
-g\rangle_n=\lim_{t\downarrow 0}\frac1t\langle T_{n,t}f-f\, ,\, g 
\rangle_n\vphantom{\left(\int\right)} \\ 
&&\hspace{.5cm}=\lim_{t\downarrow 0}\frac1t\left(\int\int\left(\t f- 
h_{\t f}\right)(y)\, p^{D^n}(t,x,y)\t g(x)\t m_n(x)\, dx\, dy\right. 
 \\ 
&&\hspace{2.0cm}\left.-\int\left(\t f-h_{\t f}\right)(x)\t g(x)\t 
m_n(x)\, dx\right) \\ 
&&\hspace{.5cm}=\left(\t f\, ,\, \Delta\left(\t g\t m_n\right)\right 
)_{L^2(D^n)}\vphantom{\left(\int\right)} \\ 
&&\hspace{1.0cm}-\lim_{t\downarrow 0}\frac1t\left(\int h_{\t f}(y)\, 
\int p^{D^n}(t,x,y)\t g(x)\t m_n(x)\, dx\, dy-\int h_{\t f}\, \t g\t 
m_n\, dx\right) \\ 
&&\hspace{.5cm}=\left(\t f\, ,\, {\T \frac12}\Delta\left(\t g\t m_n 
\right)\right)_{L^2(D^n)} \vphantom{\left(\int\right)}
\end{eqnarray*} 
where, in the last line, we have applied the Green formula. Taking 
into consideration that $\{f:\t f=G^n_\beta u\in C_r(D^n),\ u\in 
L_s^\infty(D^n)\}$ is dense in $D(A_n)$ with respect to the graph norm 
$(\langle f\, ,\, f\rangle_n+\langle A_nf\, ,\, A_nf\rangle_n)^{1/2}$ 
according to the second part of Step 3, we obtain 
\begin{eqnarray*} 
\left\{g:\t g\in C^2(D^n),\ \t g=0\ \mbox{\rm on}\ \partial D^n\ \mbox 
{\rm and}\ \nabla\t g\in C_b(D^n;{\Bbb R}^{n\cdot d})\right\}\subseteq 
D(A_n')\, . 
\end{eqnarray*} 

According to the first part of Step 3, it remains to show that 
\begin{eqnarray*} 
{\cal C}:=\left\{f\in D(A_n):\widetilde{A_nf}\in C^2(D^n),\ \widetilde 
{A_nf}=0\ \mbox{\rm on}\ \partial D^n\ \mbox {\rm and}\ \nabla 
\widetilde{A_nf}\in C_b(D^n;{\Bbb R}^{n\cdot d})\right\} 
\end{eqnarray*} 
is dense in $L^2(E,\bnu_n)$. We observe that the set 
\begin{eqnarray*} 
{\cal D}:=\left\{f:\Delta\t f\in C^2(D^n),\ \Delta\t f=0\ \mbox{\rm 
on}\ \partial D^n\ \mbox {\rm and}\ \nabla\Delta\t f\in C_b(D^n;{\Bbb 
R}^{n\cdot d})\right\} 
\end{eqnarray*} 
is dense in $L^2(E,\bnu_n)$. Among other things, this implies $\{\t 
f:f\in {\cal D}\}\subseteq C(\overline{D^n})$ in the sense that for 
$f\in {\cal D}$ and $y\in\partial D^n$ and $D^n\ni z\to y$ we have 
$\t f(z)\to\t f(y)$. By the result of Step 2 it just remains to show 
that 
\begin{eqnarray*} 
\hat{\cal C}:=\left\{f\in {\cal D}:\int_{z_i\in D}\left(\t f(z)-\t f 
(y)\right)\eta(z_i)\, dz_i=0\quad\mbox{\rm for}\ \sigma\mbox{\rm-a.e.} 
\ y\in\partial D^n \right\}
\end{eqnarray*} 
is dense in $L^2(E,\bnu_n)$. Any finite signed measure $\t \bomega$ 
on $\left(\overline{D^n},{\cal B}\left(\overline{D^n}\right)\right)$ 
that satisfies 
\begin{eqnarray*} 
\int\t f\, d\t \bomega=0\quad\mbox{\rm for all}\ f\in\hat{\cal C} 
\end{eqnarray*} 
belongs to the closed linear hull of ${\cal M}_{r,c}=\{(\eta_{n,z\{i\}} 
\, dx-\delta_{y_i})\cdot\chi_{z_1\times\ldots\times z_{i-1}\times 
\overline{D}\times z_{i+1}\times\ldots\times z_n} :i\in\{1,\ldots ,n\}, 
\ z_1,\ldots ,z_n\in D,\ y_i\in\partial D\}$ with respect to the 
convergence $\t \bomega_n\stack{n\to\infty}{\lra}\t \bomega$ if $\int\t 
\vp\, d\t \bomega_n\stack{n\to\infty}{\lra}\int\t \vp\, d\t \bomega$, 
$\vp\in {\cal D}$. Since $\t \bnu_n$ has no mass on $\partial D^n$ there 
is no $g\in L^2(E,\bnu_n)\setminus\{0\}$ such that 
\begin{eqnarray*} 
\int fg\, d\bnu_n=\int\t f\t g\, d\t \bnu_n=0\quad\mbox{\rm for all}\ 
f\in\hat{\cal C}\, .  
\end{eqnarray*} 
\qed
\medskip 

\begin{proposition}\label{Proposition4.5} 
Suppose (k)-(kw). \\ 
(a) For $f\in C_b(E)$ we have 
\begin{eqnarray*} 
\frac1n\int_{x\in D^n}\int_{y\in\partial D^n}\t f(y)\, \mu^n_x(dy)\, \t {\bf 
m}_n(dx)\stack{n\to\infty}{\lra}-\int f\, z'(\, \cdot\, ,0)\, d\bnu\, .  
\end{eqnarray*}  
For $f\in C_b(E)$ such that, for all $n\in {\Bbb N}$, it holds that $\t f\equiv 
\widetilde{f|_{E_n}}\in C(\overline{D^n})$ as well as $\t f=0$ on $\partial D^n$ 
and, in particular, for $f\in\t C_0^2(E)$ we have 
\begin{eqnarray}\label{4.37}  
\int f\, z'(\, \cdot\, ,0)\, d\bnu=0\, . 
\end{eqnarray}  
(b) For $f\in\t C_b^2(E)\cup\t C_0^2(E)$ we have 
\begin{eqnarray*}
\int\t f\, d\int_{x\in D^n}(m^n_x-\mu^n_x)\, \t {\bf m}_n(dx)\stack{n\to\infty} 
{\lra}\langle Cf\, ,\, \1\rangle\, , 
\end{eqnarray*} 
(c) For $f\in\t C_b^2(E)\cup\t C_0^2(E)$ we have 
\begin{eqnarray*}
\int\t f\, d\left(\int_{x\in D^n}\, \mu^n_x\, \t {\bf m}_n(dx)-\t {\bf m}_n 
\right)\stack{n\to\infty}{\lra}\langle Bf\, , \, \1\rangle\, . 
\end{eqnarray*} 
\end{proposition} 
Proof. {\it Step 1 } We verify (a). For $y\in\partial D^n$ we keep using 
the notation $y=(z_1,\ldots ,z_{i-1},y_i,$ $z_{i+1},\ldots ,z_n)$ where 
$y_i\in\partial D$, $z_1,\ldots ,z_n\in D$, $i\in\{1,\ldots ,n\}$. In 
addition, for $x\in D^n$ we use the representation $x=(x_1,\ldots ,x_n)$ 
with $x_1,\ldots ,x_n\in D$. As in Lemma \ref{Lemma4.4}, for the Brownian 
motion $((B_t)_{t\ge 0},$ $(P_x)_{x\in {\Bbb R}^{n\cdot d}})$ on ${\Bbb 
R}^{n\cdot d}$ we will also consider its components $\left((B_{i,t})_{t\ge 
0},(P_{x_i})_{x_i\in {\Bbb R}^d}\right)$ on ${\Bbb R}^d$, $i\in \{1,\ldots 
,n\}$. Let $\tau_i$ denote the first exit time of $B_{\cdot ,i}$ from $D$, 
$i\in\{1,\ldots ,n\}$.  Let us use the notation 
\begin{eqnarray*} 
\int_{i,y_i}\hspace{-0.1cm}\int\equiv\int\limits_{z_1,\ldots ,z_{i-1} 
\in D,}\int\limits_{y_i\in\partial D,}\int\limits_{z_{i+1},\ldots, z_n 
\in D}
\end{eqnarray*} 
It holds that
\begin{eqnarray}\label{4.38} 
&&\hspace{-.5cm}\lim_{n\to\infty}\int_{x\in D^n}\frac1n\int_{y\in\partial 
D^n}\t f(y)\, \mu^n_x(dy)\, \t {\bf m}_n(dx)\nonumber \\ 
&&\hspace{.5cm}=\lim_{n\to\infty}\int_{x\in D^n}\frac1n\sum_{i=1}^n\int_{i, 
y_i}\hspace{-0.1cm}\int\t f(y)\int_{t=0}^\infty\prod_{j\neq i}P_{x_j}\left( 
\tau_j>t,\, B_{j,t}\in dz_j\right)\times\nonumber \\ 
&&\hspace{1.0cm}\times P_{x_i}\left(\tau_i\in dt,\, B_{i,t}\in dy_i\right)\, 
\t {\bf m}_n(dx)\vphantom{\sum_{i=1}^n}\nonumber \\ 
&&\hspace{.5cm}=\lim_{n\to\infty}\int_{x\in D^n}\frac1n\sum_{i=1}^n\int_{i,y_i} 
\hspace{-0.1cm}\int f\left({\T\frac1n\sum_{j\neq i}\delta_{z_j} 
+\frac1n\delta_{y_i}}\right)\times\nonumber \\ 
&&\hspace{1.0cm}\times\int_{t=0}^\ve\prod_{j\neq i}P_{x_j}\left(\tau_j>t,\, 
B_{j,t}\in dz_j\right)P_{x_i}\left(\tau_i\in dt,\, B_{i,t}\in dy_i\right)\, 
\t {\bf m}_n(dx)\vphantom{\sum_{i=1}^n}\, . \nonumber \\ 
\end{eqnarray}  
for any $\ve>0$ and all $f\in C_b(E)$. We recall that, by definition, $E$ is 
compact. Thus, $f$ is uniformly continuous on $E$. Thus 
\begin{eqnarray*}
\lim_{n\to\infty}\sup_{z_1,\ldots ,z_n\in D,\, r\in\partial D}\left|f 
\left({\T\frac1n\sum_{j=1}^{n-1}\delta_{z_j}+\frac1n\delta_r}\right)-f 
\left({\T\frac1n\sum_{j=1}^n\delta_{z_j}}\right)\right|=0\, . 
\end{eqnarray*} 
The last two relations  imply  
\begin{eqnarray}\label{4.39} 
&&\hspace{-.5cm}\lim_{n\to\infty}\int_{x\in D^n}\frac1n\int_{y\in\partial 
D^n}\t f(y)\, \mu^n_x(dy)\, \t {\bf m}_n(dx)\nonumber \\ 
&&\hspace{.5cm}=\lim_{\ve\to 0}\lim_{n\to\infty}\int_{x\in D^n}\frac1n 
\sum_{i=1}^n\int_{z\in D^n}f\left({\T\frac1n\sum_{j=1}^n\delta_{z_j}} 
\right)\times\nonumber \\ 
&&\hspace{1.0cm}\times\int_{t=0}^\ve\prod_{j\neq i}P_{x_j}\left(\tau_j 
>t,\, B_{j,t}\in dz_j\right)P_{x_i}\left(\tau_i>t,\, B_{i,t}\in dz_i 
\right)\times\nonumber \\ 
&&\hspace{1.0cm}\times\int_{y_i\in\partial D}P_{x_i}\left(\tau_i\in dt, 
\, B_{i,t}\in dy_i\right)\, \t {\bf m}_n(dx)\vphantom{\sum_{i=1}^n} 
\nonumber \\ 
&&\hspace{.5cm}=\lim_{\ve\to 0}\lim_{n\to\infty}\int_{x\in D^n}\int_{z 
\in D^n}f\left({\T\frac1n\sum_{j=1}^n\delta_{z_j}}\right)\times 
\nonumber \\ 
&&\hspace{1.0cm}\times\int_{t=0}^\ve\prod_{j=1}^nP_{x_j}\left(\tau_j>t, 
\, B_{j,t}\in dz_j\right)\frac1n\sum_{i=1}^nP_{x_i}\left(\tau_i\in dt 
\right)\, \t {\bf m}_n(dx)\vphantom{\sum_{i=1}^n}\, . \nonumber \\ 
\end{eqnarray} 
Now we recall (\ref{2.4}) and (\ref{2.5}) as well as  $\bnu(H)=1$ (cf. 
condition (i) and $H(1)\subset H$), $\t \bnu_n(dz)=\int K^{D^n}(x,z)\, 
\t {\bf m}_n(dx)\, dz$ (cf. condition (kkk)), and $\bnu_n\stack{n\to 
\infty}\Ra\bnu$ (cf. condition (kw)). Let $\mu\in E\cap H$ and $\frac1n 
\sum_{j=1}^n\delta_{z_j}\stack{n\to\infty}{\Ra}\mu$. For the next 
calculations we give particular attention to 
\begin{eqnarray}\label{4.40} 
\lim_{t\downarrow 0}\int_{v\in D} g(v)\frac{P_v(\tau_i\le t)}{t}\, dv= 
-\int{\frac12}(\Delta g)(v)\, dv\, , \quad g\, dv=\mu\in E\cap H,  
\end{eqnarray} 
cf. (\ref{2.5}), which yields  
\begin{eqnarray*} 
\left.\left(\lim_{n\to\infty}\frac1n\sum_{i=1}^nP_{z_i}\left(\tau_i\in 
dt\right)\right)\right|_{t=0}=\left.P_{\mu}\left(\tau_{i_0}\in dt\right) 
\right|_{t=0}=-z'(\mu,0)\, dt\, ,  
\end{eqnarray*} 
for some fixed $i_0\in\{1,\ldots ,n\}$. Furthermore, 
\begin{eqnarray}\label{4.41} 
\left(K^{D^n}(x,z)\, \t {\bf m}_n(dx)\, dz-\int_{t=0}^\ve\prod_{j=1}^n 
P_{x_j}\left(\tau_j>t,\, B_{j,t}\in dz_j\right)\, \t {\bf m}_n(dx)\right) 
\stack{n\to\infty}\Ra 0 \hspace{-.5cm}\nonumber \\ 
\end{eqnarray} 
for any $\ve>0$, both items seen as mappings from measures over empirical 
probability measures to measures over $D^n\times D^n$. 
From (\ref{4.39})-(\ref{4.41}) and stochastic continuity of all $\left( 
(B_{i,t})_{t\ge 0},(P_{x_i})_{x_i\in {\Bbb R}^d}\right)$ we obtain 
\begin{eqnarray*} 
&&\hspace{-.5cm}\lim_{n\to\infty}\int_{x\in D^n}\frac1n\int_{y\in\partial 
D^n}\t f(y)\, \mu^n_x(dy)\, \t {\bf m}_n(dx) \\
&&\hspace{.5cm}=\lim_{\ve\to 0}\lim_{n\to\infty}\int_{x\in D^n}\int_{z\in 
D^n}f\left({\T\frac1n\sum_{j=1}^n\delta_{z_j}}\right)\times\nonumber \\ 
&&\hspace{1.0cm}\times\int_{t=0}^\ve\prod_{j=1}^nP_{x_j}\left(\tau_j>t,\, 
B_{j,t}\in dz_j\right)\frac1n\sum_{i=1}^nP_{z_i}\left(\tau_i\in dt\right) 
\, \t {\bf m}_n(dx)\vphantom{\sum_{i=1}^n} \\ 
&&\hspace{.5cm}=-\int f\, z'(\, \cdot\, ,0)\, d\bnu\, . 
\end{eqnarray*}  

For the proof of (b) and (c) assume $f\in\t C_b^2(E)$ with $f=\vp((h_1, 
\cdot),(h_2,\cdot),\ldots\, ,$ $(h_r,\cdot))$. The case $f\in\t C_0^2(E)$ 
can be handled in a similar fashion. 
\medskip 

\nid
{\it Step 2 } Let us prepare the proof of (b) and (c). Let $i\in\{1,\ldots , 
n\}$ and $y\equiv y^n\in\partial^{(1)}D^n$ such that $y=(z_1,\ldots ,z_{i-1}, 
y_i,z_{i+1},$ $\ldots ,z_n)$ with $z_1,\ldots ,z_n\in D$, $y_i\in\partial D$. 
In addition for $\mu=\frac1n\sum_{j=1}^n\delta_{z_j}$, let $\mu^{i,y_i}\equiv 
(\mu^{i,y_i})^n:=\frac1n\sum_{j\neq i}^n\delta_{z_j}+\frac1n\delta_{y_i}$, $i 
\in\{1,\ldots ,n\}$. Moreover, set  
\begin{eqnarray*}
r^n_{i,k}(y)&&\hspace{-.5cm}:=\int_0^1\left(\frac{\partial}{\partial x_k}\vp 
\left(\vphantom{\dot{f}}(h_1,\mu^{i,y_i})+{\textstyle\frac1n}(h_1(z_i)-h_1( 
y_i)),\ldots ,\right.\right. \\ 
&&\hspace{-1cm}\left.\left.(h_{k-1},\mu^{i,y_i})+{\textstyle\frac1n}(h_{k-1}( 
z_i)-h_{k-1}(y_i)),(h_k,\mu^{i,y_i})+{\textstyle\frac{t}{n}}(h_k(z_i)-h_k(y_i 
)),\vphantom{\int}\right.\right. \\ 
&&\hspace{-1cm}\left.\left.(h_{k+1},\mu^{i,y_i}),\ldots ,(h_r,\mu^{i,y_i}) 
\vphantom{\dot{f}}\right)-\frac{\partial}{\partial x_k}\vp(\ldots ,(h_k,\mu^{ 
i,y_i}),\ldots )\right)\, dt\, ,  
\end{eqnarray*}
$k\in\{1,\ldots ,r\}$. We have 
\begin{eqnarray*}
&&\hspace{-.5cm}n\int_{z_i\in D}\left(\t f(z)-\t f(y)\right)\eta(z_i)\, 
dz_i\\ 
&&\hspace{.5cm}=n\int_{z_i\in D}\left(\vp\left(\ldots ,(h_k,\mu^{i,y_i}) 
+{\textstyle\frac{1}{n}}(h_k(z_i)-h_k(y_i)),\ldots\right)\vphantom{\t f} 
\right. \\ 
&&\hspace{5.0cm}\left.-\vp(\ldots ,(h_k,\mu^{i,y_i}),\ldots )\vphantom 
{\t f}\right)\eta(z_i)\, dz_i \\ 
&&\hspace{.5cm}=\sum_{k=1}^r\left(\frac{\partial}{\partial x_k}\vp(\ldots 
,(h_k,\mu^{i,y_i}),\ldots )\right)\cdot\int_{z_i\in D}(h_k(z_i)-h_k(y_i)) 
\, \eta(z_i)\, dz_i \\ 
&&\hspace{1.0cm}+\sum_{k=1}^r\int_{z_i\in D}r^n_{i,k}(y)\, (h_k(z_i)-h_k 
(y_i))\, \eta(z_i)\, dz_i\, .   
\end{eqnarray*}
Since $\max_{j\in\{1,\ldots ,n\}\setminus\{i\}}\left|r^n_{j,k}(y)\right| 
\stack{n\to\infty}{\lra}0$ for any sequence $y\equiv y^n\in\partial^{(1)} 
D^n$ with $(\mu^{i,y_i})^n\Ra\mu$ and every $k\in\{1,\ldots ,r\}$, we 
obtain in this case 
\begin{eqnarray}\label{4.42}
n\int_{z_i\in D}\left(\t f(z)-\t f(y)\right)\, \eta(z_i)\, dz_i\stack{n\to 
\infty}{\lra}\sum_{i=1}^r\frac{\partial\vp}{\partial u_i}\cdot (h_i,\mu)  
\end{eqnarray} 
uniformly bounded with respect to $\bnu$-a.e. $\mu$, recall also condition 
(k2). The following is a preparation for the proof of Lemma \ref{Lemma4.12} 
below. Let $a$ be a real measurable bounded function on $E$ which is 
continuous on each $E_n$ such that for all $\frac1n\sum_{i=1}^n\delta_{z_i 
}\Ra\mu$ and $\frac1n\sum_{i=1}^n\delta_{z_i'}\Ra\mu$ we have $\left|a\left 
(\frac1n\sum_{i=1}^n\delta_{z_i}\right)-a\left(\frac1n\sum_{i=1}^n\delta_{ 
z_i'}\right)\right|\stack{n\to\infty}{\lra}0$. Assume furthermore 
$\sum_{i=1}^r\frac{\partial\vp}{\partial u_i}\cdot (h_i,\cdot)=0$ $\bnu 
$-a.e. Then replacing in the calculations of the present step $\eta(z_i)$ 
by $\eta(z_i)\cdot\t a(z)$ it turns out that 
\begin{eqnarray*}
n\int_{z_i\in D}\left(\t f(z)-\t f(y)\right)\, \eta(z_i)\cdot\t a(z)\, d 
z_i\stack{n\to\infty}{\lra}0   
\end{eqnarray*} 
uniformly bounded with respect to $\bnu$-a.e. $\mu$ if as above $(\mu^{i, 
y_i})^n\Ra\mu$ as $n\to\infty$. 
\medskip 

\nid 
{\it Step 3 } We prove (b). Relations (\ref{4.42}) and part (a) imply 
\begin{eqnarray*}
&&\hspace{-.5cm}\int\t f\, d\int_{x\in D^n}(m^n_x-\mu^n_x)\, \t {\bf m 
}_n(dx) \\ 
&&\hspace{.5cm}=\int_{x\in {D^n}}\int_{y\in\partial {D^n}}\sum_{i=1}^n 
\chi_{(y)_i\in\partial D}(y)\int_{z_i\in D}\left(\t f(z)-\t f(y)\right) 
\eta(z_i)\, dz_i\, \mu^n_x(dy)\, \t {\bf m}_n(dx) \\ 
&&\hspace{.1cm}\stack{n\to\infty}{\lra}-\sum_{i=1}^r\int \frac{\partial 
\vp}{\partial u_i}\cdot (h_i,\, \cdot\, )z'(\, \cdot\, ,0)\, d\bnu \\
&&\hspace{.5cm}=\langle Cf,\1\rangle\vphantom{\int}\, .  
\end{eqnarray*} 
{\it Step 4 } We prove (c). Recalling (\ref{4.38}) we find 
\begin{eqnarray*} 
&&\hspace{-.5cm}\lim_{n\to\infty}\left(\int_{x\in D^n}\int_{y\in\partial 
D^n}\t f\, d\mu^n_x\, \t {\bf m}_n(dx)-\int_{D^n}\t f\, d\t {\bf m}_n 
\right) \\ 
&&\hspace{.5cm}=\lim_{n\to\infty}\left(\int_{x\in D^n}\sum_{i=1}^n\int_{i 
,y_i}\hspace{-0.1cm}\int f\left({\T\frac1n\sum_{j\neq i}\delta_{z_j}+\frac 
1n\delta_{y_i}}\right)\times\right. \\ 
&&\hspace{1.0cm}\times\int_{t=0}^\ve\prod_{j\neq i}P_{x_j}\left(\tau_j 
>t,\, B_{j,t}\in dz_j\right)P_{x_i}\left(\tau_i\in dt,\, B_{i,t}\in d 
y_i\right)\, \t {\bf m}_n(dx)\left.-\int_{D^n}\t f\, d\t {\bf m}_n
\vphantom{\sum_1^1}\right)\, .  
\end{eqnarray*}  
We keep using the notation of Step 2 of the proof of Lemma \ref{Lemma4.4} 
and continue 
\begin{eqnarray*} 
&&\hspace{-.5cm}\lim_{n\to\infty}\left(\int_{x\in D^n}\int_{y\in\partial 
D^n}\t f\, d\mu^n_x\, \t {\bf m}_n(dx)-\int_{D^n}\t f\, d\t {\bf m}_n 
\right) \\ 
&&\hspace{.5cm}=\lim_{n\to\infty}\int_{x\in D^n}\int_{t=0}^\ve\sum_{i= 
1}^n\int_{y_i\in\partial D}\int_i\int\left(f\left({\T\frac1n\sum_{j\neq 
i}\delta_{z_j}+\frac1n\delta_{y_i}}\right)\right.\nonumber \\ 
&&\hspace{6.4cm}\left.-f\left({\T\frac1n\sum_{j\neq i}^n\delta_{z_j}+ 
\frac1n\delta_{x_i}}\right)\right)\times\nonumber \\ 
&&\hspace{2.0cm}\times\prod_{j\neq i}P_{x_j}\left(\tau_j>t,\, B_{j,t} 
\in dz_j\right)\, P_{x_i}\left(B_{i,t}\in dy_i,\, \tau_i\in dt\right)\, 
\t {\bf m}_n(dx)\vphantom{\sum_{i=1}^n} 
\end{eqnarray*}  
for any $\ve>0$. Proceeding as in Step 1 of this proof but replacing 
$z_i$ by $x_i$ in the definition of $r^n_{i,k}(y)$ we obtain 
\begin{eqnarray*}
&&\hspace{-0.5cm}\lim_{n\to\infty}\left(\int_{x\in D^n}\int_{\partial 
D^n}\t f\, d\mu^n_x\, \t {\bf m}_n(dx)-\int_{D^n}\t f\, d\t {\bf m}_n 
\right) \\ 
&&\hspace{0.0cm}=\lim_{n\to\infty}\frac1n\sum_{i=1}^n\int_i\int\sum_{k 
=1}^r\frac{\partial}{\partial u_k}\, \vp(\ldots ,(h_k,\mu^{i,y_i}), 
\ldots)\times \\ 
&&\hspace{1.0cm}\times\int_{x\in D^n}\int_{y_i\in\partial D}(h_k(x_i)- 
h_k(y_i))\times \\ 
&&\hspace{1.0cm}\times\int_{t=0}^\ve\prod_{j\neq i}P_{x_j}\left(\tau_j 
>t,\, B_{j,t}\in dz_j\right)\, P_{x_i}\left(\tau_i\in dt,\, B_{i,t}\in 
dy_i\right)\, \t {\bf m}_n(dx) \\ 
&&\hspace{.5cm}+\lim_{n\to\infty}\frac1n\sum_{i=1}^n\int_i\int\sum_{k 
=1}^r\int_{x\in D^n}\int_{y_i\in\partial D}r^n_{i,k}(y)\, (h_k(x_i)-h_k 
(y_i))\times \\ 
&&\hspace{1.0cm}\times\int_{t=0}^\ve\prod_{j\neq i}P_{x_j}\left(\tau_j 
>t,\, B_{j,t}\in dz_j\right)\, P_{x_i}\left(\tau_i\in dt,\, B_{i,t}\in 
dy_i\right)\, \t {\bf m}_n(dx)\, . 
\end{eqnarray*} 
We note that all $h_k(y_i)=0$ and recall that $\max_{j\in\{1,\ldots , 
n\}\setminus\{i\}}\left|r^n_{j,k}(y)\right|\stack{n\to\infty}{\lra}0$ 
for any sequence $y\equiv y^n\in\partial^{(1)}D^n$ with $(\mu^{i,y_i} 
)^n\stack{n\to\infty}\Ra\mu$ and every $k\in\{1,\ldots,r\}$. We get 
\begin{eqnarray*} 
&&\hspace{-0.5cm}\lim_{n\to\infty}\left(\int_{x\in D^n}\int_{\partial 
D^n}\t f\, d\mu^n_x\, \t {\bf m}_n(dx)-\int_{D^n}\t f\, d\t {\bf m}_n 
\right) \\ 
&&\hspace{0.0cm}=\lim_{n\to\infty}\frac1n\sum_{i=1}^n\int_{y_i\in 
\partial D}\int_i\int\sum_{k=1}^r\frac{\partial}{\partial u_k}\, \vp( 
\ldots ,(h_k,\mu^{i,y_i}),\ldots)\times \\ 
&&\hspace{1.0cm}\times\int_{x\in D^n}\int_{t=0}^\ve\prod_{j\neq i}^n 
P_{x_j}\left(\tau_j>t,\, B_{j,t}\in dz_j\right)\, h_k(x_i)P_{x_i}\left 
(\tau_i\in dt,\, B_{i,t}\in dy_i\right)\, \t {\bf m}_n(dx)\, . 
\end{eqnarray*} 
By using (\ref{4.40}), (\ref{4.41}), $P_{x_i}\left(\tau_j>t,\, B_{i,t} 
\in\cdot\right)\stack{t\to 0}\Ra\delta_{x_i}$, $\mu^{i,y_i}-\mu\stack{ 
n\to\infty}\Ra 0$, and $\t {\bf m}_n(dx)=-\frac12\Delta\t m_n(x)dx$ we 
conclude that 
\begin{eqnarray*} 
&&\hspace{-0.5cm}\lim_{n\to\infty}\left(\int_{x\in D^n}\int_{\partial 
D^n}\t f\, d\mu^n_x\, \t {\bf m}_n(dx)-\int_{D^n}\t f\, d\t {\bf m}_n 
\right) \\ 
&&\hspace{0.5cm}=\lim_{n\to\infty}\int_{z\in D^n}\sum_{k=1}^r\frac{ 
\partial}{\partial u_k}\, \vp(\ldots ,(h_k,\mu),\ldots)\times \\ 
&&\hspace{2.0cm}\times\frac1n\sum_{i=1}^n\lim_{t\downarrow 0}\int_{ 
x_i\in D} h_k(x_i)\frac{P_{x_i}(\tau_i\le t)}{t}\, dx_i\times \\ 
&&\hspace{2.0cm}\times\int_{x\in D^n}\int_{t=0}^\infty\prod_{j=1}^n 
P_{x_j}\left(\tau_j>t,\, B_{j,t}\in dz_j\right)\, dt\, \t {\bf m}_n 
(dx)\, . 
\end{eqnarray*} 
With $\bnu(H)=1$, cf. condition (i) and $H(1)\subset H$, we arrive at
\begin{eqnarray*} 
&&\hspace{-0.5cm}\lim_{n\to\infty}\left(\int_{x\in D^n}\int_{\partial 
D^n}\t f\, d\mu^n_x\, \t {\bf m}_n(dx)-\int_{D^n}\t f\, d\t {\bf m}_n 
\right) \\ 
&&\hspace{0.5cm}=\int\sum_{k=1}^r\frac{\partial}{\partial u_k}\, \vp 
(\ldots,(h_k,\mu),\ldots )\cdot\left(-{\T\frac12\Delta}h_k,\mu\right) 
\, \bnu(d\mu)=\langle Bf,\1\rangle\, . 
\end{eqnarray*} 
\qed 

\section{A Class of Measures Satisfying (i)-(iii),(j) of Subsection 2.3 
and (k)-(kw) of Section 3} 

\setcounter{equation}{0}

Next we are going to introduce the class of initial configurations we are 
particularly interested in. 
\begin{itemize} 
\item[(l)] There exists $c>1$ such that $\bnu$ is concentrated on the 
class $\hat{H}$ consisting of all probability measures $d(x)\, dx\in E$ 
with $d\in C^\infty(\overline{D})$ satisfying 
\begin{eqnarray*} 
\frac{h_1}{c}\le d\le ch_1\quad\mbox{\rm and}\quad (-\lambda_1)\cdot 
\frac{h_1}{c}\le\left(-{\T\frac12}\Delta\right)d\le (-\lambda_1)\cdot ch_1 
\end{eqnarray*} 
such that $d':=\sum_{i=1}^\infty e^{-\lambda_i}(d,h_i)\cdot h_i\in L^1(D)$ 
and $(d'/\|d'\|_{L^1})\, dx\in E\cap H$. 

\item[(ll)] With the notation $z=(z_1,\ldots ,z_n)$ and $z_1,\ldots ,z_n\in 
D$, for all $n\in {\Bbb N}$, 
\begin{eqnarray*} 
\t \bnu_n(dz):=\int_{\{\mu=d(x)\, dx\in\hat{H}\}}\prod_{i=1}^nd(z_i)\, \bnu 
(d\mu)\, dz_1\ldots dz_n\, . 
\end{eqnarray*} 
\end{itemize}
\begin{proposition}\label{Proposition4.6} 
Suppose (l) and (ll). Then we have (i)-(iii),(j) of Subsection 2.3 and (kk)-(kw) 
of Section 3.
\end{proposition} 
Proof. {\it Step 1 } It is obvious that (l) implies (ii),(iii), and (j). 
Given $\frac{h_1}{c}\le d\le ch_1$ and $(-\lambda_1)\cdot\frac{h_1}{c}\le 
\left(-{\T\frac12}\Delta\right)d\le (-\lambda_1)\cdot ch_1$ $\bmu$-a.e. for 
some $c>1$, by the definitions of $H$ and $H(1)$, the last line of (l) is 
just another way to formulate (i). 
\medskip 

\nid 
{\it Step 2 } Conditions (kk) and (kkk) follow directly from (l) and (ll). 
We continue with the verification of (kw). With $f\in\t C_b^2(E)$ with $f 
(\mu)=\vp((h_1,\mu),\ldots,(h_r,\mu))$, $\mu\in E$, as in given (\ref{4.3}) 
and $\t f(z)=f\left({\T\frac1n\sum_{j=1}^n\delta_{z_j}}\right)$ on $z=(z_1, 
\ldots ,z_n)$ on $D^n$, we have by (l) and (ll) 
\begin{eqnarray*} 
\int_E f\, d\bnu_n&&\hspace{-.5cm}=\int_{D^n}\t f(z)\int_{\{\mu=d(x)\, dx 
\in\hat{H}\}}\prod_{i=1}^nd(z_i)\, \bnu(d\mu)\, dz_1\ldots dz_n \\ 
&&\hspace{-.5cm}=\int_{\{\mu=d(x)\, dx\in\hat{H}\}}\int_{D^n}f\left({\T\frac 
1n\sum_{j=1}^n\delta_{z_j}}\right)\prod_{j=1}^nd(z_j)\, dz\bnu(d\mu)\, . 
\end{eqnarray*} 
Now, suppose we are give a sequence $Z_1,Z_2,\ldots $ of independent $D 
$-valued random variables, all with probability distribution $d(x)\, dx\in 
\hat{H}$. Then the weak law of large numbers shows that 
\begin{eqnarray*} 
&&\hspace{-.5cm}\int_{D^n}f\left({\T\frac1n\sum_{j=1}^n\delta_{z_j}}\right) 
\prod_{j=1}^nd(z_j)\, dz \\ 
&&\hspace{.5cm}=\int_{D^n}\vp\left(\left(h_1,{\T\frac1n\sum_{j=1}^n\delta_{ 
z_j}}\right),\ldots,\left(h_r,{\T\frac1n\sum_{j=1}^n\delta_{z_j}}\right) 
\right)\prod_{j=1}^n\, P(Z_j\in dz_j) \\ 
&&\hspace{.5cm}={\Bbb E}\vp\left({\T\frac1n\sum_{j=1}^nh_1(Z_j)},\ldots,{\T 
\frac1n\sum_{j=1}^nh_r(Z_j)}\right)\vphantom{\int_{D^n}}\\ 
&&\hspace{.1cm}\stack{n\to\infty}{\lra}\vp\left({\T\frac1n\sum_{j=1}^n{\Bbb 
E}h_1(Z_j)},\ldots,{\T\frac1n\sum_{j=1}^n{\Bbb E}h_r(Z_j)}\right)\vphantom 
{\int_{D^n}} \\ 
&&\hspace{.5cm}=\vp\left((h_1,d),\ldots ,(h_r,d)\right)=f\left(d(x)\, dx 
\vphantom{l^1}\right)\vphantom{\int_{D^n}}
\end{eqnarray*} 
boundedly in $d(x)\, dx\in\hat{H}$. The last two relations combined with 
condition (l) give finally 
\begin{eqnarray*} 
\int_Ef\, d\bnu_n\stack{n\to\infty}{\lra}\int_{\{\mu=d(x)\, dx\in\hat{H}\}} 
f(\mu)\, d\bnu=\int_Ef\, d\bnu\, . 
\end{eqnarray*} 
We have verified (kw). 
\qed
\medskip 

For $z_1,\ldots ,z_n,v\in {\Bbb R}^d$ and $i\in \{1,\ldots ,n\}$, let $z^{ 
(i)}(v)$ denote the vector $(z_1,\ldots ,z_{i-1},v,$ $z_{i+1},\ldots,z_n)$. 
For a real function $f$ on ${\Bbb R}^{n\cdot d}$, let $\Delta_if$ denote 
the $d$-dimensional Laplace operator with respect to $i$-th $d$-dimensional 
argument applied to $f$. For example, for $f$ with argument $z^{(i)}(v)$, 
$\Delta_if$ is the Laplace operator with respect to $v$ applied to $f$. 
For $i\in\{1,\ldots ,n\}$, let us define $\mu^{n,i}:=\frac1n\sum_{j\neq 
i}\delta _{z_j}$, 
\begin{eqnarray*} 
{\cal L}^i(z\{i\},\mu):=\int_D\left.\left(-{\T\frac12}\Delta d\right)(x) 
\right/d(x)\, \mu^{n,i}(dx) 
\end{eqnarray*} 
for $z\{i\}=\{z_1,\ldots ,z_{i-1},z_{i+1},\ldots ,z_n\}$ and $\mu=d(x)\, 
dx\in\hat{H}$, and 
\begin{eqnarray*} 
{\cal K}^i(z\{i\}):=\int_{v\in D}\left(-{\T\frac12}\Delta_i\t m_n\right) 
\left(z^{(i)}(v)\right)\, dv\, ,  
\end{eqnarray*} 
Choose 
\begin{itemize} 
\item[(lll)] 
\begin{eqnarray*} 
\eta_{n,z\{i\}}(z_i)&&\hspace{-.5cm}:=\frac{1}{{\cal K}^i(z\{i\})}\D 
\int_{\{\mu=d(x)\, dx\in \hat{H}\}}{\cal L}^i(z\{i\},\mu)\cdot\prod_{j=1 
}^nd(z_j)\, \bnu(d\mu) 
\end{eqnarray*} 
\end{itemize}
The subsequent representation of $\eta_{n,z\{i\}}$ turns out to be useful 
in the technical calculations of the present subsection, 
\begin{eqnarray*} 
\eta_{n,z\{i\}}(z_i)=\frac{\D\int_{\{\mu=d(x)\, dx\in\hat{H}\}}\int_D 
\left.\left(-{\T\frac12}\Delta d\right)(x)\right/d(x)\, \mu^{n,i}(dx) 
\cdot d(z_i)\cdot\prod_{j\neq i}d(z_j)\, \bnu(d\mu)}{\D\int_{\{\mu=d(x) 
\, dx\in \hat{H}\}}\int_D\left(-{\T\frac12}\Delta d\right)(x)\, dx\cdot 
\prod_{j\neq i}d(z_j)\, \bnu(d\mu)}\, . 
\end{eqnarray*} 
\begin{proposition}\label{Proposition4.7} 
(a) Suppose (l)-(lll). Then we have (k1) and (k2). \\ 
(b) Let $\Phi:\hat{H}\to {\Bbb R}$ be bounded on $\hat{H}$ and continuous 
in $\mu^0\equiv d^0(x)\, dx\in\hat{H}$ in the sense of $\lim_{\|d-d^0\|_{ 
L^1(D)}\to 0}\left|\Phi(\mu)-\Phi(\mu^0)\right|=0$ where $\mu\equiv d(x) 
\, dx\in\hat{H}$. Then 
\begin{eqnarray*} 
\frac{\D\int_{\{\mu=d(x)\, dx\in\hat{H}\}}\Phi(\mu)\prod_{j=1}^nd(z_j)\, 
\bnu(d\mu)}{\D\int_{\{\mu=d(x)\, dx\in\hat{H}\}}\prod_{j=1}^nd(z_j)\, \bnu 
(d\mu)}\stack{}{\lra}\Phi(\mu^0)\quad\mbox{\rm as}\quad{\T\frac1n\sum_{j=1 
}^n}\delta_{z_j}\stack{n\to\infty}{\Ra}\mu^0. 
\end{eqnarray*} 
\end{proposition} 
Proof. The above definition yields directly (k1). We verify (k2). For $\mu 
\equiv d(x)\, dx\in\hat{H}$, $i\in\{1,\ldots ,n\}$, and $z_i\in D$ set 
\begin{eqnarray*} 
\hat{\bnu}(d\mu)\equiv\hat{\bnu}_{z_i}(d\mu):=\int_{x\in D}\left(-{\T\frac12} 
\Delta d\right)(x)\, dx\cdot\frac{1}{d(z_i)}\, \bnu(d\mu) 
\end{eqnarray*} 
and choose $\mu^0\equiv d^0(x)\, dx\in\hat{H}$. We note that (l) implies 
$\hat{\bnu}(\hat{H})<\infty$. Let $\Phi:\hat{H}\to {\Bbb R}$ be as in the 
formulation of part (b) of this proposition. Define 
\begin{eqnarray*} 
g_\Phi\left({\T\frac1n\sum_{j=1}^n\delta_{z_j}}\right):=\int_{\{\mu=d(x)\, 
dx\in\hat{H}\}}\Phi(\mu)\prod_{j=1}^nd(z_j)\, \hat{\bnu}(d\mu)\, . 
\end{eqnarray*} 
According to the weak law of large numbers, for all $f\in\t C^2_b(E)$
\begin{eqnarray*} 
&&\hspace{-.5cm}\int_{D^n}g_\Phi\left({\T\frac1n\sum_{j=1}^n\delta_{z_j}} 
\right)\cdot f\left({\T\frac1n\sum_{j=1}^n\delta_{z_j}}\right)\, dz \\ 
&&\hspace{.5cm}=\int_{\{\mu=d(x)\, dx\in\hat{H}\}}\int_{D^n}f\left({\T\frac 
1n\sum_{j=1}^n\delta_{z_j}}\right)\prod_{j=1}^nd(z_j)\, dz\cdot\Phi(\mu)\, 
\hat{\bnu}(d\mu) \\ 
&&\hspace{0.1cm}\stack{n\to\infty}{\lra}\int_{\{\mu=d(x)\, dx\in\hat{H}\}} 
f(\mu)\Phi(\mu)\, \hat{\bnu}(d\mu)\, . 
\end{eqnarray*} 
In particular, we have demonstrated that 
\begin{eqnarray*} 
\frac{\D\int_{\{\mu=d(x)\, dx\in\hat{H}\}}\Phi(\mu)\prod_{j=1}^nd(z_j)\, 
\hat{\bnu}(d\mu)}{\D\int_{\{\mu=d(x)\, dx\in\hat{H}\}}\prod_{j=1}^nd(z_j) 
\, \hat{\bnu}(d\mu)}=\frac{g_\Phi\left({\T\frac1n\sum_{j=1}^n\delta_{z_j}} 
\right)}{g_{\1}\left({\T\frac1n\sum_{j=1}^n\delta_{z_j}}\right)}\stack{} 
{\lra}\Phi(\mu^0) 
\end{eqnarray*} 
as $\frac1n\sum_{j=1}^n\delta_{z_j}\stack{n\to\infty}{\Ra}\mu^0$. Repeating 
the previous calculations with $\bnu$ instead of $\hat{\bnu}$, we get part 
(b) of the present lemma. Let us continue to examine the ``$\hat{\bnu} 
$-case". The last limit shows that the measure $\prod_{j=1}^nd(z_j)\, \hat 
{\bnu}(d\mu)/\int\prod_{j=1}^nd(z_j)\, \hat{\bnu}(d\mu)$ converges weakly 
to the one point probability measure concentrated on $\mu^0$ as $\frac1n 
\sum_{j=1}^n\delta_{z_j}\stack{n\to\infty}{\Ra}\mu^0$. This implies also 
\begin{eqnarray*} 
\frac{\D\int_{\{\mu=d(x)\, dx\in\hat{H}\}}\Phi(\mu)\cdot\frac{\int_D\left({ 
\T-\frac12}\Delta d\right)(x)/d(x)\, \mu^{n,i}(dx)}{\int_D\left(-{\T\frac12} 
\Delta d\right)(x)\, dx}\cdot\prod_{j=1}^nd(z_j)\, \hat{\bnu}(d\mu)}{\D\int_{ 
\{\mu=d(x)\, dx\in\hat{H}\}}\prod_{j=1}^nd(z_j)\, \hat{\bnu}(d\mu)}\stack{} 
{\lra}\Phi(\mu^0) 
\end{eqnarray*} 
as $\frac1n\sum_{j=1}^n\delta_{z_j}\stack{n\to\infty}{\Ra}\mu^0$. 
Choosing now $\vp\in C(\overline{D})$, 
\begin{eqnarray*} 
\Phi(\mu):=\int_{x\in D}\vp(x)\, d(x)\, dx\, ,\quad\mu\equiv d(x)\, dx 
\in\hat{H}, 
\end{eqnarray*} 
and recalling the definition of $\hat{\bnu}$ we obtain 
\begin{eqnarray*} 
&&\hspace{-.5cm}\int_D\vp\, \eta_{n,z\{i\}}\, dx \\ 
&&\hspace{.5cm}=\frac{\D\int_{\{\mu=d(x)\, dx\in\hat{H}\}}\int_D\left. 
\left(-{\T\frac12}\Delta d\right)(x)\right/d(x)\, \mu^{n,i}(dx)\cdot\int_{ 
z_i\in D}\vp(z_i)d(z_i)\, dz_i\cdot\prod_{j\neq i}d(z_j)\, \bnu(d\mu)}{\D 
\int_{\{\mu=d(x)\, dx\in\hat{H}\}}\int_D\left(-{\T\frac12}\Delta d\right) 
(x)\, dx\cdot\prod_{j\neq i}d(z_j)\, \bnu(d\mu)} \\ 
&&\hspace{0.1cm}\stack{}{\lra}\int_D\vp\, d^0\, dx 
\end{eqnarray*} 
as $\frac1n\sum_{j=1}^n\delta_{z_j}\stack{n\to\infty}{\Ra}\mu^0$ with  
$\mu^0(dx)=d^0(x)\, dx$. We have verified (k2). 
\qed 
\medskip 

\nid 
{\bf Remarks }(1) on condition (jj). One goal of the paper is to 
demonstrate Mosco type convergence for the Fleming-Viot type system. 
In order to establish the integration by parts formula in Corollary 
\ref{Corollary4.3} which we will use for this, it is quite natural 
to require the differentiability in the form of (jj) of the related 
reference measure 
$\bnu$. 
\medskip

\nid
(2) on conditions (jjj). We recall also $T_tf(\mu)=f(U^+(t,\mu) 
)$ and that therefore $T_t:L^2(E,\bnu)\to L^2(E,\bnu)$ is equivalent 
with $d\bnu\circ U^+(-t,\mu)/d\bnu\in L^\infty(E,\bnu)$. 
Condition (jjj) is motivated by 
Corollary \ref{Corollary4.3} together with condition (c6(ii)) of 
Section 2. So far, we have used (jjj) just in Corollary 
\ref{Corollary4.3}. However it also will become important below, from 
Proposition \ref{Proposition4.10} (b) on. 
\medskip

\nid
(3) on condition (l). Looking at the spectral representation of the 
solution to the equation $\frac12\Delta u=\frac{\partial u}{\partial 
t}$ with initial condition $u(0,\cdot)=d'/\|d'\|_{L^1}$ it turns out 
that (l) is an implicit condition on the decay of $(d,h_i)$ with respect 
to $i\in {\Bbb N}$. 

\section{Integration by Parts on $E_n$} 
\setcounter{equation}{0}

\begin{lemma}\label{Lemma4.8} 
Assume (l)-(lll) and let $c$ be the constant specified in condition (l). 
Then we have for all $n\in {\Bbb N}$ and $z=(z_1,\ldots ,z_n)$ with $z_1, 
\ldots ,z_n\in D$ 
\begin{eqnarray*} 
\left|\frac1n\, \frac{\left(-\frac12\Delta\t m_n\right)(z)}{\t m_n(z)}- 
\frac{\D\eta_{n,z\{i\}}(z_i)\int_{v\in D}\left({\T-\frac12}\Delta_i\t m_n 
\right)\left(z^{(i)}(v)\right)\, dv}{\D\t m_n(z)}\right|\le\frac{-\lambda_1 
c^2}{n}\, . 
\end{eqnarray*} 
\end{lemma} 
Proof. Let $i\in\{1,\ldots ,n\}$. According to (ll), we have 
\begin{eqnarray}\label{4.43} 
&&\hspace{-.5cm}\frac1n\, \frac{\D\left({\T-\frac12}\Delta\t m_n\right) 
(z)}{\D\t m_n(z)}=\frac{\D\frac1n\sum_{j=1}^n\int\left({\T-\frac12}\Delta 
d\right)(z_j)\cdot\prod_{j'\neq j}d(z_{j'})\, d\bnu}{\D\t m_n(z)}\nonumber 
 \\ 
&&\hspace{.5cm}=\frac{\D\int\left(\int_D\left(\left({\T-\frac12}\Delta d 
\right)/d\vphantom{\dot{f}}\right)\, d\mu^{n,i}\cdot d(z_i)+\frac1n\left( 
{\T-\frac12}\Delta d\right)(z_i)\right)\cdot\prod_{j\neq i}d(z_j)\, d\bnu} 
{\D\t m_n(z)}\nonumber \\  
&&\hspace{.5cm}=\frac{\D\eta_{n,z\{i\}}(z_i)\int_{v\in D}\left({\T-\frac12 
\Delta_i}\t m_n\right)\left(z^{(i)}(v)\right)\, dv}{\D\t m_n(z)}+\frac1n\, 
\frac{\D\int\left(\left.\left({\T-\frac12}\Delta d\right)(z_i)\right/d(z_i 
)\right)\cdot\prod_{j=1}^nd(z_j)\, d\bnu}{\D\t m_n(z)} \nonumber \\ 
\end{eqnarray} 
where we have used (lll) for the last equality sign. Recalling (l) and 
(ll), the lemma follows now from 
\begin{eqnarray*} 
\left|\frac{\int\left(\left.\left({\T-\frac12}\Delta d\right)(z_i)\right/d 
(z_i)\right)\cdot\prod_{j=1}^nd(z_j)\, d\bnu}{\t m_n(z)}\right|\le\left\| 
\frac{\left({\T-\frac12}\Delta d\right)}{d}\right\|\le -\lambda_1c^2\, . 
\end{eqnarray*} 
\qed
\medskip 

We continue with an auxiliary lemma. For this we recall the definitions of 
$\t C_b^2(E)$ in (\ref{4.3}) and $\t C_0^2(E)$ in (\ref{4.4}). 
\begin{lemma}\label{Lemma4.9} 
(a) The set $\t C_b^2(E)$ is dense in $C_b(E)$. Both sets, $\t C_b^2(E)$ 
and $\t C_0^2(E)$, are dense in $L^2(E,\bnu_n)$, $n\in {\Bbb N}$, as well 
as $L^2(E,\bnu)$. \\ 
(b) Let $h_n\in {\cal C}$ such that $\langle h_n\, ,\, h_n\rangle_n$ is 
bounded in $n\in {\Bbb N}$. Furthermore, let $h\in L^2(E,\bnu)$ and assume 
$\langle f\, ,\, h_n \rangle_n\stack{n\to\infty}{\lra}\langle f\, ,\, h 
\rangle$ for all $f\in\t C_b^2(E)$ or for all $f\in\t C_0^2(E)$. Then $h_n 
\wstack{n\to\infty}{\lra}h$. \\ 
(c) For every $f\in C(E'_n)$ there exists $\vp\in {\cal C}$ with $\vp =f$ 
on $E_n$. \\ 
(d) Suppose (l) and (ll). For every $f\in L^2(E,\bnu)$ there exists $\vp\in 
{\cal C}$ with $\vp =f$ $\bnu$-a.e. on $E$. 
\end{lemma} 
Proof. Density of $\t C_b^2(E)$ in $C_b(E)$ is a consequence of the 
Stone-Weierstrass theorem. Anything else in part (a) is now obvious. To 
verify (b), let $\ve >0$ and $\vp\in {\cal C}$, and choose $f\in\t C_b^2 
(E)$ or $f\in\t C_0^2(E)$ with $\langle\vp-f\, ,\, \vp-f\rangle<\ve$. We 
have 
\begin{eqnarray*} 
\langle\vp-f\, ,\, \vp-f\rangle_n&&\hspace{-.5cm}=\langle\vp\, ,\, \vp 
\rangle_n+\langle f\, ,\, f\rangle_n-2\langle\vp\, ,\, f\rangle_n \\ 
&&\hspace{-1.0cm}\stack{n\to\infty}{\lra}\langle\vp\, ,\, \vp\rangle+\langle 
f\, ,\, f\rangle-2\langle\vp\, ,\, f\rangle \\ 
&&\hspace{-.5cm}=\langle\vp-f\, ,\, \vp-f\rangle<\ve
\end{eqnarray*} 
since $\langle\vp\, ,\, \vp\rangle_n\stack{n\to\infty}{\lra}\langle\vp\, , 
\, \vp\rangle$ by $\vp\in {\cal C}\subseteq {\cal D}$, $\langle f\, ,\, f 
\rangle_n\stack{n\to\infty}{\lra}\langle f\, ,\, f\rangle$ by $\bnu_n\stack 
{n\to\infty}{\Ra}\bnu$ according to (kw), and $\langle\vp\, ,\, f\rangle_n 
\stack{n\to\infty}{\lra}\langle\vp\, ,\, f\rangle$ by (c2). Recall also the 
introduction to Subsection 2.2, especially ${\cal F}=C_b(E)$ and (\ref{2.2}). 
Claim (b) follows now from 
\begin{eqnarray*} 
\left|\langle\vp\, , \, h_n\rangle_n-\langle f\, , \, h_n\rangle_n 
\vphantom{l^1}\right|\le\langle\vp-f\, ,\, \vp-f\rangle_n^{1/2}\cdot 
\langle h_n\, ,\, h_n\rangle_n^{1/2}\, , 
\end{eqnarray*} 
$\left|\langle\vp\, ,\, h\rangle-\langle f\, ,\, h\rangle\vphantom{l^1} 
\right|\le\langle\vp-f\, ,\, \vp-f\rangle^{1/2}\cdot\langle h\, ,\, h 
\rangle^{1/2}$, and the hypotheses of (b). 

For part (c) we mention that $E'_n$ is a closed set of the metric space 
$E$. By the Tietze extension theorem there is a continuous extension $F$ 
of $f$ to $E$, i. e., $F\in C_b(E)={\cal F}$, recall for this again the 
introduction to Subsection 2.2. The rest is trivial since $\vp:=f$ on $E_n$ 
and $\vp=g$ on $E\setminus E_n$ for any $g\in C_b(E)$ yield an element 
of ${\cal C}$ which proves claim (c). 

Let us focus on (d). We set $\vp\equiv\vp_f:=f$ on $E\setminus 
\bigcup_{n=1}^\infty E_n$. Furthermore, we define $\vp\equiv\vp_f$ 
on $E_n$, $n\in {\Bbb N}$, by 
\begin{eqnarray}\label{4.44} 
\quad\vp\left({\T\frac1n\sum_{i=1}^n}\delta_{z_i}\right):=\frac{\D 
\int_{\{\mu=d(x)\, dx\in\hat{H}\}}f(\mu)\prod_{j=1}^nd(z_j)\, \bnu(d 
\mu)}{\D\int_{\{\mu=d(x)\, dx\in\hat{H}\}}\prod_{j=1}^nd(z_j)\, \bnu 
(d\mu)} 
\end{eqnarray} 
where $z=(z_1,\ldots ,z_n)$ and $z_1,\ldots ,z_n\in D$. It follows from 
(l) and (ll) that $\vp$ is bounded and continuous on $E_n$, $n\in {\Bbb 
N}$. In order to show $\vp\in {\cal C}$ we have to verify (c2') of Subsection 
2.2 and $\vp\in {\cal D}$. Using now the particular structure of $\bnu_n$ 
and $\bnu$ this means we have to show
\begin{eqnarray}\label{4.45} 
\langle\vp\, ,\, \psi\rangle_n&&\hspace{-.5cm}=\int_{\{\mu=d(x)\, dx 
\in\hat{H}\}}\int_{z\in D^n}\vp\left({\T\frac1n\sum_{i=1}^n}\delta_{z_i 
}\right)\psi\left({\T\frac1n\sum_{i=1}^n}\delta_{z_i}\right)\prod_{j=1 
}^n d(z_j)\, dz_1\ldots dz_n\, d\bnu\nonumber \\ 
&&\hspace{-1.0cm}\stack{n\to\infty}{\lra}\int_{\hat{H}}\vp\psi\, d\bnu 
\, ,\quad\psi\in C_b(E), 
\end{eqnarray}
and
\begin{eqnarray}\label{4.46} 
\langle\vp\, ,\, \vp\rangle_n&&\hspace{-.5cm}=\int_{\{\mu=d(x)\, dx 
\in\hat{H}\}}\int_{z\in D^n}\vp^2\left({\T\frac1n\sum_{i=1}^n}\delta_{ 
z_i}\right)\prod_{j=1}^n d(z_j)\, dz_1\ldots dz_n\, d\bnu\nonumber \\ 
&&\hspace{-1.0cm}\stack{n\to\infty}{\lra}\int_{\hat{H}}\vp^2\, d\bnu 
\, .
\end{eqnarray}

We recall Proposition \ref{Proposition4.7} (b) and assume for a moment 
that $f$ is bounded and lower semicontinuous. Since $f$ is in this case 
the supremum of some sequence of bounded and continuous functions on $E$ 
we obtain 
\begin{eqnarray*} 
\liminf\vp\left({\T\frac1n\sum_{i=1}^n}\delta_{z_i}\right)\ge\vp(\mu) 
\quad\mbox{\rm as}\quad{\T\frac1n\sum_{j=1}^n}\delta_{z_j}\stack{n\to 
\infty}{\Ra}\mu 
\end{eqnarray*} 
for all $\mu\in\hat{H}$. Similarly, if $f$ was upper semicontinuous we 
would get 
\begin{eqnarray*} 
\limsup\vp\left({\T\frac1n\sum_{i=1}^n}\delta_{z_i}\right)\le\vp(\mu) 
\quad\mbox{\rm as}\quad{\T\frac1n\sum_{j=1}^n}\delta_{z_j}\stack{n\to 
\infty}{\Ra}\mu. 
\end{eqnarray*} 
We recall definition (\ref{4.44}) and keep the last two relations in 
mind. If $f$ is just bounded and measurable, the Vitali-Carath\'{e}odory 
theorem yields $\lim\vp\left({\T\frac1n\sum_{i=1}^n}\delta_{z_i}\right)= 
\vp(\mu)$ as $\frac1n\sum_{j=1}^n\delta_{z_j}$ $\stack{n\to\infty}{\Ra} 
\mu$ for $\bnu$-almost all $\mu\in\hat{H}$. For bounded measurable $f$ 
we get (\ref{4.45}) and (\ref{4.46}) immediately. For unbounded $f$ we 
use $f_r:=((-r)\vee f)\wedge r$, $r\in {\Bbb N}$, first and verify 
(\ref{4.45}) and (\ref{4.46}) for $f_r$. Then we observe that for all 
$n\in {\Bbb N}$
\begin{eqnarray*} 
&&\hspace{-.5cm}\langle\vp-\vp_r\, ,\, \vp-\vp_r\rangle_n\le\int_{\{ 
\mu=d(x)\, dx\in\hat{H}\}}\int_{z\in D^n}\D\int_{\{\mu'=d'(x)\, dx\in 
\hat{H}\}}(f(\mu')-f_r(\mu'))^2\vphantom{\frac{\prod_{j=1}^n}{\prod_{j 
=1}^n}}\times \\ 
&&\hspace{1.0cm}\times\left(\frac{\prod_{j=1}^n d'(z_j)}{\int_{\{\mu'' 
=d''(x)\, dx\in\hat{H}\}}\prod_{j=1}^nd''(z_j)\, \bnu(d\mu'')}\right) 
\, \bnu(d\mu')\prod_{j=1}^nd(z_j)\, dz\, \bnu(d\mu) \\ 
&&\hspace{.5cm}=\int_{\{\mu'=d'(x)\, dx\in\hat{H}\}}(f(\mu')-f_r(\mu') 
)^2\int_{z\in D^n}\prod_{j=1}^n d'(z_j)\, dz\, \bnu(d\mu')=\left\|f-f_r 
\right\|^2_{L^2(E,\sbnu)}\, . 
\end{eqnarray*} 
This and (\ref{4.45}) and (\ref{4.46}) for $f_r$ give now (\ref{4.45}) 
and (\ref{4.46}) for $f$ by letting $r\to\infty$. 
\qed
\bigskip

For $z_1,\ldots ,z_n\in D$, $z:=(z_1,\ldots ,z_n)$, $y:=(z_1,\ldots , 
z_{i-1},y_i,z_{i+1},\ldots ,z_n)$, $y_i\in\partial D$, $i\in\{1,\ldots 
,n\}$, and $g\in\t C_b^2(E)$ define 
\begin{eqnarray*} 
\zeta_g(z)&&\hspace{-.5cm}:=\int_{x\in D^n}\sum_{i=1}^n\int_{y_i\in 
\partial D}\frac{d\mu^n_x}{d\sigma}(y)\, s(dy_i)\cdot\eta(z_i)\cdot 
\left(-{\T\frac12}\Delta(\t g\t m_n)\right)(x)\, dx \\ 
&&\hspace{-.5cm}=\int_{x\in D^n}\frac{m_x^n(dz)}{dz}\left(-{\T\frac12} 
\Delta (\t g\t m_n)\right)(x)\, dx\, .  
\end{eqnarray*} 
\begin{proposition}\label{Proposition4.10} 
Suppose (l)-(lll). \\ 
(a) (Partial integration) All $g\in\t C_0^2(E)$ restricted to $E_n$ belong to 
$D(A_n')$ and $g=\1$ belongs also to $D(A_n')$ and we have 
\begin{eqnarray*}
\widetilde{A_n'g}=-\frac{\left(-\frac12\Delta(\t g\t m_n)\right)}{\t m_n}+ 
\frac{\zeta_g}{\t m_n}\in C(\overline{D^n})\, . 
\end{eqnarray*} 
Furthermore, $\|A_n'\1\|_{C_b(E_n)}$ is uniformly bounded in $n\in {\Bbb N}$ 
by $-\lambda_1c^2$ where $c$ is the constant specified in condition (l). For 
every $g\in \t C_0^2(E)$, the norm $\|A_n'g\|_{C_b(E_n)}$ is uniformly bounded 
in $n\in {\Bbb N}$. \\ 
(b) Assume, in addition, (jj) and (jjj). Let $g\in{\cal C}$, $f\in\t C_b^2(E)$, 
and $f_n\in D(\hat{A}_n)$, $n\in {\Bbb N}$, such that there exists a sequence 
$\vp_n\in {\cal C}$, $n\in {\Bbb N}$, with $\vp_n\wstack{n\to\infty}{\lra}f$ 
and $\langle f_n-\vp_n\, ,\, f_n-\vp_n\rangle_n\stack{n\to\infty}{\lra}0$. 
Then we have 
\begin{eqnarray*}
\langle\hat{A}_nf_n\, ,\, g\rangle_n\stack{n\to\infty}{\lra}\langle\hat{A}f\, 
,\, g\rangle\, . 
\end{eqnarray*}
(c) Suppose (jj) and (jjj). Let $g\in D(\hat{A})\cap {\cal C}$, $f\in L^2(E, 
\bnu)$, and $f_n\in D(\hat{A}_n)$, $n\in {\Bbb N}$, such that $\langle A_nf_n 
\, ,\, A_nf_n\rangle_n$ is bounded in $n\in {\Bbb N}$ and there exists a 
sequence $\vp_n\in {\cal C}$, $n\in {\Bbb N}$, with $\vp_n\wstack{n\to\infty} 
{\lra}f$ and $\langle f_n-\vp_n\, ,\, f_n-\vp_n\rangle_n\stack{n\to\infty} 
{\lra}0$. Then we have 
\begin{eqnarray*}
\langle-\hat{A}_nf_n\, ,\, g\rangle_n\stack{n\to\infty}{\lra}\langle f\, ,\, 
\hat{A}g\rangle\, . 
\end{eqnarray*}
\end{proposition} 
Proof. {\it Step 1 } We prepare the verification of part (a). Let $B^D_j$ 
denote the process obtained from $B_j$ by killing upon hitting $\partial D$. 
Furthermore, let $p^D(t,x_j,z_j)$ denote the transition density of $B^D_j$ 
from $x_j\in D$ to $z_j\in D$ in $t>0$ units of time, $j\in\{1,\ldots ,n\}$. 
We have 
\begin{eqnarray*}
\zeta_1(z)&&\hspace{-.5cm}=\int_{x\in D^n}\sum_{i=1}^n\int_{y_i\in\partial 
D}\frac{d\mu^n_x}{d\sigma}(y)\, s(dy_i)\cdot\eta(z_i)\cdot\left(-{\T\frac12} 
\Delta\t m_n\right)(x)\, dx \\ 
&&\hspace{-.5cm}=\int_{x\in D^n}\sum_{i=1}^n\int_{t=0}^\infty\prod_{j\neq i} 
\frac{P_{x_j}\left(\tau_j\ge t,\, B_{j,t}\in dz_j\right)}{dz_j}\, P_{x_i} 
\left(\tau_i\in dt\right)\cdot\eta(z_i)\cdot\left(-{\T\frac12}\Delta\t m_n 
\right)(x)\, dx \\ 
&&\hspace{-.5cm}=\sum_{i=1}^n\eta(z_i)\int_{x\in D^n}\int_{t=0}^\infty\prod_{ 
j\neq i} p^D(t,x_j,z_j)\, P_{x_i}\left(\tau_i\in dt\right)\left(-{\T\frac12} 
\Delta\t m_n\right)(x)\, dx\, .  
\end{eqnarray*} 
From, for example, the spectral resolution of $P_{x_i}\left(\tau_i\in dt 
\right)$ with respect to the eigenfunctions of $\Delta_i$ we obtain 
\begin{eqnarray}\label{4.47}
\zeta_1(z)&&\hspace{-.5cm}=\sum_{i=1}^n\eta(z_i)\int_{x\in D^n}\int_{t=0 
}^\infty\prod_{j\neq i}p^D(t,x_j,z_j)\int_{v\in D}p^D(t,x_i,v)\, dv\, dt 
\left(-{\T\frac12}\Delta_i\right)\left(-{\T\frac12}\Delta\t m_n\right)(x)\, 
dx\nonumber \\ 
&&\hspace{-.5cm}=\sum_{i=1}^n\eta(z_i)\int_{v\in D}\int_{x\in D^n}K^{D^n} 
\left(x,z^{(i)}(v)\right)\left(-{\T\frac12}\Delta_i\right)\left(-{\T\frac12} 
\Delta\t m_n\right)(x)\, dx\, dv\nonumber \\ 
&&\hspace{-.5cm}=\sum_{i=1}^n\eta(z_i)\int_{v\in D}\int_{x\in D^n}K^{D^n} 
\left(x,z^{(i)}(v)\right)\left(-{\T\frac12}\Delta\right)\left(-{\T\frac12} 
\Delta_i\t m_n\right)(x)\, dx\, dv\nonumber \\ 
&&\hspace{-.5cm}=\sum_{i=1}^n\eta_{n,z\{i\}}(z_i)\int_{v\in D}\left(-{\T\frac 
12}\Delta_i\t m_n\right)\left(z^{(i)}(v)\right)\, dv 
\end{eqnarray} 
where, for the third equality sign, we have used (l). This implies 
\begin{eqnarray}\label{4.48}
&&\hspace{-.5cm}-\frac{\left(-\frac12\Delta\t m_n\right)(z)}{\t m_n(z)}+ 
\frac{\zeta_1(z)}{\t m_n(z)}\nonumber \\ 
&&\hspace{.5cm}=\sum_{i=1}^n\left(-\frac1n\frac{\left(-\frac12\Delta\t m_n 
\right)(z)}{\t m_n(z)}+\frac{\D\eta_{n,z\{i\}}(z_i)\int_{v\in D}\left(-{\T 
\frac12}\Delta_i\t m_n\right)\left(z^{(i)}(v)\right)\, dv}{\D\t m_n(z)} 
\right)\, . 
\end{eqnarray} 
{\it Step 2 } We prove part (a). For this we keep Lemma \ref{Lemma4.4} (b) 
and (d) in mind. Let $f\in D^2_b(A_n)\cup\t C_b^2(E)\cup\t C_0^2 (E)$ and 
$g\in\t C_b^2(E)\cup\t C_0^2(E)$. Denoting by ${\bf n}_n$ the inner normal 
vector on $\partial D^n$ we obtain  
\begin{eqnarray}\label{4.49}
\left\langle{\T\frac12}\Delta\t f\, ,\, \t g\right\rangle_n&&\hspace{-.5cm} 
=\left\langle\t f\, ,\, \frac{{\T\frac12}\Delta (\t g\t m_n)}{\t m_n}\right 
\rangle_n+\frac12\int_{\partial D^n}\t f\cdot\frac{d(\t g \t m_n)}{d{\bf n 
}_n}\, d\sigma\nonumber \\ 
&&\hspace{-.5cm}=\left\langle\t f\, ,\, {\T\frac12}\Delta\t g\right\rangle_n 
+\left\langle\t f\, ,\, \frac{\nabla\t g\cdot\nabla\t m_n}{\t m_n}\right 
\rangle_n\nonumber \\ 
&&\hspace{0.0cm}-\int_{D^n}\t g\t f\, d\t {\bf m}_n+\int_{x\in D^n}\int_{ 
\partial D^n}\t g\t f\, d\mu^n_x\, \t {\bf m}_n(dx)\, . 
\end{eqnarray} 
By means of Lebesgue's theorem on monotone convergence we get (\ref{4.25}) 
and (\ref{4.26}) for all $f\in D^2_b(A_n)$. With this (\ref{4.49}) specifies 
for $g=\1$ and $f\in D^2_b(A_n)$ to 
\begin{eqnarray}\label{4.50}
\langle A_nf\, ,\, \1\rangle_n&&\hspace{-.5cm}=\left\langle{\T\frac12}\Delta 
\t f\, ,\, \t \1\right\rangle_n\nonumber \\ 
&&\hspace{-.5cm}=-\int_{D^n}\t f\, d\t {\bf m}_n+\int_{x\in D^n}\int_{\partial 
D^n}\t f\, d\mu^n_x\, \t {\bf m}_n(dx)\nonumber \\ 
&&\hspace{-.5cm}=-\int_{D^n}\t f\, d\t {\bf m}_n+\int_{x\in D^n}\int_{D^n}\t 
f\, dm^n_x\, \t {\bf m}_n(dx)\nonumber \\ 
&&\hspace{-.5cm}=-\left\langle\t f\, ,\, \frac{\left(-\frac12\Delta\t m_n 
\right)}{\t m_n}\right\rangle_n+\left\langle\t f\, ,\, \frac{\zeta_1}{\t m_n} 
\right\rangle_n\, . 
\end{eqnarray} 
According to (l)-(lll) and (\ref{4.47}) we have $\left(-\frac12\Delta\t m_n 
\right)/\t m_n\in C(\overline{D^n})$ and $\zeta_1/\t m_n\in C(\overline{D^n})$. 
By (\ref{4.50}), the result of Step 1, and Lemma \ref{Lemma4.8} we find that 
$g=\1$ belongs to $D(A_n')$, that we have 
\begin{eqnarray*}
\widetilde{A_n'\1}=-\frac{\left(-\frac12\Delta(\t m_n)\right)}{\t m_n}+\frac{ 
\zeta_1}{\t m_n}\in C(\overline{D^n})\, ,  
\end{eqnarray*} 
and that $\|A_n'\1\|_{C_b(E_n)}$ is uniformly bounded in $n\in {\Bbb N}$ by 
$-\lambda_1c^2$. Now, let $g\in\t C_0^2(E)$. Replacing $\t m_n$ with $\t g\t 
m_n$ relation (\ref{4.50}) implies  
\begin{eqnarray*}
\langle A_n f\, ,\, g\rangle_n&&\hspace{-.5cm}=\left\langle{\T\frac12}\Delta 
\t f\, ,\, \t g\right\rangle_n\nonumber \\ 
&&\hspace{-.5cm}=-\left\langle\t f\, ,\, \frac{\left(-\frac12\Delta(\t g 
\t m_n)\right)}{\t m_n}\right\rangle_n+\left\langle\t f\, ,\, \frac{\zeta_g} 
{\t m_n}\right\rangle_n\, ,\quad f\in D^2_b(A_n). 
\end{eqnarray*} 
We mention that the verification of $\zeta_g/\t m_n\in C(\overline{D^n})$ 
is similar to $\zeta_1/\t m_n\in C(\overline{D^n})$ above, $n\in {\Bbb N}$. 
\medskip 

It remains to show in the rest of the present Step 2 that, for $g\in\t C^2_0 
(E)$, $\|A_n'g\|_{C_b(E_n)}$ is uniformly bounded in $n\in {\Bbb N}$. For 
this, we replace in (\ref{4.47}) $\t m_n$ by $\t g\t m_n$ and obtain for all 
$z\in D^n$
\begin{eqnarray*}
&&\hspace{-.5cm}\widetilde{A_n'g}(z)=-\frac{\left(-\frac12\Delta(\t g\t m_n) 
\right)(z)}{\t m_n(z)}+\frac{\zeta_g(z)}{\t m_n(z)} \\ 
&&\hspace{0.5cm}=({\T\frac12}\Delta\t g)(z)+\sum_{i=1}^n\frac{\D\eta_{n,z 
\{i\}}(z_i)\int_{v\in D}\left(-{\T\frac12}\Delta_i\t g\right)\left(z^{(i)} 
(v)\right)\cdot\t m_n\left(z^{(i)}(v)\right)\, dv}{\D\t m_n(z)} \\ 
&&\hspace{1.0cm}+\frac{\nabla\t g(z)\cdot\nabla\t m_n(z)}{\t m_n(z)}-\sum_{ 
i=1}^n\frac{\D\eta_{n,z\{i\}}(z_i)\int_{v\in D}\nabla_i\t g\left(z^{(i)}(v) 
\right)\cdot\nabla_i\t m_n\left(z^{(i)}(v)\right)\, dv}{\D\t m_n(z)} \\ 
&&\hspace{1.0cm}+\t g(z)\cdot\widetilde{A_n'\1}(z)+\sum_{i=1}^n\frac{\D 
\eta_{n,z\{i\}}(z_i)\int_{v\in D}\left(\t g\left(z^{(i)}(v)\right)-\t g(z) 
\right)\cdot\left(-{\T\frac12}\Delta_i\t m_n\right)\left(z^{(i)}(v)\right) 
\, dv}{\D\t m_n(z)} \\ 
&&\hspace{0.5cm}=T_1(z)+T_2(z)+T_3(z)\, ,\quad z=(z_1,\ldots ,z_n)\, ,\ 
z_1,\ldots ,z_n\in D,\vphantom{\sum^l} 
\end{eqnarray*} 
linewise, where $\nabla_i$, similar to $\Delta_i$, denotes the $d 
$-dimensional gradient with respect to $i$-th $d$-dimensional argument 
and the product sign between two gradients stands for the scalar product 
in the corresponding Euclidean vector spaces. Furthermore, we note that 
$|\eta_{n,z\{i\}}(v)|\le c_1\cdot h_1(v)$ for some positive constant 
$c_1$ which is by the second line of (lll) and (l) uniformly bounded in 
$n\in {\Bbb N}$ and $z\{i\}$. From (\ref{4.4}) we take the representation 
$g(\mu)=\psi((h_1,\mu),\ldots,(h_r,\mu))\cdot\psi_0((k,\mu))\in\t C^2_0(E 
)$. With $\Psi(z):=\psi\left((h_1,\frac1n\sum_{j=1}^n\delta_{z_j}),\ldots 
,\right.$ $\left.(h_r,\frac1n\sum_{j=1}^n\delta_{z_j})\right)$, $\Psi_0(z 
):=\psi_0\left((k,\frac1n\sum_{j=1}^n\delta_{z_j})\right)$, 
\begin{eqnarray*}
&&\hspace{-.5cm}T(z,v)=\sum_{j=1}^r\lambda_j\cdot\frac{\partial\Psi} 
{\partial x_j}(z)\cdot\Psi_0(z)\cdot\frac1nh_j(v)+\Psi(z)\cdot\frac{\partial 
\Psi_0}{\partial x}(z)\cdot\frac1n\left({\T\frac12}\Delta k\right)(v) \\ 
&&\hspace{.5cm}+\frac{1}{2n}\sum_{j,j'=1}^r\frac{\partial^2\Psi}{\partial 
x_j\partial x_{j'}}(z)\cdot\Psi_0(z)\cdot\frac1n\left(\nabla h_j\cdot\nabla 
h_{j'}\right)(v) \\ 
&&\hspace{.5cm}+\frac{1}{2n}\left(2\sum_{j=1}^r\frac{\partial\Psi}{\partial 
x_j}(z)\, \frac{\partial\Psi_0}{\partial x}(z)\cdot\frac1n\left(\nabla h_j 
\cdot\nabla k\right)(v)+\Psi(z)\cdot\frac{\partial^2\Psi_0}{\partial x^2}(z) 
\cdot\frac1n\left(\nabla k\cdot\nabla k\right)(v)\right)
\end{eqnarray*} 
it holds that $({\T\frac12}\Delta\t g)(z)=\sum_{l=1}^nT(z,z_l)$ and 
$\left(-{\T\frac12}\Delta_i\t g\right)\left(z^{(i)}(v)\right)=T(z,v)$ 
and thus $\left\|({\T\frac12}\Delta\t g)\right\|\le c_2$ and $\left\| 
-{\T\frac12}\Delta_i\t g\right\|\le c_3/n$ with positive constants $c_2 
,c_3$ not depending on $n\in {\Bbb N}$ and $i\in\{1,\ldots ,n\}$. This 
implies 
\begin{eqnarray*}
|T_1(z)|&&\hspace{-.5cm}\le c_2+c_1c_3\,\frac{\frac1n\sum_{i=1}^n\left( 
\int_{v\in D}\t m_n\left(z^{(i)}(v)\right)\, dv\cdot h_1(z_i)\right)}{\D 
\t m_n(z)} \\ 
&&\hspace{-.5cm}\le c_2+c_1c_3\, \frac{\frac1n\sum_{i=1}^n\left(\int_{ 
\{\mu=d(x)\, dx\in\hat{H}\}}\prod_{j\neq i}^nd(z_j)\, \bnu(d\mu)\cdot 
h_1(z_i)\right)}{\D\t m_n(z)}
\end{eqnarray*} 
and with condition (l), $\|T_1\|_{C_b(E_n)}\le c_2+c_1c_3c$. We continue 
with the investigation of $T_2$. We have 
\begin{eqnarray*}
\nabla_i\t g(z)=\Psi_0(z)\cdot\frac1n\sum_{j=1}^r\frac{\partial\Psi} 
{\partial x_j}(z)\cdot\nabla h_j(z_i)+\Psi(z)\cdot\frac1n\frac{\partial 
\Psi_0}{\partial x}(z)\cdot\nabla k(z_i)\, ,\quad i\in \{1,\ldots ,n\}, 
\end{eqnarray*} 
and $\nabla\t g(z)=\left(\left(\nabla_1\t g(z)\right)^T,\ldots ,\left( 
\nabla_n\t g(z)\right)^T\right)^T$. Furthermore, we note that according 
to (\ref{4.4}) supp$\, \t g=K^n$ for some compact subset $K$ of $D$. 
Therefore, $\|\nabla_i\t g\|\equiv\|\nabla_i\t g\|_{C(\overline{D^n}; 
{\Bbb R}^d)}=\|\nabla_i\t g\|_{C(K^n;{\Bbb R}^d)}\le c_4/n$, $i\in\{1, 
\ldots ,n\}$, and $\|\nabla\t g\|\equiv\|\nabla\t g\|_{C(\overline{D^n} 
;{\Bbb R}^{n\cdot d})}=\|\nabla\t g\|_{C(K^n;{\Bbb R}^{n\cdot d})}\le 
c_4/\sqrt{n}$ for some $c_4>0$ which is independent on $n\in {\Bbb N}$. 

In addition, by condition (l), $\left(-\frac12\Delta d\right)$ is 
bounded on $D$ uniformly in $\mu=d(x)\, dx\in\hat{H}$. Accordingly, 
Gauss' theorem in thin tubes (cylinders) shows that $\nabla d$ is also 
bounded on $D$ uniformly in $\mu=d(x)\, dx\in\hat{H}$. In other words, 
there is $c_5>0$ such that $\|\nabla d\|\equiv \|\nabla d\|_{C(\overline 
{D};{\Bbb R}^d)}<c_5$, uniformly in $\mu=d(x)\, dx\in\hat{H}$. In 
particular, we get from (ll) 
\begin{eqnarray*}
\nabla_i\t m_n(z)=\int_{\{\mu=d(x)\, dx\in\hat{H}\}}\nabla d(z_i)\prod_{ 
j\neq i}^nd(z_j)\, \bnu (d\mu)\, ,\quad i\in \{1,\ldots ,n\}, 
\end{eqnarray*} 
and $\nabla\t m_n(z)=\left(\left(\nabla_1\t m_n(z)\right)^T,\ldots ,\left 
(\nabla_n\t m_n(z)\right)^T\right)^T$. With (l) and (ll) it follows that  
\begin{eqnarray*}
\left|\frac{\nabla\t g(z)\cdot\nabla\t m_n(z)}{\t m_n(z)}\right|\le c_4 
\cdot\sup_{v\in K\atop \mu=d(x)\, dx\in\hat{H}}\frac{|\nabla d(v)|}{d(v)} 
\le c_4c_5c\cdot\sup_{v\in K}h_1^{-1}(v)\, ,\quad z\in D^n, 
\end{eqnarray*} 
and 
\begin{eqnarray*}
\frac{\sum_{i=1}^n\left|\int_{v\in D}\nabla_i\t g\left(z^{(i)}(v)\right) 
\cdot\nabla_i\t m_n\left(z^{(i)}(v)\right)\, dv\cdot h_1(z_i)\right|}{\t 
m_n(z)}\le c_4c\cdot\sup_{v\in K\atop \mu=d(x)\, dx\in\hat{H}}|\nabla d(v)| 
\cdot |D|\, , 
\end{eqnarray*} 
$z\in D^n$, which gives $\|T_2\|_{C_b(E_n)}\le c_4c_5c\cdot\sup_{v\in K} 
h_1^{-1}(v)+c_1c_4c_5c\cdot |D|$. We turn to the estimation of $T_3$. It 
holds that 
\begin{eqnarray*}
\left|\t g\left(z^{(i)}(v)\right)-\t g(z)\right|\le\|\nabla_i\t g\|\cdot 
\left|z^{(i)}(v)-z\right|\le \frac1n\, c_4\cdot\mbox{\rm diam}(D)\, ,\quad 
z\in D^n,\ v\in D, 
\end{eqnarray*} 
where diam$(D)$ denotes the diameter of $D$. Together with conditions (l) 
and (ll) this yields 
\begin{eqnarray*}
\frac{\sum_{i=1}^n\left|\int_{v\in D}\left(\t g\left(z^{(i)}(v)\right)-\t 
g(z)\right)\cdot\left(-{\T\frac12}\Delta_i\t m_n\right)\left(z^{(i)}(v) 
\right)\, dv\cdot h_1(z_i)\right|}{\t m_n(z)}\le c_4c^3\cdot\mbox 
{\rm diam}(D) 
\end{eqnarray*} 
for all $z\in D^n$, $v\in D$, and $i\in\{1,\ldots ,n\}$. Therefore $\|T_3 
\|_{C_b(E_n)}\le\|g\|_{C_b(E_n)}\cdot (-\lambda_1)c^2+c_1c_4c^3\cdot 
$diam$(D)$. We have shown that $\|A_n'g\|_{C_b(E_n)}$ is uniformly bounded 
in $n\in {\Bbb N}$. 
\medskip

\nid 
{\it Step 3 } We demonstrate part (b). We get from Proposition 
\ref{Proposition4.5} (c) and (\ref{4.27}), (\ref{4.49}) for $f\in\t C_b^2 
(E)$ and $g\in\t C_0^2(E)$ 
\begin{eqnarray}\label{4.51}
\lim_{n\to\infty}\left\langle\t f\, ,\, \frac{\nabla\t g\cdot\nabla\t m_n} 
{\t m_n}\right\rangle_n=-2\langle f,Bg\rangle\, .  
\end{eqnarray} 
With this (\ref{4.49}), and (\ref{4.50}) we find for $f\in\t C_b^2(E)$ and 
$g\in\t C_0^2(E)$ 
\begin{eqnarray*}
&&\hspace{-.5cm}\left\langle f,A_n'g\right\rangle_n-{\T\frac12}\left\langle 
f,A_n'\1\cdot g\right\rangle_n=\left\langle\t f\, ,\, {\T\frac12}\Delta\t g 
\right\rangle_n+\left\langle\t f\, ,\, \frac{\nabla\t g\cdot\nabla\t m_n}{\t 
m_n}\right\rangle_n \\ 
&&\hspace{1.0cm}-\int_{D^n}\t g\t f\, d\t {\bf m}_n+\int_{x\in D^n}\int_{ 
\partial D^n}\t g\t f\, d\mu^n_x\, \t {\bf m}_n(dx) \\ 
&&\hspace{1.0cm}-\frac12\int\t g\t f\, d\left(\int_{x\in D^n}\, m^n_x\, \t 
{\bf m}_n(dx)-\t {\bf m}_n\right) \\ 
&&\hspace{.5cm}=\left\langle\t f\, ,\, {\T\frac12}\Delta\t g\right\rangle_n 
+\left\langle\t f\, ,\, \frac{\nabla\t g\cdot\nabla\t m_n}{\t m_n}\right 
\rangle_n-\frac12\int\t g\t f\, d\int_{x\in D^n}(m^n_x-\mu^n_x)\, \t {\bf 
m}_n(dx) \\ 
&&\hspace{1.0cm}+\frac12\int\t g\t f\, d\left(\int_{x\in {D^n}}\mu^n_x\, \t 
{\bf m}_n(dx)-\t {\bf m}_n\right)\, . 
\end{eqnarray*}
Using, in this order, (\ref{4.27}), (\ref{4.51}), and again Proposition 
\ref{Proposition4.5} (b),(c) we verify 
\begin{eqnarray*}
&&\hspace{-.5cm}\left\langle f,A_n'g\right\rangle_n-{\T\frac12}\left\langle 
f,A_n'\1\cdot g\right\rangle_n\vphantom{\frac12} \\ 
&&\hspace{-.0cm}\stack{n\to\infty}{\lra}\langle f,Bg\rangle-2\langle f,Bg 
\rangle-\frac12\left(\langle Cf,g\rangle+\langle f,Cg\rangle\right)+\frac12 
\left(\langle Bf,g\rangle+\langle f,Bg\rangle\vphantom{l^1}\right) \\ 
&&\hspace{.5cm}=\langle Bf,g\rangle-\frac12\left(\langle Bf,g\rangle+\langle 
f,Bg\rangle\vphantom{l^1}\right)-\frac12\left(\langle Cf,g\rangle+\langle f, 
Cg\rangle\vphantom{l^1}\right) \\ 
&&\hspace{.5cm}=\langle Bf,g\rangle-\frac12\left(\langle Af,g\rangle+\langle 
f,Ag\rangle\vphantom{l^1}\right)\, .  
\end{eqnarray*}
Now we use Corollary \ref{Corollary4.3}, Proposition \ref{Proposition4.5} (a), 
especially (\ref{4.37}), and recall the definitions of the space $\t C_0^2(E)$ 
in (\ref{4.4}) and the operator $C$ to derive 
\begin{eqnarray*}
\left\langle f,A_n'g\right\rangle_n-{\T\frac12}\left\langle f,A_n'\1\cdot g 
\right\rangle_n\, &&\hspace{-1cm}\stack{n\to\infty}{\lra}\langle Bf,g\rangle 
+\frac12\int fg\delta(A^f)\, d\bnu \\ 
&&\hspace{-.5cm}=\langle Bf,g\rangle-\frac12\langle A'\1\cdot f,g\rangle \\ 
&&\hspace{-.5cm}=\langle Af,g\rangle-\frac12\langle A'\1\cdot f,g\rangle \\ 
&&\hspace{-.5cm}=\langle\hat{A}f,g\rangle \, .  
\end{eqnarray*}
From Lemma \ref{Lemma4.2} as well as Corollary \ref{Corollary4.3} it follows 
that 
\begin{eqnarray*}
\left\langle f,A_n'g\right\rangle_n-{\T\frac12}\left\langle f,A_n'\1\cdot g
\right\rangle_n\stack{n\to\infty}{\lra}-\langle f,\hat{A}g\rangle\, ,\quad 
f\in\t C_b^2(E),\ g\in\t C_0^2(E)\, . 
\end{eqnarray*}
According to part (a) of this proposition and Lemma \ref{Lemma4.9} (c) there 
is a sequence $\psi_n\in {\cal C}$ such that $\psi_n=A_n'g-{\T\frac12}A_n'\1 
\cdot g$ on $E_n$ and $\|\psi_n\|_{C_b(E_n)}$ is uniformly bounded in $n\in 
{\Bbb N}$. Now Lemma \ref{Lemma4.9} (b) implies $\psi_n\wstack{n\to\infty} 
{\lra}-\hat{A}g$. In order to prove even $s$-convergence we observe that by 
Part (a) and (\ref{4.43})-(\ref{4.48}) 
\begin{eqnarray*} 
\widetilde{A_n'\1}(z)&&\hspace{-.5cm}=\frac{\D\frac1n\sum_{i=1}^n\int 
\left(\left.\left( {\T\frac12}\Delta d\right)(z_i)\right/d(z_i)\right) 
\cdot\prod_{j=1}^nd(z_j)\, d\bnu}{\D\t m_n(z)}\, , \quad z=(z_1,\ldots 
,z_n), 
\end{eqnarray*} 
$z_1,\ldots ,z_n\in D$ and thus for $\psi\in {\cal C}$ 
\begin{eqnarray}\label{4.52}
\left\langle A_n'\1\, ,\, \psi\right\rangle_n&&\hspace{-.5cm}=\int_{z\in 
D^n}\frac{\t \psi(z)}{n}\sum_{i=1}^n\int\left(\left.\left( {\T\frac12} 
\Delta d\right)(z_i)\right/d(z_i)\right)\cdot\prod_{j=1}^nd(z_j)\, d\bnu 
\, dz_1\ldots dz_n\nonumber \\ 
&&\hspace{-.5cm}=\int\int_{z\in D^n}\t \psi(z)\left(\frac{{\T\frac12} 
\Delta d}{d}\, ,\, {\T\frac1n\sum_{i=1}^n\delta_{z_i}}\right)\cdot\prod_{ 
j=1}^nd(z_j)\, dz_1\ldots dz_n\, d\bnu\nonumber \\ 
&&\hspace{-1.0cm}\stack{n\to\infty}{\lra}\int_{\{\mu=d(x)\, dx\in\hat{H} 
\}}\psi(\mu)\left({\T\frac12}\Delta d,\1\right)\, \bnu(d\mu)
\end{eqnarray} 
where, for the second line, Fubini's theorem applies because of 
condition (l) and, for the third line, we use $\bnu_n\Ra\bnu$, cf. 
(kw) and Proposition \ref{Proposition4.6}, Lemma \ref{Lemma2.1} 
(a), (d), and (\ref{2.2}) together with the remark underneath. 
Furthermore, 
\begin{eqnarray}\label{4.53} 
&&\hspace{-.5cm}\left\langle A_n'\1\cdot g\, ,\, A_n'\1\cdot g\right 
\rangle_n=\int\t g^2(z)\, \left(\frac{\frac1n\sum_{i=1}^n\int\left( 
\left.\left( {\T\frac12}\Delta d\right)(z_i)\right/d(z_i)\right)\cdot 
\prod_{j=1}^nd(z_j)\, d\bnu}{\t m_n(z)}\right)^2\, \t \bnu_n(dz) 
\nonumber \\ 
&&\hspace{.5cm}=\int_{z\in D^n}\t g^2(z)\left(\int\left(\frac{{\T 
\frac12}\Delta d}{d}\, ,\, {\T\frac1n\sum_{i=1}^n\delta_{z_i}}\right) 
\cdot\frac{\prod_{j=1}^nd(z_j)}{\t m_n(z)}\, d\bnu\right)^2\t m_n(z)\, 
dz_1\ldots dz_n\nonumber \\ 
&&\hspace{.5cm}=\int_{\{\mu=d'(x)\, dx\in\hat{H}\}}\int_{z\in D^n}\t 
g^2(z)\left(\vphantom{\frac{\prod_{j=1}^n}{\prod_{j=1}^n}}\int_{\{ 
\mu=d(x)\, dx\in\hat{H}\}}\left(\frac{{\T\frac12}\Delta d}{d}\, ,\, 
{\T\frac1n\sum_{i=1}^n\delta_{z_i}}\right)\times\right.\nonumber \\ 
&&\hspace{1.0cm}\left.\times\frac{\prod_{j=1}^nd(z_j)}{\int_{\{\mu=d(x) 
\, dx\in\hat{H}\}}\prod_{j=1}^nd(z_j)\, \bnu(d\mu)}\, \bnu(d\mu)\right 
)^2\cdot\prod_{j=1}^nd'(z_j)\, dz_1\ldots dz_n\, \bnu(d\mu')\nonumber 
 \\ 
&&\hspace{0.0cm}\stack{n\to\infty}{\lra}\int_{\{\mu'=d'(x)\, dx\in\hat 
{H}\}}g^2(\mu')\left({\T\frac12}\Delta d',\1\right)^2\, \bnu(d\mu')= 
\left\langle A'\1\cdot g\, ,\, A'\1\cdot g\right\rangle\, , 
\end{eqnarray} 
the limit by Proposition \ref{Proposition4.7} (b). In order to show 
$\left\langle A_n'g\, ,\, \psi\right\rangle_n\stack{n\to\infty}{\lra} 
\left\langle A'g\, ,\, \psi\right\rangle$ as well as $\left\langle 
A_n'g\, ,\, A_n'g\right\rangle_n\stack{n\to\infty}{\lra}\left\langle 
A'g\, ,\, A'g\right\rangle$ one may proceed as in Step 2 and individually 
prove convergence of $\left\langle T_i\, ,\, \psi\right\rangle_n$ and 
$\left\langle T_i\, ,\, T_i\right\rangle_n$ as $n\to\infty$ by plugging 
in (ll) as well as (lll) and repeating the way we have concluded in 
(\ref{4.52}) and (\ref{4.53}). We omit this long calculation here in the 
paper. 

One obtains finally $\psi_n\sstack{n\to\infty}{\lra}-\hat{A}g$. Combining 
this with ${\cal C}\ni\vp_n\wstack{n\to\infty}{\lra}f\in\t C_b^2(E)$, 
Proposition \ref{Proposition2.3} (c) gives 
\begin{eqnarray*}
\left\langle\vp_n,A_n'g\right\rangle_n-{\T\frac12}\left\langle\vp_n,A_n' 
\1\cdot g\right\rangle_n=\langle \vp_n\, ,\, \psi_n\rangle_n\stack{n\to 
\infty}{\lra}-\langle f,\hat{A}g\rangle\, ,\quad g\in\t C_0^2(E)\, . 
\end{eqnarray*}
For $f_n\in D(\hat{A}_n)$ with $\langle f_n-\vp_n\, ,\, f_n-\vp_n\rangle_n 
\stack{n\to\infty}{\lra}0$ we have therefore 
\begin{eqnarray*}
\left\langle f_n,A_n'g\right\rangle_n-{\T\frac12}\left\langle f_n,A_n'\1\cdot 
g\right\rangle_n=\langle f_n\, ,\, \psi_n\rangle_n\stack{n\to\infty}{\lra}- 
\langle f,\hat{A}g\rangle\, ,\quad g\in\t C_0^2(E)\, . 
\end{eqnarray*}
According to Corollary \ref{Corollary4.3} this yields 
\begin{eqnarray}\label{4.54}
\langle\hat{A}_nf_n\, ,\, g\rangle_n\stack{n\to\infty}{\lra}\langle\hat{A}f, 
g\rangle\, . 
\end{eqnarray} 
for all $g\in\t C_0^2(E)$. Recalling Lemma \ref{Lemma4.9} (a) one may choose 
a sequence $F_n\in\t C_0^2(E)$, $n\in {\Bbb N}$, with 
\begin{eqnarray}\label{4.55} 
\langle F_n-\hat{A}_nf_n,F_n-\hat{A}_nf_n\rangle_n\stack{n\to\infty}{\lra} 
0 
\end{eqnarray} 
and therefore 
\begin{eqnarray}\label{4.56}
\left\langle F_n\, ,\, g\right\rangle_n\stack{n\to\infty}{\lra}\langle\hat 
{A}f,g\rangle\, ,\quad g\in\t C_0^2(E)\, . 
\end{eqnarray} 
If one chooses in \cite{Lo13} $\, {\cal C}=\t C_0^2(E)$ then it follows 
from Proposition 3.3 (a) of the same reference that $\sup_{n\in{\Bbb N}} 
\langle F_n,F_n\rangle_n<\infty$.  Now (\ref{4.56}) and Lemma \ref{Lemma4.9} 
(b) imply that $F_n\wstack{n\to\infty}{\lra}\hat{A}f$. This and (\ref{4.55}) 
yield (\ref{4.54}) for all $g\in {\cal C}$; here ${\cal C}$ in the sense of 
the present paper. We have proved part (b) of the proposition. 
\medskip 

\nid
{\it Step 4 } We prove part (c). First let $g\in\t C_0^2(E)$ and note that 
in this case $g\in D(\hat{A})\cap{\cal C}$ by Lemma \ref{Lemma4.2} (b) and 
($c_1$'), ($c_2$') in the introduction to Subsection 2.2. As verified in the 
previous step, there is a sequence $\psi_n\in {\cal C}$ such that $\psi_n= 
A_n'g-{\T\frac12}A_n'\1\cdot g$ on $E_n$ and $\psi_n\sstack{n\to\infty} 
{\lra}-\hat{A}g$. Furthermore, $\sup_{n\in {\Bbb N}}\|\psi_n\|_{C_b(E_n)}< 
\infty$ by Step 2 and $\langle f_n-\vp_n\, ,\, f_n-\vp_n\rangle_n\stack{n 
\to\infty}{\lra}0$ by hypothesis which show $\langle f_n-\vp_n\, ,\, \psi_n 
\rangle_n\stack{n\to\infty}{\lra}0$. Using Proposition \ref{Proposition2.3} 
(c) and the hypothesis $\vp_n\wstack{n\to\infty}{\lra}f$ we obtain 
\begin{eqnarray*}
-\left\langle\hat{A}_nf_n\, ,\, g\right\rangle_n&&\hspace{-.5cm}=-\left 
\langle f_n\, ,\, A_n'g-{\T\frac12}A_n'\1\cdot g\right\rangle_n \\ 
&&\hspace{-.5cm}=-\langle\vp_n\, ,\, \psi_n\rangle_n -\langle f_n-\vp_n 
\, ,\, \psi_n\rangle_n \\ 
&&\hspace{-1.0cm}\stack{n\to\infty}{\lra}\langle f\, ,\, \hat{A}g\rangle 
\, . \vphantom{\left\langle\hat{A}_n\right\rangle_n}
\end{eqnarray*} 
For general $g\in D(\hat{A})\cap {\cal C}$ the claim is a consequence 
of approximation as follows. First we remind of the facts that by 
definition (\ref{2.15}), condition (jjj) together with Corollary 
\ref{Corollary4.3}, i. e., $A'\1\in L^\infty(E,\bnu)$, and Lemma 
\ref{Lemma4.2} (c) the set $\t C_0^2(E)$ is dense in $D(\hat{A})$ 
with respect to the graph norm $(\langle g\, ,\, g\rangle+\langle \hat 
{A}g\, ,\, \hat{A}g\rangle)^{1/2}$. In particular, $\t C_0^2(E)\ni g_r 
\stack{r\to\infty}{\lra}g\in D(\hat{A})\cap {\cal C}$ in the graph 
norm yields $\hat{A}g_r\stack{r\to\infty}{\lra}\hat{A}g$ as well as 
$g_r\stack{r\to\infty}{\lra}g$, both in $L^2(E,\bnu)$. Furthermore, 
$\sup_{n\in {\Bbb N}}\langle A_nf_n\, ,\, A_nf_n\rangle_n<\infty$ by 
hypothesis. Now, for every $r\in {\Bbb N}$ 
\begin{eqnarray*}
&&\hspace{-.5cm}\left|-\left\langle\hat{A}_nf_n\, ,\, g\right\rangle_n- 
\left\langle f\, ,\, \hat{A}g\right\rangle\right| \\ 
&&\hspace{.5cm}\le\left|\left\langle\hat{A}_nf_n\, ,\, g-g_r\right 
\rangle_n\right|+\left|-\left\langle\hat{A}_nf_n\, ,\, g_r\right 
\rangle_n-\left\langle f\, ,\, \hat{A}g_r\right\rangle\right|+\left| 
\left\langle f\, ,\, \hat{A}g_r-\hat{A}g\right\rangle\right| \\ 
&&\hspace{.5cm}\le\sup_{n\in {\Bbb N}}\left\langle\hat{A}_nf_n\, ,\, 
\hat{A}_nf_n\right\rangle_n^{1/2}\left\langle g-g_r\, ,\, g-g_r\right 
\rangle_n^{1/2}+\left|-\left\langle\hat{A}_nf_n\, ,\, g_r\right\rangle_n 
-\left\langle f\, ,\, \hat{A}g_r\right\rangle\right| \\ 
&&\hspace{1.0cm}+\left\langle f\, ,\, f\right\rangle^{1/2}\left\langle 
\hat{A}g_r-\hat{A}g\, ,\, \hat{A}g_r-\hat{A}g\right\rangle^{1/2} \\ 
&&\hspace{-.0cm}\stack{n\to\infty}{\lra}\sup_{n\in {\Bbb N}}\left\langle 
\hat{A}_nf_n\, ,\, \hat{A}_nf_n\right\rangle_n^{1/2}\left\langle g-g_r\, 
,\, g-g_r\right\rangle_n^{1/2}+\left\langle f\, ,\, f\right\rangle^{1/2} 
\left\langle\hat{A}g_r-\hat{A}g\, ,\, \hat{A}g_r-\hat{A}g\right\rangle^{ 
1/2}
\end{eqnarray*} 
\qed 

\section{Verifying Mosco Type Convergence} 
\setcounter{equation}{0}

We recall that under (jj),(jjj) and (l)-(lll), $C:={\textstyle\frac 
12}\|A'\1\|_{L^\infty(E,\sbnu)}\, \vee\, \sup_{n\in {\Bbb N}}{\T\frac 
12}\|A_n'\1\|_{L^\infty (E,\sbnu_n)}<\infty$, cf. Corollary 
\ref{Corollary4.3} and Proposition \ref{Proposition4.10} (a). 
\begin{proposition}\label{Proposition4.11} 
Suppose (jj),(jjj) and (l)-(lll). (a) For $\mu=d(x)\, dx\in\hat{H}$, i. e., 
for $\bnu$-a.e. $\mu\in E$, we have the representation $-\delta (A^f)(\mu)= 
A'\1(\mu)=\left({\T\frac12}\Delta d,\1\right)$.  \\ 
(b) For all $g\in L^2(E,\bnu)$, all sequences $g_n\in {\cal C}$ $s$-converging 
to $g$, and all $\beta>C$, we have $G_{n,\beta}g_n\sstack{n\to\infty}{\lra} 
G_\beta g$. 
\end{proposition}
Proof. Part (a) follows from (\ref{4.52}) and (\ref{4.57}) below where we 
also recall Corollary \ref{Corollary4.3}. In Step 1 we verify condition (c6). 
In Step 2 condition (c3') is verified. An immediate consequence is (c3) (for 
$G_{n,\beta}$, $n\in {\Bbb N}$, and $G_\beta$) and, by Lemma \ref{Lemma2.13}, 
(c3) for $\hat{G}_{n,\beta}$, $n\in {\Bbb N}$, and $\hat{G}_\beta$. In Steps 
3-5, we demonstrate that the forms $\hat{S}_n$ converge to the form $\hat{S}$ 
as $n\to\infty$ in the sense of Definition \ref{Definition2.4}. What we claim 
is then a consequence of Theorem \ref{Theorem2.14} and Remark (1) of Section 2. 
\medskip 

\nid 
{\it Step 1 } Let us verify condition (c6). Lemma \ref{Lemma4.2} (a) 
together with Proposition \ref{Proposition4.1} (a) implies 
\begin{eqnarray*} 
T_t'\1=\frac{d\bnu\circ U^+(-t,\mu)}{d\bnu}=\exp\left\{-\int_0^t\delta (A^f) 
(U^+(-s,\mu))\, ds\right\}\, . 
\end{eqnarray*} 
From this and (jjj), we obtain (c6$(ii)$). We also figure that $\1\in 
D(A')$, the first part of (c6$(i)$). Anything else required in (c6$(i)$) 
has been shown in Proposition \ref{Proposition4.10} (a). Condition 
(c6$(iv)$) has been verified in Lemma \ref{Lemma4.4} (d). 
\medskip

In the remainder of the present step, let us focus on (c6$(iii)$). 
Recalling (c1'), (c2'), and $\bnu_n\stack{n\to\infty}{\Ra}\bnu$ by (kw) 
and Proposition \ref{Proposition4.6}, we get the implication $\vp,\psi\in 
{\cal C}$ yields $\vp\cdot\psi\in {\cal C}$. Therefore, with Proposition 
\ref{Proposition4.10} (b), 
\begin{eqnarray*} 
\langle A_n'\1\cdot\vp\, ,\, \psi\rangle_n\ &&\hspace{-.5cm}=\ \ -2 
\left\langle\hat{A}_n\1\, ,\, \vp\cdot\psi\right\rangle_n \\ 
&&\hspace{-1.0cm}\stack{n\to\infty}{\lra}-2\left\langle\hat{A}\1\, , 
\, \vp\cdot\psi\right\rangle=\langle A'\1\cdot\vp\, ,\, \psi\rangle 
\, , \quad \psi\in{\cal C}. 
\end{eqnarray*} 
In other words, 
\begin{eqnarray}\label{4.57}
A_n'\1\cdot\vp\wstack{n\to\infty}{\lra}A'\1\cdot\vp\quad\mbox{\rm for all} 
\quad\vp\in {\cal C}\, . 
\end{eqnarray} 
In order to show that even $A_n'\1\cdot\vp\sstack{n\to\infty}{\lra} 
A'\1\cdot\vp$, $\vp\in {\cal C}$, we re-verify (\ref{4.53}) with $\vp$ 
instead of $g$. This completes the verification of (c6$(iii)$). 
\medskip 

\nid 
{\it Step 2 } According to (c1') in the introduction to Subsection 2.2, for 
$n\in {\Bbb N}$ we have ${\cal C}\subseteq L^\infty (E,\bnu_n)$. Condition 
(c3'$(ii)$), that is ${\cal G}_n=\{G_{n,\beta} g:g\in L^\infty (E,\bnu_n) 
\, ,\ \beta >0\}\subseteq {\cal C}$, is because of (c1') in Subsection 2.2 and 
Lemma \ref{Lemma4.9} (c) equivalent to $G_{n,\beta}g\in C(E'_n)$, $g\in 
L^\infty (E,\bnu_n)$. The latter is a repetition of the arguments in 
\cite{Lo09}, Proposition 1 and Step 2 of the proof of \cite{Lo09}, Lemma 
3 (a). Condition (c3'$(i)$) is a direct consequence of Lemma \ref{Lemma4.9} 
(d). 
\medskip 

\nid
{\it Step 3 } The objective of this step is to verify condition (i) of 
Definition \ref{Definition2.4} for $\hat{S}_n$, $n\in {\Bbb N}$, and $\hat
{S}$. Let $\vp\in L^2(E,\bnu)$. We assume that there is a  
subsequence $n_k$, $k\in {\Bbb N}$, of indices and $\vp_{n_k}\in D(A_{n_k 
})\cap {\cal C}$, $k\in {\Bbb N}$, such that $\vp_{n_k}\wstack{k\to\infty} 
{\lra}\vp$ and $\sup_{k\in {\Bbb N}}\left\langle A_{n_k}\vp_{n_k}\, ,\, 
A_{n_k}\vp_{n_k}\right\rangle_{n_k}<\infty$. Let $g\in D(A)\cap {\cal C}$. 
We note that Proposition \ref{Proposition4.10} (c) holds also for 
subsequences $n_k$, $k\in {\Bbb N}$, of indices. Thus, we have 
\begin{eqnarray*} 
\left\langle\hat{A}_{n_k}\vp_{n_k}\, ,\, g\right\rangle_{n_k}\stack{k\to 
\infty}{\lra}-\left\langle\vp\, ,\, \hat{A}g\right\rangle 
\end{eqnarray*}
which is equivalent to 
\begin{eqnarray*} 
\left\langle A_{n_k}\vp_{n_k}\, ,\, g\right\rangle_{n_k}-\left\langle 
{\T\frac12}A'_{n_k}\1\cdot \vp_{n_k}\, ,\, g\right\rangle_{n_k}\stack{k 
\to\infty}{\lra}-\left\langle\vp\, ,\, Ag\right\rangle+\left\langle{\T 
\frac12}A'\1\cdot\vp\, ,\, g\right\rangle 
\end{eqnarray*}
and 
\begin{eqnarray}\label{4.58}
\left\langle A_{n_k}\vp_{n_k}\, ,\, g\right\rangle_{n_k}\stack{k\to\infty} 
{\lra}-\left\langle\vp\, ,\, Ag\right\rangle+\left\langle A'\1\cdot\vp\, 
,\, g\right\rangle 
\end{eqnarray}
because of $A_n'\1\cdot g\sstack{n\to\infty}{\lra}A'\1\cdot g$ demonstrated 
in Step 1. By $\sup_{k\in{\Bbb N}}\left\langle A_{n_k}\vp_{n_k}\, ,\, A_{n_k} 
\vp_{n_k}\right\rangle_{n_k}<\infty$ and Proposition \ref{Proposition2.3} (b) 
for such a subsequence $\vp_{n_k}$, $k\in{\Bbb N}$, there is a (sub)subsequence 
$\vp_{n_l}$, $l\in {\Bbb N}$, and an element $F\in L^2(E,\bnu)$ with 
\begin{eqnarray}\label{4.59}
\left\langle A_{n_l}\vp_{n_l}\, ,\, g\right\rangle_{n_l}\stack{l\to\infty} 
{\lra}\left\langle F\, ,\, g\right\rangle\, . 
\end{eqnarray}
By Lemma \ref{Lemma4.2} (b), $\t C_b^2\subseteq D(A)\cap {\cal C}$ which 
together with Lemma \ref{Lemma4.2} (c) says that $D(A)\cap {\cal C}$ is 
dense in $D(A)$ with respect to the graph norm $(\langle f\, ,\, f\rangle 
+\langle Af\, ,\, Af\rangle)^{1/2}$. By (\ref{4.58}) and (\ref{4.59}) we 
obtain for $\vp\in L^2(E,\bnu)$, 
\begin{eqnarray*}
\left\langle\vp\, ,\, Ag\right\rangle=\left\langle A'\1\cdot\vp-F\, ,\, g 
\right\rangle\, ,\quad g\in D(A). 
\end{eqnarray*}
In other words, we have $\vp\in D(A')$ and $A'\vp=A'\1\cdot\vp-F$. Taking 
into consideration Corollary \ref{Corollary4.3} (b), we even get $\vp\in D 
(A)$. Recalling also Corollary \ref{Corollary4.3} (a), we obtain 
\begin{eqnarray*} 
\hat{S}(\vp,\vp)=\langle -A\vp\, ,\, \vp\rangle-{\textstyle\frac12}\left 
\langle\, \delta (A^f)\cdot\vp\, ,\, \vp\right\rangle=0 
\end{eqnarray*}
which verifies condition (i) of Definition \ref{Definition2.4}. 
\medskip 

\nid 
{\it Step 4 } In order to verify condition (ii) of Definition 
\ref{Definition2.4} for $\hat{S}_n$, $n\in {\Bbb N}$, and $\hat{S}$, let 
us first restrict ourselves to $\psi\in\t C_b^2(E)\left(\subseteq D(\hat 
{A})\right)$. For $n\in {\Bbb N}$, choose $\psi_n\in D(\hat{A}_n)$ with 
$\psi_n\sstack{n\to\infty}{\lra}\psi$. For this recall also $D(\hat{A}_n) 
=D(A_n)$, Lemma \ref{Lemma4.4} (d) and that all functions in ${\cal C}$ 
satisfy $(c1')$. In particular, the above choice guarantees
\begin{eqnarray}\label{4.60}
\langle\psi_n-\psi\, ,\, \psi_n-\psi\rangle_n\stack{n\to\infty}{\lra} 
0\, . 
\end{eqnarray} 
By Proposition \ref{Proposition4.10} (b) we have 
\begin{eqnarray*}
\langle-\hat{A}_n\psi_n\, , \, \vp\rangle_n&\stack{n\to\infty}{\lra}& 
\langle-\hat{A}\psi\, , \, \vp\rangle\, , \quad \vp\in {\cal C}.   
\end{eqnarray*} 
Lemma \ref{Lemma4.9} (c) yields the existence of a sequence $\Psi_n\in 
{\cal C}$ with 
\begin{eqnarray}\label{4.61}
\left\langle -\hat{A}_n\psi_n-\Psi_n\, ,\, -\hat{A}_n\psi_n-\Psi_n\right 
\rangle_n\stack{n\to\infty}{\lra}0\, .  
\end{eqnarray} 
It follows from the last two relations that $-\Psi_n\wstack{n\to\infty} 
{\lra}\hat{A}\psi$ and thus $\sup_{n\in {\Bbb N}}\langle\hat{A}_n\psi_n 
\, ,\, \hat{A}_n\psi_n\rangle_n<\infty$, cf. \cite{Lo14-1}, 
Proposition 2.3 (a). 

Furthermore, according to (\ref{4.60}) we have $\sup_{n\in{\Bbb N}} 
\langle\psi_n\, ,\, \psi_n\rangle_n<\infty$. With (\ref{4.61}) this shows 
$\left\langle -\hat{A}_n\psi_n-\Psi_n\, ,\, \psi_n\right\rangle_n\stack{n 
\to\infty}{\lra}0$. Proposition \ref{Proposition2.3} (c) implies now 
\begin{eqnarray*} 
\hat{S}_n(\psi_n,\psi_n)&&\hspace{-0.5cm}=\langle\Psi_n\, ,\, \psi_n 
\rangle_n+\left\langle-\hat{A}_n\psi_n-\Psi_n\, ,\, \psi_n\right\rangle_n 
 \\ 
&&\hspace{-1.0cm}\stack{n\to\infty}{\lra}\langle-\hat{A}\psi\, ,\, \psi 
\rangle=\hat{S}(\psi,\psi)\, .  
\end{eqnarray*} 
For general $\psi\in D(\hat{S})=D(\hat{A})$, we remind of the fact that 
by definition (\ref{2.15}), $A'\1\in L^\infty(E,\bnu)$ cf. condition (jjj) 
as well as Corollary \ref{Corollary4.3}, and Lemma \ref{Lemma4.2} (c), the 
set $\t C_b^2(E)$ is dense in $D(\hat{A})$ with respect to the graph norm 
$(\langle f\, ,\, f\rangle+\langle\hat{A}f\, ,\, \hat{A}f\rangle)^{1/2}$. 
Condition (ii) of Definition \ref{Definition2.4} follows from the above 
particular case $\psi\in\t C_b^2(E)(\subseteq D(\hat{A}))$ by 
approximation from Lemma \ref{Lemma4.9} (c) and Proposition 
\ref{Proposition4.10} (b) which we use as before. In particular we get 
$\sup_{n\in {\Bbb N}}\langle\hat{A}_n\psi_n\, ,\, \hat{A}_n\psi_n\rangle_n 
<\infty$ for general $\psi\in D(\hat{S})=D(\hat{A})$. 
\medskip 

\nid
{\it Step 5 } Let us verify condition (iv) for $\hat{A}_n$, $n\in {\Bbb 
N}$, and $\hat{A}$ which is because of Lemma \ref{Lemma2.6} sufficient 
for (iii) of Definition \ref{Definition2.4} for $\hat{S}_n$, $n\in {\Bbb 
N}$, and $\hat{S}$. For this, let $D(\hat{S}_n)\cap {\cal C}=D(S_n)\cap 
{\cal C}\ni u_n\wstack{n\to\infty}{\lra}v$ for some $v\in L^2(E,\bnu)$. 
Let $u\in D(S)=D(\hat{S})$ and suppose $\hat{A}_nu_n\in {\cal C}$ as well 
as $\beta u_n-\hat{A}_nu_n\wstack{n\to\infty}{\lra}\beta u-\hat{A}u$. One 
consequence is $\hat{A}_nu_n\wstack{n\to\infty}{\lra}\beta(v-u)+\hat{A}u$ 
which, by Proposition \ref{Proposition2.3} (a), says that $\sup_{n\in 
{\Bbb N}}\langle\hat{A}_nu_n\, ,\, \hat{A}_nu_n\rangle_n<\infty$. We have 
\begin{eqnarray*} 
\langle u\, ,\, \beta g+\hat{A}g\rangle &=&\langle\beta u-\hat{A}u\, ,\, 
g\rangle \\ 
&=&\lim_{n\to\infty}\langle\beta u_n\, ,\, g\rangle_n+\lim_{n\to\infty} 
\langle-\hat{A}_nu_n\, , \, g\rangle_n \\ 
&=&\langle v\, ,\, \beta g+\hat{A}g\rangle \vphantom{\lim_t}\, ,\quad g 
\in D(\hat{A})\cap {\cal C}. 
\end{eqnarray*} 
Here, the first line is true because of (\ref{4.23}) and the third 
line is true because of Proposition \ref{Proposition4.10} (c). Choose 
$g=\hat{G}_\beta f$, $f\in {\cal C}$. Recalling that by Step 2 of 
this proof and Lemma \ref{Lemma2.13} (a) it holds that $g\in D(\hat 
{S})\cap {\cal C}=D(\hat{A})\cap {\cal C}$ it turns out that 
\begin{eqnarray*} 
0&=&\langle u-v\, , \, 2\beta \hat{G}_\beta f-f\rangle\nonumber \\ 
&=&\langle 2\beta \hat{G}_\beta '(u-v)-(u-v)\, , \, f\rangle\, . 
\end{eqnarray*} 
By Lemma \ref{Lemma2.1} (c), $2\beta \hat{G}_\beta '(u-v)=u-v$ which 
implies $u-v\in D(\hat{A}')$ and $-\beta (u-v)=\hat{A}'(u-v)$. Since 
$\hat{A}'$ is an antisymmetric operator, $-\beta <0$ is not an 
eigenvalue. Consequently, $u=v$.  
\qed

\section{Relative Compactness} 
\setcounter{equation}{0}

In this subsection, we establish relative compactness of the family of 
processes ${\bf X}^n=((X^n_t)_{t\ge 0}$, $P_{{\sbnu}_n})$, $n\in {\Bbb 
N}$. Throughout the whole subsection we assume (jj),(jjj) and (l)-(lll), 
and $\beta>\|A'\1\|_{L^\infty(E,\sbnu)}\, \vee\, \sup_{n\in {\Bbb N}}\| 
A_n'\1\|_{L^\infty (E,\sbnu_n)}=2C$, recall for this Proposition 
\ref{Proposition4.10} (a) and the assumptions of the previous subsection. 
Let 
\begin{eqnarray*}
\t C_b(E):=\left\{f=\vp\left((h_1,\cdot),\ldots ,(h_r,\cdot)\right)\in\t 
C_b^2(E):\sum_{i=1}^r\frac{\partial\vp}{\partial x_i}\cdot (h_i,\cdot)=0 
\right\}\, . 
\end{eqnarray*}
For such an element $f\in \t C_b(E)$ we note that 
\begin{eqnarray}\label{4.62}
f-\frac1\beta Af=f-\frac1\beta Bf=:g\equiv g(f)\in {\cal C} 
\end{eqnarray}
by Lemma \ref{Lemma4.2}. For this we have recalled conditions (c1') and 
(c2') from the beginning of Subsection 2.2 and (\ref{2.2}). For $n\ge 2$, $F\in 
C_r(D^n)$, and $y=(z_1,\ldots z_{i-1},y_i,z_{i+1},\ldots ,z_n)\in\partial^{ 
(1)}D^n$ with $i\in\{1,\ldots ,n\}$, $y_i\in\partial D$, $z_1,\ldots ,z_n\in 
D$, set 
\begin{eqnarray*}
e_F(y):=\int_{z_i\in D}\left(F(z)-F(y)\right)\eta_{n,z\{i\}}(z_i)\, dz_i\, . 
\end{eqnarray*} 
For $n\ge 2$ identify $\partial^{(1)}D^n$ with $\partial D^n$. Let $f\in\t 
C_b(E)$ and $g\equiv g(f)$ as in (\ref{4.62}). Define 
\begin{eqnarray}\label{4.63}
b(n):=\frac{1}{(\beta -2C)\cdot\left\|f-\frac1\beta Bf\right\|^2}\cdot\int_{ 
x\in D^n}\int_{y\in\partial{D^n}}\left|e_{(\t f-\beta\widetilde{G_{n,\beta}g 
})^2}(y)\right|\, \mu^n_x(dy)\, \t {\bf m}_n(dx)\, .  
\end{eqnarray} 
Furthermore, for $B$ defined in (\ref{3.1}) and $n\in {\Bbb N}$, put $\hat 
{B}\equiv \hat{B}_n(\psi):=\{\psi^2\ge {\Bbb E}_\cdot e^{-\beta\tau_{B^c}} 
\}$ where $\psi$ is some bounded everywhere on $E_n$ defined function which 
we specify below. Let us denote by $E_\cdot$, $A$, $\hat{A}$, $\tau_{A^c} 
\equiv\tau_{A^c}^n$ the images of ${\Bbb E}_\cdot $, $B$, $\hat{B}$, $\tau_{ 
B^c}\equiv\tau_{B^c}^n$ under the map $\frac1n\sum_{i=1}^n\delta_{z_i}\to( 
z_1,\ldots ,z_n)$, $z_1,\ldots ,z_n\in\overline{D}$, taking into 
consideration invariance under permutations of $(z_1,\ldots ,z_n)$. 
\begin{lemma}\label{Lemma4.12} 
For $n\ge 2$, $f\in\t C_b(E)$, and $g\equiv g(f)$ the integral (\ref{4.63}) 
is finite and we have 
\begin{eqnarray*}
\lim_{n\to 0}b(n)=0\, .  
\end{eqnarray*} 
\end{lemma}
Proof. The proof uses ideas of the proof of Lemma 5 in \cite{Lo09}, 
Lemma 5. However the jump mechanism and the measures $\bnu_n$, $\bnu$ 
used in the present paper differ from those used in \cite{Lo09} and 
lead to different technical details which are necessary to point out. 
Note also that the notation in \cite{Lo09} is slightly different. 
\medskip 

\nid 
{\it Step 1 } Fix $n\in{\Bbb N}$ in this step and define $\Psi:=f- 
\beta G_{n,\beta}g$. As above, let $y=(z_1,\ldots z_{i-1},y_i,$ $z_{i+ 
1},\ldots ,z_i)\in\partial^{(1)}D^n$ with $i\in\{1,\ldots ,n\}$, $y_i 
\in\partial D$, $z_1,\ldots ,z_n\in D$. We have 
\begin{eqnarray*}
e_{\t\Psi^2}(y)&&\hspace{-.5cm}=\int_{z_i\in D}\left(\t\Psi(z)-\t\Psi(y) 
\right)\left(\t\Psi(z)+\t\Psi(y)\right)\eta_{n,z\{i\}}(z_i)\, dz_i \\ 
&&\hspace{-.5cm}=\int_{z_i\in D}\left(\t f(z)-\t f(y)\right)\left(\t f(z) 
+\t f(y)-\widetilde{\beta G_{n,\beta}g}(z)-\beta\widetilde{G_{n,\beta}g} 
(y)\right)\eta(z_i)\, dz_i \\ 
&&\hspace{-.0cm}-\int_{z_i\in D}\left(\beta\widetilde{G_{n,\beta}g}(z)- 
\beta\widetilde{G_{n,\beta}g}(y)\right)\left(\t f(z)-\beta\widetilde{G_{ 
n,\beta}g}(z)\right)\eta(z_i)\, dz_i 
\end{eqnarray*} 
because of (\ref{4.25}) which reads as $\int_{z_i\in D}\left(\beta 
\widetilde{G_{n,\beta}g}(z)-\beta\widetilde{G_{n,\beta}g}(y)\right)\eta(z_i 
)\, dz_i=0$. Using this again, we continue by 
\begin{eqnarray}\label{4.64}
e_{\t\Psi^2}(y)&&\hspace{-.5cm}=e_{\t f^2}(y)-2\int_{z_i\in D}\left(\t 
f(z)-\t f(y)\right)\cdot\beta\widetilde{G_{n,\beta}g}(z)\eta(z_i)\, dz_i 
\nonumber \\ 
&&\hspace{-0.0cm}+\int_{z_i\in D}\left(\left(\widetilde{G_{n,\beta}g}(z) 
\right)^2-\left(\widetilde{G_{n,\beta}g}(y)\right)^2\right)\eta(z_i)\, 
dz_i\nonumber \\ 
&&\hspace{-.5cm}=e_{\t f^2}(y)-2\int_{z_i\in D}\left(\t f(z)-\t f(y) 
\right)\cdot\beta\widetilde{G_{n,\beta}g}(z)\eta(z_i)\, dz_i+e_{ 
(\widetilde{G_{n,\beta}g})^2}(y)\, . \qquad
\end{eqnarray} 
{\it Step 2 } Let us consider the particular items of (\ref{4.64}). 
Since $f\in\t C_b(E)$ yields $f^2\in\t C_b(E)$, we obtain by Proposition 
\ref{Proposition4.5} (a) and (\ref{4.42}) 
\begin{eqnarray}\label{4.65} 
\int_{x\in D^n}\int_{y\in\partial{D^n}}\left|e_{\t f^2}(y)\right|\, 
\mu^n_x(dy)\, \t {\bf m}_n(dx)\stack{n\to\infty}{\lra}0\, . 
\end{eqnarray} 

For the second item in (\ref{4.64}) we recall Proposition 
\ref{Proposition4.5} (a) and the appendix below (\ref{4.42}). We use 
the inequality $\left|\widetilde{G_{n,\beta}g}(x)-\widetilde{G_{n, 
\beta}g}(y)\right|\le\frac1\beta\sup_{\left\{z\in D^n:z+y-x\in D^n 
\right\}}\left|g(z)-g(z+y-x)\vphantom{\displaystyle l^1}\right|$, $n 
\in {\Bbb N}$, and the fact that $g=f-\frac1\beta Af=f-\frac1\beta Bf$ 
is a bounded continuously differentiable cylindric function because of 
$f\in\t C_b(E)$. We obtain 
\begin{eqnarray}\label{4.66} 
&&\hspace{-.5cm}\int_{x\in D^n}\int_{y\in\partial{D^n}}\sum_{i=1}^n 
\chi_{(y)_i\in\partial D}(y)\left|\int_{z_i\in D}\left(\t f(z)-\t f(y) 
\right)\times\right.\nonumber \\  
&&\hspace{5.0cm}\left.\times\beta\widetilde{G_{n,\beta}g}(z)\eta(z_i) 
\, dz_i\vphantom{\int_{z_i\in D}}\right|\, \mu^n_x(dy)\, \t {\bf m}_n 
(dx)\stack{n\to\infty}{\lra}0\, . \qquad
\end{eqnarray} 

In order to handle the third item in (\ref{4.64}) we observe that by 
(\ref{4.25}) 
\begin{eqnarray*}
e_{(\beta\widetilde{G_{n,\beta}g})^2}=\int_{z_i\in D}\left(\beta\widetilde 
{G_{n,\beta}g}(z)-\beta\widetilde{G_{n,\beta}g}(y)\right)^2\eta(z_i)\, dz_i 
\ge 0\, . 
\end{eqnarray*} 
Furthermore, by Lemma \ref{Lemma4.4} (c) and multidimensional differential 
calculus, 
\begin{eqnarray*}
0&&\hspace{-.5cm}\le\int_{x\in D^n}\int_{\partial D^n}\left|e_{(\beta 
\widetilde{G_{n,\beta}g})^2}\right|\, d\mu^n_x\, \t {\bf m}_n(dx) \\ 
&&\hspace{-.5cm}=\int_{x\in D^n}\int_{\partial D^n}e_{(\beta\widetilde 
{G_{n,\beta}g})^2}\, d\mu^n_x\, \t {\bf m}_n(dx) \\ 
&&\hspace{-.5cm}=\int\left(\beta\widetilde{G_{n,\beta}g}\right)^2\, d 
\int_{x\in D^n}(m^n_x-\mu^n_x)\left({\T -\frac12\Delta}\t m_n\right)(x) 
\, dx \\ 
&&\hspace{-.5cm}=\int\left({\T-\frac12}\Delta\left(\beta\widetilde{G_{ 
n,\beta}g}\right)^2\right)\, d\t \bnu_n+\lim_{t\downarrow 0}\frac1t\int 
\left(\t T_{n,t}\left(\beta\widetilde{G_{n,\beta}g}\right)^2-\left( 
\beta\widetilde{G_{n,\beta}g}\right)^2\right)\, d\t \bnu_n \\ 
&&\hspace{-.5cm}=-\int\Delta\beta\widetilde{G_{n,\beta}g}\cdot\beta 
\widetilde{G_{n,\beta}g}\, d\t \bnu_n-\int\left(\nabla\beta\widetilde{ 
G_{n,\beta}g}\right)^2\, d\t \bnu_n+\int\left(\beta\widetilde{G_{n, 
\beta}g}\right)^2\, \widetilde{A'_n\1}\, d\t \bnu_n \\ 
&&\hspace{-.5cm}\le 2\beta\int\left(\t g-\beta\widetilde{G_{n,\beta}g} 
\right)\cdot\beta\widetilde{G_{n,\beta}g}\, d\t \bnu_n+\int\left(\beta 
\widetilde{G_{n,\beta}g}\right)^2\, \widetilde{A'_n\1}\, d\t \bnu_n 
\end{eqnarray*} 
By the Mosco type convergence of Proposition \ref{Proposition4.11} and 
Proposition \ref{Proposition2.3} (c), we get the limit 
\begin{eqnarray*} 
\int\left(\t g-\beta\widetilde{G_{n,\beta}g}\right)\cdot\beta\widetilde 
{G_{n,\beta}g}\, d\t \bnu_n\stack{n\to\infty}{\lra}\langle g-\beta 
G_\beta g\, ,\, \beta G_\beta g\rangle\, . 
\end{eqnarray*} 
In addition, $\vp:E\to{\Bbb R}$ given by $\vp=\beta G_{n,\beta}g$ on 
$E_n$, $n\in {\Bbb N}$, and $\vp=\beta G_\beta g$ $\bnu$-a.e. on $E 
\setminus\bigcup_{n=1}^\infty E_n$ is by Proposition \ref{Proposition4.11} 
an element of ${\cal C}$, cf. also (c1') and (c2') of Subsection 2.2. Thus, 
by (\ref{4.57}), Proposition \ref{Proposition4.11} and Proposition 
\ref{Proposition2.3} (c), 
\begin{eqnarray*}
&&\hspace{-.5cm}\int\left(\beta\widetilde{G_{n,\beta}g}\right)^2\, 
\widetilde{A'_n\1}\, d\t \bnu_n=\left\langle\beta G_{n,\beta}g\, ,\, 
A'_n\1\cdot\vp\right\rangle_n\nonumber \\ 
&&\hspace{1.0cm}\stack{n\to\infty}{\lra}\langle\beta G_\beta g\, ,\, 
A'\1\cdot\vp\rangle=\langle\beta G_\beta g\, ,\, A'\1\cdot\beta G_\beta 
g\rangle\, .  
\end{eqnarray*} 
Putting the last three relations together we conclude that 
\begin{eqnarray}\label{4.67} 
0&&\hspace{-.5cm}\le\int_{x\in D^n}\int_{\partial D^n}\left|e_{(\beta 
\widetilde{G_{n,\beta}g})^2}\right|\, d\mu^n_x\, \t {\bf m}_n(dx) 
\nonumber \\ 
&&\hspace{-.5cm}\le 2\beta\int\left(\t g-\beta\widetilde{G_{n,\beta}g} 
\right)\cdot\beta\widetilde{G_{n,\beta}g}\, d\t \bnu_n+\int\left(\beta 
\widetilde{G_{n,\beta}g}\right)^2\, \widetilde{A'_n\1}\, d\t \bnu_n 
\nonumber \\ 
&&\hspace{-1.0cm}\stack{n\to\infty}{\lra}2\beta\langle g-\beta G_\beta 
g\, ,\, \beta G_\beta g\rangle+\langle\beta G_\beta g\, ,\, A'\1\cdot 
\beta G_\beta g\rangle\vphantom{\int}\nonumber \\ 
&&\hspace{-.5cm}=-2\langle A\beta G_\beta g\, ,\, \beta G_\beta g 
\rangle +\langle\beta G_\beta g\, ,\, A'\1\cdot\beta G_\beta g\rangle=0 
\vphantom{\int}
\end{eqnarray} 
where, for the last equality sign, we recall $f=\beta G_\beta g\in\t 
C_b(E)\subseteq\t C_b^2(E)$ which implies $f^2\in \t C_b^2(E)\subseteq 
D(A)$ by Lemma \ref{Lemma4.2} (b). The conclusion ``$=0$" follows now 
from Corollary \ref{Corollary4.3}. The lemma is a consequence of 
(\ref{4.64}) and the estimates of its items, 
(\ref{4.65})-(\ref{4.67}). 
\qed 
\medskip 

Define $\nn\frac12\Del g\nn:=\sup_{\mu\in\bigcup_{n\in {\Bbb N}}E_n}| 
\frac12\Del g(\mu)|$ and set 
\begin{eqnarray*}
C_1=\|f\|+2\nn{\T\frac12}\Del f\nn+3\beta\|g\| 
\end{eqnarray*} 
as well as  
\begin{eqnarray*}
C_2:=\frac{C_1}{(\beta -2C)\cdot\left\|f-\frac1\beta Bf\right\|^2} 
\end{eqnarray*} 
where, for the constant $C$, we recall the introduction to the present 
subsection. Furthermore, for $n\ge 2$, let 
\begin{eqnarray*}
\ve_n:=\left(C_2\left(\int\left(f-\beta G_{n,\beta}g\right)^2\, d\bnu_n 
\right)^{1/2}+b(n)\right)^{1/3}\, . 
\end{eqnarray*} 
We observe $\ve_n\stack{n\to\infty}{\lra}0$ by the Mosco type convergence 
of Proposition \ref{Proposition4.11} and by Lemma \ref{Lemma4.12}. For $f$ 
and $g$ as in (\ref{4.62}) we specify 
\begin{eqnarray*}
\psi:=\frac{1}{\ve_n\left\|f-\frac1\beta Bf\right\|}\left(f-\beta G_{n, 
\beta}g\right)=\frac{1}{\ve_n\|g\|}\left(f-\beta G_{n,\beta}g\right)\, , 
\quad n\in {\Bbb N}. 
\end{eqnarray*} 
\begin{theorem}\label{Theorem4.13} 
The family of processes ${\bf X}^n=((X^n_t)_{t\ge 0}$, $P_{{\sbnu}_n})$, 
$n\in {\Bbb N}$, is relatively compact with respect to the topology 
of weak convergence of probability measures over the Skorokhod space $D_E 
([0,\infty))$. 
\end{theorem}
Proof. Recalling Theorem \ref{Theorem3.2} everything to verify is (c7) 
in Step 1 below and (c8) in Step 2 below. 
\medskip 

\nid
{\it Step 1 } The set $\t C_b(E)$ is an algebra containing the constant 
functions. Furthermore, we observe that 
\begin{eqnarray*}
&&\hspace{-.5cm}\left\{\Phi\left({\T\sum_{l=1}^{2k}}a_l\cdot\ln (h_l, 
\cdot)^2\right):\Phi\in C^2({\Bbb R})\ \mbox{\rm with both, }\vphantom 
{\sum_{l=1}^{2k}}\right. \\ 
&&\hspace{.5cm}\left.\lim_{x\to\infty}e^{bx}/\Phi(x)\ \mbox{\rm and } 
\lim_{x\to-\infty}e^{-bx}/\Phi(x)\ \mbox{\rm being finite for some } 
\right. \\ 
&&\hspace{.5cm}\left. b>\frac12,\ a_l\in\{-1,1\},\ \sum_{l=1}^{2k}a_l 
=0,\ k\in {\Bbb N}\right\} 
\end{eqnarray*}
can be represented by 
a subset of $\t C_b(E)$ that separates points in $E$. The rest of 
(c7) is a direct consequence of the particular choice of $\t C_b(E)$ 
which implies $f-\frac1\beta Af=f-\frac1\beta Bf=g\equiv g(f)\in {\cal 
C}$ if $f\in\t C_b(E)$ as already mentioned in (\ref{4.62}). 
\medskip 

\nid
{\it Step 2 } By $\beta>2C$ and (\ref{4.28}) we have 
\begin{eqnarray}\label{4.68} 
&&\hspace{-.5cm}(\beta-2C){\Bbb E}_{\sbnu_n}e^{-\beta\tau_{B^c}}= 
(\beta-2C)\int_{D^n}\t\psi^2\vee E_\cdot e^{-\beta\tau_{A^c}}\, d\t 
\bnu_n\nonumber \\ 
&&\hspace{.5cm}\le\beta\int_{D^n}\t \psi^2\vee E_\cdot e^{-\beta\tau_{ 
A^c}}\, d\t \bnu_n-\int_{D^n}\t \psi^2\vee E_\cdot e^{-\beta\tau_{A^c}} 
\cdot A'_n\t \1\, d\t \bnu_n\nonumber \\ 
&&\hspace{.5cm}=\beta\int_{D^n}\t \psi^2\vee E_\cdot e^{-\beta\tau_{ 
A^c}}\, d\t \bnu_n\nonumber \\ 
&&\hspace{1.0cm}-\lim_{t\downarrow 0}\frac1t\int\left(\t T_{n,t}\left( 
\psi^2\vee{\Bbb E}_\cdot e^{-\beta\tau_{B^c}}\right)-\t\psi^2\vee 
E_\cdot e^{-\beta\tau_{A^c}}\right)\, d\t \bnu_n\nonumber \\ 
&&\hspace{.5cm}\le\beta\int_{D^n}\t\psi^2\vee E_\cdot e^{-\beta\tau_{ 
A^c}}\, d\t \bnu_n\nonumber \\ 
&&\hspace{1.0cm}-\int_{\hat{A}}{\T \frac12}\Delta\left(\t\psi^2\vee 
E_\cdot e^{-\beta\tau_{A^c}}\right)\, d\t \bnu_n-\int_{\hat{A}^c}{\T 
\frac12}\Delta\left(\t \psi^2\vee E_\cdot e^{-\beta\tau_{A^c}}\right) 
\, d\t \bnu_n\nonumber \\ 
&&\hspace{1.0cm}-\int\left(\t\psi^2\vee E_\cdot e^{-\beta\tau_{A^c}} 
\right)\, d\int_{x\in D^n}(m^n_x-\mu^n_x)\left({\T -\frac12\Delta}\t 
m_n\right)(x)\, dx
\end{eqnarray} 
where, for the last $``\, \le\, "$ sign, we have taken into consideration 
that ${\t \psi^2\vee E_\cdot e^{-\beta\tau_{A^c}}}$ is convex from below 
in direction of every vector orthogonal to $\partial\{\t \psi^2\ge E_\cdot 
e^{-\beta\tau_{A^c}}\}$. Recalling that $E_\cdot e^{-\beta\tau_{A^c}}$ is 
$\beta$-harmonic on $\hat{A}^c$ with respect to ${\T\frac12}\Delta$ we 
obtain 
\begin{eqnarray}\label{4.69} 
&&\hspace{-.5cm}(\beta-2C){\Bbb E}_{\sbnu_n}e^{-\beta\tau_{B^c}}\le\beta 
\int_{\hat{A}}\t\psi^2\, d\t \bnu_n-\int_{\hat{A}}{\T \frac12}\Delta\t 
\psi^2\, d\t\bnu_n-\int_{x\in {D^n}}\int_{y\in\partial{D^n}}\sum_{i=1}^n 
\chi_{(y)_i\in\partial D}(y)\times\nonumber \\ 
&&\hspace{.5cm}\times\int_{z_i\in D}\left(\t \psi^2\vee E_\cdot e^{-\beta 
\tau_{A^c}}(z)-\t \psi^2\vee E_\cdot e^{-\beta\tau_{A^c}}(y)\right)\eta 
(z_i)\, dz_i\, \mu^n_x(dy)\, \t {\bf m}_n(dx)\, . \qquad
\end{eqnarray} 
From the jump mechanism of the process $X$ we observe that, for any $y\in 
\partial D^n\setminus\partial A$ with $y=(z_1,\ldots ,z_{i-1},y_i,z_{i+1}, 
\ldots ,z_n)\in\partial^{(1)}D^n$ where $z_1,\ldots ,z_{i-1},z_{i+1}, 
\ldots ,z_n\in D$ and $y_i\in\partial D$, $i\in\{1,\ldots ,n\}$, we have 
\begin{eqnarray*}
E_y e^{-\beta\tau_{A^c}}=\int_{z_i\in D}E_ze^{-\beta\tau_{A^c}}\eta(z_i) 
\, dz_i\, .  
\end{eqnarray*} 
Here $z_1,\ldots ,z_{i-1},z_{i+1}, \ldots ,z_n$ have been fixed in the 
notation $z=(z_1,\ldots ,z_n)$. With the same notation it follows that 
\begin{eqnarray}\label{4.70} 
\int_{z_i\in D}\t\psi(z)^2\vee E_ze^{-\beta\tau_{A^c}}\eta (z_i)\, dz_i- 
\t\psi(y)^2\vee E_ye^{-\beta\tau_{A^c}}\ge 0\quad\mbox{\rm on}\quad y\in 
\hat{A}^c\cap\partial D^n\, .   
\end{eqnarray} 
Putting (\ref{4.69}) and (\ref{4.70}) together we find 
\begin{eqnarray*} 
&&\hspace{-.5cm}(\beta-2C){\Bbb E}_{\sbnu_n}e^{-\beta\tau_{B^c}}\le\int_{ 
\hat{A}}\left(\beta\t \psi^2-2\left({\T\frac12}\Delta\t\psi\right)\t \psi 
\right)\, d\t\bnu_n-\int_{D^n}\int_{y\in\hat{A}\cap\partial{D^n}}\sum_{i= 
1}^n\chi_{(y)_i\in\partial D}(y)\times \\ 
&&\hspace{1.0cm}\times\int_{z_i\in D}\left(\t\psi^2(z)\vee E_ze^{-\beta 
\tau_{A^c}}-\t\psi^2(y)\right)\eta(z_i)\, dz_i\, \mu^n_\cdot(dy)\, \t{\bf 
m}_n(dx) \\ 
&&\hspace{.5cm}\le\int\left|\beta\t\psi-2\left({\T\frac12}\Delta\t\psi 
\right)\right|\, |\t\psi|\, d\t\bnu_n-\int_{x\in D^n}\int_{y\in 
\hat{A}\cap\partial{D^n}}\sum_{i=1}^n\chi_{(y)_i\in\partial D}(y)\times 
 \\ 
&&\hspace{1.0cm}\times\int_{z_i\in D}\left(\t\psi^2(z)-\t\psi^2(y)\right) 
\eta(z_i)\, dz_i\, \mu^n_x(dy)\, \t {\bf m}_n(dx) \\ 
&&\hspace{.5cm}\le\int\left|\beta\t\psi-2\left({\T\frac12}\Delta\t\psi 
\right)\right|\, |\t\psi|\, d\t\bnu_n+\int_{x\in D^n}\int_{y\in\partial 
{D^n}}\sum_{i=1}^n\chi_{(y)_i\in\partial D}(y)\times \\ 
&&\hspace{1.0cm}\times\left|\int_{z_i\in D}\left(\t\psi^2(z)-\t\psi^2(y) 
\right)\eta(z_i)\, dz_i\right|\, \mu^n_x(dy)\, \t {\bf m}_n(dx) \\ 
&&\hspace{.5cm}=\frac{1}{\ve_n^2\left\|f-\frac1\beta Bf\right\|^2} 
\int\left|\beta (f-2g)-2\left({\T\frac12}\Del f\right)+\beta^2 G_{n, 
\beta}g\vphantom{\dot{f}}\right|\cdot\left|f-\beta G_{n,\beta}g\right|\, 
d\bnu_n \\ 
&&\hspace{1.0cm}+\frac{(\beta-2C)\cdot b(n)}{\ve_n^2}\vphantom{\frac{\D 
\int}{\int}} \\ 
&&\hspace{.5cm}\le\frac{\beta\|f\|+2\nn{\T\frac12}\Del f\nn+3\beta\|g\|} 
{\ve_n^2\left\|f-\frac1\beta Bf\right\|^2}\left(\int\left(f-\beta G_{n, 
\beta}g\right)^2\, d\bnu_n\right)^{1/2}+\frac{(\beta-2C)\cdot b(n)} 
{\ve_n^2} \\ 
&&\hspace{.5cm}\le\frac{\D C_1\left(\int\left(f-\beta G_{n,\beta}g\right 
)^2\, d\bnu_n\right)^{1/2}}{\ve_n^2\left\|f-\frac1\beta Bf\right\|^2 }+ 
\frac{(\beta-2C)\cdot b(n)}{\ve_n^2}\, .  
\end{eqnarray*} 
In other words, for $n\ge 2$, 
\begin{eqnarray*} 
{\Bbb E}_{\sbnu_n}e^{-\beta\tau_{B^c}}\le\frac{\D C_2\left(\int 
\left(f-\beta G_{n,\beta}g\right)^2\, d\bnu_n\right)^{1/2}+b(n)} 
{\ve_n^2}=\ve_n  
\end{eqnarray*} 
which gives (c8) by the already verified relation $\ve_n\stack{n\to 
\infty}{\lra}0$. 
\qed 
\medskip 

\nid 
{\bf Remark }(1) We recall Remark (2) of Section 2. We recall also that 
condition (c3) has been verified in Step 2 of the proof of Proposition 
\ref{Proposition4.11} as a consequence of (c3'). Furthermore, we have 
${\cal T}=\{T_t g:g\in {\cal C},\ \beta >0\}\subseteq {\cal C}$ by Lemma 
\ref{Lemma4.9} (d). For ${\cal T}_n=\{T_{n,t}g :g\in {\cal C},\ \beta >0\} 
\subseteq {\cal C}$, $n\in {\Bbb N}$, we recall the proof of Lemma 
\ref{Lemma4.4}, in particular (\ref{4.34}). Indeed, $g\in {\cal C}$ implies 
by (c1') that $\t g$ bounded and continuous on $D^n$ for which we have 
\begin{eqnarray*}
&&\hspace{-.5cm}\t T_{n,t}\t g(x)=\int\left(\t g-h_{\t g}\right)(y)\, 
P_x(B^{D^n}_t\in dy)+h_{\t g}(x) \\ 
&&\hspace{.5cm}+\sum_{i=1}^n\int_i\int\int_{z_i\in D}\int_{y_i\in\partial D} 
\left(\t g(z)-\t g(y)\right)\frac{P_x(B_\tau\in dy,\, \tau <t)}{\sigma(dy)}\, 
s(dy_i)\, \eta(z_i)\, dz\, .
\end{eqnarray*} 
For $g\in {\cal C}$, ${\cal C}\ni g_n\sstack{n\to\infty}{\lra}g$, and $\beta> 
{\T\frac12}\|A'\1\|_{L^\infty(E,\sbnu)}\vee\, \sup_{n\in {\Bbb N}}{\T\frac12} 
\|A_n'\1\|_{L^\infty (E,\sbnu_n)}$ we have by Theorem \ref{Theorem2.14} and 
Remark (2) of Section 2 $\ T_{n,t}g_n\sstack{n\to\infty}{\lra} T_tg$.

\end{document}